\documentclass[final,pagebackref,coverfontpercent=100]{adsphd}

\usepackage[english,dutch]{babel}

\usepackage{nomencl}   
\usepackage[style=long,number=none]{glossary} 
\usepackage{imakeidx}
\makeindex[intoc]

\title{Extensions of submanifold theory to non-real settings, with applications}

\author{Victor}{Pessers}

\supervisor{Prof.~dr.~J.~Van~der~Veken}{}
\president{Prof.~dr.~J.~Deprez}
\jury{
Prof.~dr.~F.~Cools\\
Prof.~dr.~W.~Goemans\\
Prof.~dr.~L.~Vrancken\\
Prof.~dr.~M.~Zambon
}
\externaljurymember{Prof.~dr.~A.~Carriazo}{Universidad de Sevilla}

\phddegree{Science: Mathematics} 
\faculty{Faculty of Science}
\department{Department of Mathematics}

\researchgroup{Geometry Section}
\website{}
\email{victor.pessers@wis.kuleuven.be} 

\address{Celestijnenlaan 200B box 2400}
\addresscity{B-3001 Leuven} 

\date{September 2016}
\copyyear{2016}

\renewcommand{\nomname}{List of Symbols}
\makeatletter
\let\@printnomenclatureorig\@printnomenclature
\def\@printnomenclature[#1]{%
  \cleardoublepage%
  \chaptermark{\nomname}
  \@printnomenclatureorig[#1]
}
\makeatother
\makenomenclature

\newcommand{\glossname}{Abbreviations}
\makeglossary

\let\printglossaryorig\printglossary
\renewcommand{\printglossary}{
  \renewcommand{\glossaryname}{\glossname}
  \cleardoublepage
  \printglossaryorig\chaptermark{\glossname}}

\bibliography{allpapers}

\bibstyle{abbrv}

\usepackage{amssymb,amsthm}
\usepackage{amsmath}
\usepackage{mathtools}
\usepackage{twoopt}
\usepackage{thmtools}
\usepackage{theoremref}
\usepackage[mathscr]{euscript}
\usepackage{babel}
\usepackage{tikz}
\usetikzlibrary{babel}
\usepackage{tikz-cd}
\usepackage{graphicx}
\usepackage{hyperref}
\usepackage{amsthm}
\usepackage{enumerate}
\usepackage{todonotes}
\presetkeys{todonotes}{fancyline}{}

\theoremstyle{plain}
\newtheorem{theorem}{Theorem}[section]
\newtheorem{lemma}[theorem]{Lemma}
\newtheorem{proposition}[theorem]{Proposition}
\newtheorem{corollary}[theorem]{Corollary}

\theoremstyle{definition}
\newtheorem{definition}[theorem]{Definition}
\newtheorem{example}[theorem]{Example}
\newtheorem{remark}[theorem]{Remark}

\declaretheoremstyle[
headfont=
\bfseries, 
spaceabove = 0.2cm, 
spacebelow = 0.2cm
]{mydef}
\theoremstyle{mydef}
\newtheorem{p}[theorem]{}

\newcommand{\titleone}{Submanifold theory with Rinehart spaces}
\newcommand{\titletwo}{Submanifold theory and Wick-relations}
\newcommand{\emindex}[1]{{\em #1}\index{#1}}

\begin{document}


\makefrontcoverXII

\maketitle

\frontmatter 

\includepreface{preface}
\includeabstract{abstract}
\includeabstractnl{abstractnl}



\tableofcontents



\mainmatter 

\instructionschapters\cleardoublepage

  \graphicspath{{\chaptersdir/introduction/image/}}%
  \chapter{Introduction}\label{ch:introduction}

{

\def\lrs{Lie-Rinehart algebra}
\def\lrm{Rinehart space}

This thesis consists of two different contributions to Riemannian geometry and in particular submanifold theory, which both have in common that they relate Riemannian geometries over different ground rings. The first contribution, which will be the topic of Chapter \ref{ch:chapter1}, aims to provide a general framework for doing Riemannian geometry over an arbitrary ring, and it does so in terms of {\lrm s}, as we refer to them, rather than classical manifolds. Among others, this framework unifies variations of Riemannian geometry over ground rings like the real, complex and split-complex numbers. The second contribution, which is the subject of Chapter \ref{ch:chapter2}, also has a unifying character, but in this case we show how Riemannian geometry over the complex numbers (more particularly holomorphic Riemannian geometry) can be used to relate the real pseudo-Riemannian geometries of so-called Wick-related spaces.

In the following of this chapter, we provide introductions for the topics mentioned above. Part of this chapter originates from the introduction of \cite{PeVa}, which is joint work with Joeri Van der Veken.

\section{\titleone}

The concept of {\lrs s} makes it possible to describe the main objects in differential geometry, such as manifolds and certain derived structures, in a formal algebraic language. The notion of {\lrs s} already exists for quite a while in the literature (its earliest implicit presence can be traced back to as early as 1944 \cite{Jac2,Hu04}) but it has not been until relatively recently  that the topic is receiving some serious research interest. Now {\lrs s} have already found many applications in various branches of differential geometry. A good overview of the history and some of the applications of {\lrs s}, as well as many more references, can be found in \cite{Hu04}. 

There are a couple of reasons why one might consider the use of {\lrs s} beneficial. First, because the notion of {\lrs s} is very general in nature, there are many theorems in differential geometry which are stated in different settings, but which can be unified in a single statement on {\lrs s}. Also, in some situations where singularities are involved, {\lrs s} may provide an appropriate framework to reason in a clear and valid manner about them, whereas this is not always possible or practical in the classical language of manifolds with charts. As we will demonstrate in this thesis, it even opens up the possibility to study Riemannian geometry in settings that would not make clear sense otherwise, such as Riemannian geometry over finite fields (cf. Example \ref{ex:finite field}). A disadvantage of the involvement of {\lrs s}, is perhaps that the proofs may geometrically be less intuitive to follow. Note that in all of the advantages and disadvantages listed above, the potential benefits of {\lrs s} for differential geometry might be somewhat similar to the benefits which the ring-theoretic approach has brought to algebraic geometry.

For our particular goals we aim to fulfil, we will restrict our attention to the notion of {\lrm s}, which can be regarded as a special kind of {\lrs s}, endowed with some extra structure. In a way, {\lrm s} provide a generalization (in the form of an algebraic analogue) of the concept of differentiable manifolds, whereas {\lrs s} also entail derived structures such as Lie-algebroids. We would like to emphasize, that our aim is to use {\lrs s} (in the form of {\lrm s}) as a tool to enrich the study of Riemannian geometry, rather than to enrich the study of {\lrs s} with the concept of Riemannian metrics, as has been done in \cite{Bou,Bru} for the particular case of real Lie-algebroids. We will touch once more on how {\lrm s} compare to {\lrs s} in the conclusion of this thesis.



\section{\titletwo}

In Chapter \ref{ch:chapter2} of this thesis, we show how certain problems in \mbox{(pseudo-)}\linebreak[0]Riemannian submanifold theory, that are situated in different ambient spaces, can be related to each other by translating the problem to an encompassing holomorphic Riemannian space. This approach seems new in the area of submanifold theory, although it incorporates several existing insights, such as the theory on analytic continuation, complex Riemannian geometry and real slices, as well as the method of Wick rotations, which is often used in physics.

The relation between pseudo-Riemannian geometry and complex analysis can be traced back to the very birth of pseudo-Riemannian geometry. In the early publications on Lorentzian geometry by Poincar\'e and Minkowski (cf. \cite{Po,Mi}), the fourth coordinate of space-time was represented as $it$ (or $ict$, $c$ being the speed of light), so that space-time was essentially modelled as $\mathbb{R}^3\times i\mathbb{R}$, where the standard complex bilinear form played the role of metric. Likely due to a later reformulation by Minkowksi himself, this point of view soon fell in abeyance in favor of the nowadays more common presentation in terms of an indefinite real inner product. Admittedly, as long as ones attention is kept restricted to four-dimensional Minkowski space alone, the use of complex numbers to deal with the signature bears little advantage.

This relationship between space-time geometry and complex numbers received renewed attention, when it was shown by Wick how problems from the Lorentzian setting are turned into problems in a Euclidean setting, after a so-called Wick rotation is applied on the time coordinate (cf. \cite{Wi}). This method of Wick rotations lays at the basis of the later theory on Euclidean quantum gravity, developed by Hawking and Gibbons among others (cf. \cite{Ha}), but also in other areas of theoretical physics it remained a valuable tool ever since.

Despite that the concept of Wick rotations is known by physicists for quite some time already, there are still many research domains where such insights have not been fully exploited yet. Our aim will be to demonstrate that also for the particular area of submanifold theory, one may benefit from taking a complex viewpoint to problems about submanifolds of pseudo-Riemannian spaces. We will show how holomorphic Riemannian geometry can be used to relate certain kinds of submanifolds in one pseudo-Riemannian space to submanifolds with corresponding geometric properties in other so-called Wick-related spaces. It should be noted that the subject of complex Riemannian manifolds, first introduced in an also physically motivated article by LeBrun (cf. \cite{Le}), has appeared in several articles on submanifold theory (e.g. \cite{BoFrVo,Sl}), though not as a tool but rather as an object of study in its own right.

\section{Further remarks}

It is important to remark that the content of Chapter \ref{ch:chapter2} is not phrased in the same theoretical language as is developed in Chapter \ref{ch:chapter1}, but in the classical language of manifolds instead. The reason for that is not that it isn't possible, but simply because the theory in Chapter \ref{ch:chapter1} is to be developed much further in order to do so. For example, in Chapter \ref{ch:chapter2} we use local arguments now and then, but in order to do so in the framework of {\lrm s}, one would have to extend this theory to sheafs of {\lrm s}. Another relevant aspect we have not investigated in Chapter \ref{ch:chapter1}, is the spectrum (and related notions) of the function algebras under consideration, which would be necessary to deal with point-wise phenomena. But as we will see, there is a lot that can be said and done even before the involvement of such tools. We will come back to these and other issues in Chapter \ref{ch:conclusion}, where we will discuss the obtained results, as well as possible directions for further research.

}



%
%
%
%
%
%
%
\cleardoublepage

%

%
  \graphicspath{{\chaptersdir/chapter1/image/}}%



\chapter{\titleone}\label{ch:chapter1}

{

\newcommand{\Hom}{\operatorname{Hom}}
\newcommand{\Der}{\operatorname{Der}}
\newcommand{\Ker}{\operatorname{Ker}}
\newcommand{\Ima}{\operatorname{Im}}
\newcommand{\D}{\operatorname{D}}
\newcommand{\ad}{\operatorname{ad}}
\newcommand{\spa}{\operatorname{span}}
\renewcommand{\d}{\mathrm{d}}
\newcommand{\di}[1]{\d^{#1}}

\renewcommand{\-}[1]{}
\newcommand{\R}{\mathbb{R}}
\newcommand{\C}{\mathbb{C}}
\newcommand{\F}{\mathbb{F}}
\newcommand{\Lie}{\mathcal{L}}
\newcommand{\defeq}{\coloneqq}
\renewcommand{\vec}[1]{#1^\sharp}
\newcommand{\form}[1]{#1^\flat}
\newcommand{\al}[1]{\mathcal{#1}}
\newcommand{\alg}{\al O}
\newcommand{\ideal}{\al I}
\newcommand{\ring}[1][\mathbb{K}]{#1}
\renewcommand{\mod}[1]{{\bf #1}}
\renewcommand{\sup}[1]{\hat{#1}}
\newcommand{\sub}[1]{\check{#1}}
\newcommand{\dual}[1]{{#1}^\vee}
\newcommand{\cycsum}[1]{\sum^{\mathrm{cycl}}}
\newcommand{\derx}{\mathfrak{X}}
\newcommand{\derf}{\Omega}
\newcommand{\lra}[1][A]{\mathscr{#1}}
\newcommand{\tander}{\derx^\top}
\newcommand{\derI}{\derx_{\ideal}}
\newcommand{\rhos}{\rho^{\sharp}}
\newcommand{\rhof}{\rho^{\flat}}
\newcommand{\fun}[1]{#1^*}
\newcommand{\rhon}{\fun \rho}
\newcommandtwoopt{\met}[2][\cdot][\cdot]{\langle #1,#2\rangle}
\newcommandtwoopt{\metA}[2][\cdot][\cdot]{\langle #1,#2\rangle_{\!\lra}}
\newcommandtwoopt{\metI}[2][\cdot][\cdot]{\langle #1,#2\rangle_{\!\lra/\ideal}}
\newcommandtwoopt{\pair}[2][\cdot][\cdot]{\met[#1][#2]}
\newcommandtwoopt{\pairpif}[2][\cdot][\cdot]{\pair[#1][#2]}
\newcommand{\derz}{\derx_{\ideal}}
\newcommand{\del}[1]{\frac{\partial}{\partial #1}}
\newcommand{\grad}[1]{\d^\sharp #1}
\newcommand{\cor}{x}
\newcommand{\corv}{X}
\renewcommand{\d}{d}
\newcommand{\con}{\al C}
\newcommand{\ort}[1]{#1^{\top}}
\newcommand{\rhoc}[1]{\rho \circ #1}
\renewcommand{\pmod}{\mod M}
\newcommand{\dmod}{\mod W}
\newcommand{\op}[1]{#1^{\mathrm{op}}}
\newcommand{\derive}[1]{\d_{#1}}
\newcommand{\dird}[2]{\delta_{#1#2}}

\newcommand{\lr}{Lie-Rinehart}
\newcommand{\lrs}{{\lr} algebra}
\renewcommand{\r}{Rinehart}
\newcommand{\lrm}{{\r} space}
\newcommand{\man}{???}
\newcommand{\lrsm}{{\r} subspace}
\newcommand{\gen}{duplicative}
\newcommand{\xspace}{tangent space}
\newcommand{\omspace}{cotangent space}
\newcommand{\inp}{inner product}
\newcommand{\musical}{musical}
\newcommand{\dm}{{\musical} \inp}
\newcommand{\rmm}{Riemannian metric}

In this chapter, partly based on joint work with Joeri van der Veken \cite{PeVa2}, we describe how the notion of what we call {\lrm s}, can be used to provide a generalized theory of Riemannian (sub)manifolds. In the first three sections we rehearse and further develop some theoretical preliminaries on respectively modules and tensors, dual pairs and inner product modules, and lastly derivations on modules. These sections are rather self-contained, and basic knowledge on common algebraic structures should be sufficient to understand the proofs, but background information can for instance be found in \cite{MiHu}, \cite{Jar} and \cite{Jac1}. However, our treatment extends the usual theory to better fit our goals, in particular to deal with settings where the underlying algebras have a rich collection of ideals. We are not aware if such a treatment can also be found in other sources. Also, it should be noted that the terminology used in sources like mentioned above (such as the meaning of an {\em inner product} in \cite{MiHu}) will not always coincide with our terminology.

Having dealt with the algebraic prerequisites, we will come to the main contribution of this chapter, by showing in Section \ref{sec:Rinehart spaces} how the notion of {\lrm s} can provide a generalized framework for doing Riemannian geometry. Throughout this chapter, we will always assume rings and algebras (except Lie algebras) to be associative and commutative.

\section{Modules and tensors}\label{sec:modules_and_tensors}

\def\A{\lra} 
\def\B{\lra[B]}
\def\c{\lra[C]}
\def\modi{\pmod_{\A}} \def\modii{\pmod_{\B}}
\def\algi{\alg_{\A}} \def\algii{\alg_{\B}} 
\def\modhom{\phi} \def\alghom{\rho}
\def\W{W}
\def\V{V}
\def\U{U}
\def\X{M}
\def\Y{L}
\def\Z{K}

\begin{p}
An algebra strictly speaking consist of the combination of two seperate data, namely a ground ring $\ring$ and the actual algebra $\alg$ over this ground ring. For example, one may whish to distinguish between an algebra of complex valued functions (on a certain domain) regarded as algebra over $\C$ or as algebra over $\R$. However, because throughout the following all algebras we consider have (per situation) a fixed common ground ring, we do not need to be very strict about this distinction. Nevertheless, when it comes to modules, we will encouter situations where two modules are the same, except for the fact that their underlying algebras are different. For this reason, by a module we will usually mean an object $\A$ consisting of an algebra $\alg = \alg_{\A}$\index{$\mathcal{O}_{\mathscr{A}}$} (over some ground ring $\ring$), and the actual module $\bf{M} = \modi$\index{$\bf{M}_{\mathscr{A}}$} over the algebra $\alg$. So throughout the following, we will either refer to a module as such a combined object $\A$, or we more explicitly introduce it as a module $\mod{M}$ over an algebra $\alg$.
\end{p}

\begin{p}
Suppose $\A$ and $\B$ are modules in the above sense. As we mentioned just before, the ground ring $\ring$ of the underlying algebras $\alg_{\A}$ and $\alg_{\B}$ is assumed to be fixed. Let us furthermore assume that we have an algebra homomorphism $\rhon: \algi \to \algii$, i.e. a $\ring$-linear map, preserving the multiplicative structure on $\algi$. In such a situation, the module $\B$ can be turned into a module over the algebra $\algi$ as well, with scalar multiplication being defined by
\[
f \X \defeq \rhon(f) \X
\]
for $f\in \algi$ and $\X\in \modii$.

By definition, a morphism $\rho:\A\to \B$ consists of an algebra homomorphism $\rhon:\algi\to \algii$\index{$\rho^*$}, and a $\rhon$-equivariant module homomorphism $\rhos:\modi\to \modii$\index{$\rho^{\sharp}$}, which means that $\rhos$ is a homomorphism of Abelian groups, which is moreover compatible with respect to scalar multiplication through the map $\rhon$, i.e.
$$\rhos(f \X) = \rhon(f) \rhos(\X)$$
for all $f\in \algi$ and $\X\in \modi$. We call a morphism $\rho:\A\to \B$ surjective (resp. injective), if $\rhos$ and $\rhon$ are both surjective (resp. injective). It is not a coincidence, by the way, that the notation $\rhon$ is reminiscent of the common notation for a pull-back map, as in concrete examples we will discuss later, the map $\rhon$ will actually turn out to be a pull-back map between the function algebras of two manifolds.
\end{p}

\def\modi{\pmod_{1}} \def\modii{\pmod_{2}}
\def\modiii{\pmod_{\lra[A]}}
\def\tmap{F} 

\begin{definition}
The {\em dimension of a module}\index{dimension of a module} $\pmod$ with respect to some underlying algebra $\alg$, 
is the largest cardinal number $\kappa$ such that there exists a subset of $\pmod$ of cardinality $\kappa$, whose elements are all linearly independent over $\alg$. In particular, the dimension of the zero module $\{0\}$ is $0$. 
\end{definition} 

\begin{p}
Given two modules $\modi$ and $\modii$ over the same algebra $\alg$, we may construct their {\em tensor product}\index{tensor product}, which we denote by $\modi\otimes_{\alg} \modii$, or even $\modi\otimes \modii$ if the underlying algebra can be understood from context. This space is uniquely characterized by the universal property that there is an $\alg$-bilinear map $\otimes:\modi\times \modii \to \modi\otimes \modii$, such that any bilinear map $\tmap$ from $\modi\times \modii$ to some $\alg$-module $\modiii$, factors as 
$$\tmap: \modi\times \modii\xrightarrow{\otimes} \modi\otimes\-{_{\alg}} \modii \xrightarrow{\tmap'} \modiii,$$ where $\tmap'$ is $\alg$-linear as well. 

The tensor product of two $\alg$-modules is again an $\alg$-module, and taking tensor products is an associative operation up to canonical isomorphism of the spaces as $\alg$-modules. The above property thus extends to multilinear maps with an arbitrary number of arguments: any $\alg$-multilinear map $\tmap:\pmod_{1}\times\ldots\times\pmod_{n} \to \modiii$ factors as 
$$\tmap: \pmod_{1}\times\ldots\times\pmod_{n}\xrightarrow{\otimes} \pmod_{1}\otimes \ldots\otimes\pmod_{n} \xrightarrow{\tmap'} \modiii.$$
Because the map $F'$ also uniquely determines the map $F$, we will usually identify both maps and simply write $F:\pmod_{1}\otimes \ldots\otimes\pmod_{n} \to \modiii$. In case $\modiii = \alg$, we refer to such a map as a {\em multilinear form}\index{form!multilinear}.

\end{p}

%

\begin{definition}\label{def:non-degenerate}
Suppose $\pmod_1,\ldots,\pmod_n,\modiii$ are all modules over a common algebra $\alg$, and that we have an $\alg$-multilinear map $\tmap:\pmod_1\otimes\ldots\otimes \pmod_n\to \modiii$. If $\X\in \pmod_k$, we define the $k$-th argument {\em contraction}\index{contraction of argument} of $\tmap$ with $\X$ as the $\alg$-linear map $\X\lrcorner_k \tmap:\pmod_1\otimes\ldots\otimes \pmod_{k-1}\otimes\pmod_{k+1}\otimes\ldots\otimes \pmod_n\to \modiii$\index{$\lrcorner_k$}, given by
$$\X\lrcorner_k \tmap(\alpha\otimes \beta) \defeq F(\alpha\otimes \X\otimes \beta),$$
for any $\alpha\in \pmod_1\otimes\ldots\otimes \pmod_{k-1}$ and $\beta\in \pmod_{k+1}\otimes\ldots\otimes \pmod_{n}$. When writing $X\lrcorner \tmap$ (so without specifying the number $k$), we will always mean contraction with the first argument. It should be noted that under canonical identification of tensor products, the notion of `$k$-th argument' depends implicitly on the particular decomposition of the domain of $\tmap$. For example, if $\pmod_1 = \pmod_1'\otimes \pmod_1''$, and we would write the domain as $\pmod_1'\otimes \pmod_1''\otimes\ldots\otimes \pmod_n$, then the meaning of the $k$-th argument contraction changes accordingly. 

If $\pmod_1\otimes\ldots\otimes \pmod_n$ is not a zero module, then a map $\tmap:\pmod_1\otimes\ldots\otimes \pmod_n\to \modiii$ is said to be {\em non-degenerate}\index{non-degenerate!map} in the $k$-th argument (where $k\in\{1,\ldots,n\}$), if for any non-zero $\X\in\pmod_k$ we have that $\X\lrcorner_k \tmap$ is not the zero map. If the map is non-degenerate in all of its arguments, we simply call it non-degenerate. Also in case $\pmod_1\otimes\ldots\otimes \pmod_n$ is {\em trivial}\index{trivial module}, i.e. the zero module, we call $\tmap$ non-degenerate by convention. Note that the notion of non-degeneracy also depends on the chosen decomposition of the domain of $\tmap$.
\end{definition}

\section{Dual pairs}\label{sec:dual_pairs}

\subsection{Dual pairs and inner products}

\def\alghom{\rho} \def\modhom{\phi} 
\def\forhom{\modhom^{\flat}} 

\def\W{W}
\def\X{M}
\def\Y{L}
\def\pairing{pairing}
\newcommandtwoopt{\epair}[2][\cdot][\cdot]{\met[#1][#2]}
\newcommandtwoopt{\emet}[2][\cdot][\cdot]{\met[#1][#2]}
\def\rhon{\fun \rho}

\begin{definition}
A {\em dual pair}\index{dual pair} $\lra$ consists of two modules, a {\em primary module}\index{primary module} $\pmod = \pmod_{\A}$ and a {\em secondary module}\index{secondary module} $\dmod = \dmod_{\A}$\index{$\bf{W}_{\mathscr{A}}$}, both over a common algebra $\alg = \alg_{\A}$ with ground ring $\ring$, for which moreover there exists
an $\alg$-bilinear form 
$$\epair_{\A}:\pmod\otimes \dmod\to \alg,\ \ \X\otimes \W \mapsto \epair[\X][\W],\index{$\langle \cdot,\cdot\rangle_{\mathscr{A}}$}$$
which we refer to as the {\em pairing}\index{pairing} of the dual pair. We call a dual pair {\em non-degenerate}\index{non-degenerate!dual pair} if its {\pairing} is non-degenerate, and {\em non-trivial}\index{non-trivial dual pair} if none of $\alg$, $\pmod$ and $\dmod$ are trivial. Moreover, we will not be strict about the order of arguments for pairings between $\pmod$ and $\dmod$, i.e. for $\X\in \pmod$ and $\W\in \dmod$, we do not distinguish between $\epair[\X][\W]$ and $\epair[\W][\X]$.
\end{definition}

\begin{remark}\label{rem:determined_by_non-deg}
Note that for a non-degenerate dual pair, any element $\X\in \pmod$ is completely determined
by the values of $\epair[\X][\W]$ for varying $\W$: if $\epair[\X][\W] = \epair[\Y][\W]$ for all $\W\in \dmod$, then it follows that $\X = \Y$.
\end{remark}

\begin{definition}\label{def:morphism_dual_pairs}
By definition, a {\em morphism of dual pairs}\index{morphism!of dual pairs} $\rho:\A\to\B$, entails a morphism between both their primary modules and their secondary modules, which preserves the \pairing. More specifically, it consists of the following:
\begin{itemize}
\item an algebra homomorphism $\rhon:\alg_{\A}\to \alg_{\B}$ between the underlying algebras of $\A$ and $\B$,
\item a $\rhon$-equivariant module homomorphism $\rhos:\pmod_{\A}\to\pmod_{\B}$,
\item a $\rhon$-equivariant module homomorphism $\rhof:\dmod_{\A}\to\dmod_{\B}$,\index{$\rho^{\flat}$}
\end{itemize}
and the additional requirement that
\begin{equation}\label{eq:dual pair morphism}
\rhon \pair[\X][\W]_{\A} = \pair[\rhos \X][\rhof \W]_{\B}
\end{equation}
for all $\X\in\pmod_{\A}$ and $\W\in\dmod_{\A}$. We call a morphism of dual pairs $\rho$ surjective (resp. injective), if $\rhon$, $\rhos$ and $\rhof$ are all surjective (resp. injective).

\end{definition}

\def\dA{\op{\A}}
\begin{example}\thlabel{ex:opposite dual pairs}
For any dual pair $\A$, there exists a unique {\em opposite dual pair}\index{opposite dual pair} $\dA$\index{$\mathscr{A}^{\mathrm{op}}$} which is given by  $\alg_{\dA} = \alg_{\A}$, $\pmod_{\dA} = \dmod_{\A}$ and $\dmod_{\dA} = \pmod_{\A}$, and such that the pairing between elements $\W$ and $\X$ of $\pmod_{\dA}$ and $\dmod_{\dA}$ respectively, is simply taken to be the pairing between $\W$ and $\X$ as elements of $\dmod_{\A}$ and $\pmod_{\A}$ respectively. 
\end{example}

\begin{example}
Any module $\pmod$ over an algebra $\alg$, forms a dual pair together with its {\em dual space}\index{dual space} $\dual \pmod \defeq \Hom_{\alg}(\pmod,\alg)$\index{$\bf{M}^\vee$}, i.e. the space of homomorphisms from $\pmod$ to $\alg$. Because a module can be a module over several different algebras, it should be noted that the concept of dual space depends on a specification of which underlying algebra is meant. Note that for a non-degenerate dual pair $\A$, there always exists a canonical injective morphism $\A\to (\pmod_{\A},\dual \pmod_{\A})$ (where the latter dual pair also has $\alg_{\A}$ as underlying algebra).

\end{example}

\subsubsection{Inner product modules}

\def\subi{\mod M}
\def\subii{\mod N}
\begin{definition}\label{def:bilinear form inner product etc}
An important special case of a dual pair, is a dual pair of the form $(\pmod_,\pmod)$, i.e. a dual pair for which the primary and secondary module are identical. In this case, the pairing can be regarded as a bilinear form on $\pmod$.
If $\met[\X][\Y] = \met[\Y][\X]$ for all $\X,\Y\in \pmod$, we say that $\met$ is {\em symmetric}\index{symmetric!bilinear form}. 
A symmetric bilinear form we will call an {\em \inp}\index{inner product}, and a module endowed with an inner product we call an {\em {\inp} module}\index{inner product module}. 
Note that $\met$ is an {\inp}, if and only if the identity morphism $\rho$ (i.e. $\rhon$, $\rhos$ and $\rhof$ all being the identity) between the dual pair $(\mod M,\mod M)$ and its opposite dual pair (cf. \thref{ex:opposite dual pairs}) yields an isomorphism.
\end{definition}

\subsubsection{Inner product pairs}

\def\pmodi{\pmod_{\A}} \def\dmodi{\dmod_{\A}}

{\def\pmodi{\pmod} \def\dmodi{\dmod}
\begin{definition}
By an {\em {\inp} pair}\index{inner product pair}, we mean a dual pair endowed with an {\inp} on its primary module $\pmodi$, i.e. a symmetric bilinear form $\emet:\pmodi\otimes \pmodi\to \alg$.\index{$\langle \cdot,\cdot\rangle_{\mathscr{A}}$} This {\inp} thus exists alongside the similarly denoted pairing $\pair:\pmodi\otimes \dmodi\to \alg$. (It can be inferred from the whether if $\emet$ refers to an inner product between elements from $\pmodi$ or a pairing between elements from $\pmodi$ and $\dmodi$.) If $\Y\in \pmodi$ and $\W\in \dmodi$ are such that $\epair[\X][\W] = \emet[\X][\Y]$ for all $\X\in\pmodi$, then we call $\Y$ and $\W$ {\em duplicates}\index{duplicate} of each other. If for every element $\Y\in\pmodi$ there exists a duplicate $\form{\Y}\in \dmodi$\index{$(\cdot)^{\flat}$}, then we say that $\dmodi$ {\em duplicates}\index{duplicate!to duplicate} $\pmodi$. 
Likewise, if for every element $\W \in \dmodi$ there exists a duplicate $\vec{\W}\in\pmodi$, then we say that $\pmodi$ duplicates $\dmodi$. If the {\inp} is such that $\pmodi$ and $\dmodi$ duplicate each other, we say that the {\inp} is {\em {\gen}}\index{duplicative!inner product}. An {\inp} which is both non-degenerate and {\gen}, we call a {\em {\dm}}\index{musical!inner product}. By a {\em {\gen}\index{duplicative!pair} (resp. {\musical}) pair}\index{musical!pair}, we mean an {\inp} pair whose {\inp} is {\gen} (resp. {\musical}).
\end{definition}

\begin{proposition}
For a {\gen} pair $\A$, there exists a unique {\inp} on $\A$'s secondary module $\dmod$, which we will again denote by $\met$, which satisfies
\begin{equation}\label{eq:dual metric}
\met[\W][\V] = \met[\vec{\W}][\vec{\V}],
\end{equation} 
for all $\W,\V\in \dmod$ with duplicates $\vec{\W}, \vec{\V}\in \pmod$ respectively. In this way a {\gen} (resp. {\musical}) {\inp} on $\A$ determines a {\gen} (resp. {\musical}) {\inp} on the opposite dual pair $\dA$.
\end{proposition}
\begin{proof}
We will prove that taking \eqref{eq:dual metric} as definition, determines a well-defined {\inp} on $\dmod$. Let us suppose that $\vec{\W_1}$, $\vec{\W_2}$ and $\vec{\V_1}$, $\vec{\V_2}$ are different duplicates of $\W$ and $\V$ respectively. Then we find that
$$
\met[\vec{\W_1}][\vec{\V_1}] = \met[\vec{\W_1}][\V] = \met[\vec{\W_1}][\vec{\V_2}] = \met[\W][\vec{\V_2}] = \met[\vec{\W_2}][\vec{\V_2}].
$$
When it comes to the final claim, it is straightforward to verify that this inner product is {\gen} if regarded as an {\inp} on the primary module of $\dA$, and moreover that it is non-degenerate if the inner product on $\A$ is non-degenerate.
\end{proof}

\begin{proposition}\label{prop:musical isomorphism}
On a {\musical} pair $\A$ there exists an isomorphism  ${\form{(\cdot)}:\pmod\to \dmod}$, such that 
\begin{equation}\label{eq:relation_inp_and_musical_isom}
\emet[\X][\Y] = \epair[\X][\form{\Y}] = \epair[\Y][\form{\X}],
\end{equation}
for all $\X, \Y\in \pmod$, which moreover implies that $\A$'s pairing is non-degenerate.

On the other hand, if $\A$ is a non-degenerate dual pair and there exists an isomorphism $\form{(\cdot)}:\pmodi\to \dmodi$ such that $\epair[\X][\form{\Y}] = \epair[\Y][\form{\X}]$ for all $\X,\Y \in \pmodi$, then there exists a unique {\dm} on $\A$, such that \eqref{eq:relation_inp_and_musical_isom} holds.
\end{proposition}
\begin{proof}
Note that from the symmetry of the inner product on $\pmodi$, the first equality in \eqref{eq:relation_inp_and_musical_isom} implies the second one. Now we claim that for any $\Y\in\pmod$ there exists exactly one element $\form{\Y}\in \dmod$, such that the first equation of \eqref{eq:relation_inp_and_musical_isom} holds for all $\X\in\pmod$. We have existence because our {\inp} is {\gen}, and uniqueness because of non-degeneracy. Hence there exists a well-defined bijective map $\form{(\cdot)}:\pmod\to \dmod$, such that $\epair[\X][\form{\Y}] = \epair[\X][\Y]$ for all $\X,\Y\in\pmod$. Further we verify that:
\begin{align*}
\emet[\X][\form{(f\Y+\Z)}] &= \epair[\X][f\Y+\Z] = f\epair[\X][\Y] + \epair[\X][\Z]\\ 
&= f\emet[\X][\form{\Y}]+\emet[\X][\form{\Z}] = \emet[\X][f\form{\Y}+\form{\Z}]
\end{align*}
for all $\X\in\pmod$, so that by the non-degeneracy of the inner product and Remark \ref{rem:determined_by_non-deg} it follows that $\form{(f\Y+\Z)} = f\form{\Y}+\form{\Z}$. So $\form{(\cdot)}$ is an isomorphism of $\alg$-modules. From this it also follows directly that the pairing of $\A$ is non-degenerate as well.

Conversely, define a bilinear form on $\pmod$ by 
$$
\emet[\X][\Y] \defeq \epair[\X][\form{\Y}].
$$ 
By our assumptions, it easily follows that this bilinear form is symmetric and thus an inner product. Its non-degeneracy is implied by the assumed non-degeneracy of the dual pair $\A$. Moreover, if we denote the inverse of the isomorphism $\form{(\cdot)}$ by $\vec{(\cdot)}:\dmod\to \pmod$, then by construction we have for any $\Y\in \pmod$ that $\form \Y\in\dmod$ satisfies $\epair[\X][\form \Y] = \emet[\X][\Y]$ for all $\X\in\pmod$, and similarly for any $\W\in \dmod$ we have that $\vec{\W}\in\pmod$ satisfies $\epair[\X][\W] = \emet[\X][\vec{\W}]$, which proves that the inner product is {\gen}. 
\end{proof}

\begin{definition}\label{def:inverse canonical isomorphism}
For a {\musical} pair, we refer to the map $\vec{(\cdot)}$ and its inverse $\form{(\cdot)}$ as the {\em musical isomorphisms}\index{musical!isomorphism}.\index{$(\cdot)^{\sharp}$}\index{$(\cdot)^{\flat}$}
\end{definition}

}

\subsection{Subpairs and quotient dual pairs}

\def\modi{\pmod_{\A}} \def\modii{\pmod_{\B}}
\def\subi{\A} \def\subii{\B}

\begin{definition}
If we have a morphism $\rho:\modi\to \modii$ between modules with {\inp s} $\met_{\subi}$ and $\met_{\subii}$ respectively, and by setting $\rhof \defeq \rhos$ we obtain a morphism between $(\modi,\modi)$ and $(\modii,\modii)$ as dual pairs, then we call $\rho$ a {\em morphism of inner product modules}\index{morphism!of inner product modules|seealso{isometry}}. In formula this means
\begin{equation}\label{eq:definition inherited bilinear form}
\rhon\met[\X][\Y]_{\subi} = \met[\rhos \X][\rhos \Y]_{\subii}
\end{equation}
for all $\X,\Y\in \modi$. We will also refer to morphisms of inner product modules as {\em isometries}\index{isometry}.

A {\em morphism of {\inp} pairs}\index{morphism!of inner product pairs} $\rho:\A\to \B$, is by definition a morphism between $\A$ and $\B$ as dual pairs, such that at the same time $\rhos:\pmod_{\A}\to\pmod_{\B}$ determines a isometry between inner product modules, i.e. \eqref{eq:dual pair morphism} and \eqref{eq:definition inherited bilinear form} both hold.
\end{definition}

The following propositions and corollaries are proven rather straightforwardly, but are nonetheless worth stating here as we will use it later on.

\def\supp{\B}
\def\subb{\A}

\begin{proposition}\label{prop:surjection_dual_pairs}
Let $\subb$ denote a dual pair. Suppose we have another dual pair $\supp$ and a surjective morphism $\rho:\subb\to \supp$, then it holds that
\begin{align}
\label{eq:Ker_rhos}\pair[\Ker \rhos][\dmod_{\subb}]_{\subb} \subset \Ker \rhon,\\
\label{eq:Ker_rhof}\pair[\pmod_{\subb}][\Ker \rhof]_{\subb} \subset \Ker \rhon.
\end{align}
On the other hand, suppose we have modules $\pmod_{\supp}$ and $\dmod_{\supp}$ over some algebra $\alg_{\supp}$, and suppose we have surjective homomorphisms $\rhon:\alg_{\subb}\to \alg_{\supp}$, $\rhos:\pmod_{\subb}\to\pmod_{\supp}$ and $\rhof:\dmod_{\subb}\to\dmod_{\supp}$, such that \eqref{eq:Ker_rhos} and \eqref{eq:Ker_rhof} hold, then by defining
\begin{equation}\label{eq:def_surj_on_dual_pair}
\pair[\rhos \X][\rhof \W]_{\supp} \defeq \rhon \pair[\X][\W]_{\subb}
\end{equation}
we end up with a dual pair $\supp$, such that $\rho$ becomes a surjective morphism.
\end{proposition}

\begin{corollary}
Let $\subb$ denote an inner product module. Suppose we have another inner product module $\supp$ and a surjective morphism $\rho:\subb\to \supp$, then it holds that
\begin{align}\label{eq:Ker_inner_product}
\pair[\Ker \rhos][\pmod_{\subb}]_{\subb} \subset \Ker \rhon.
\end{align}
On the other hand, suppose we have a module $\pmod_{\supp}$ over some algebra $\alg_{\supp}$, and surjective homomorphisms $\rhon:\alg_{\subb}\to \alg_{\supp}$ and $\rhos:\pmod_{\subb}\to\pmod_{\supp}$, such that \eqref{eq:Ker_inner_product} holds, then by defining
\begin{equation}\label{eq:def_surj_on_inp_mod}
\met[\rhos \X][\rhos \Y]_{\supp} \defeq \rhon \pair[\X][\Y]_{\subb},
\end{equation}
we end up with an inner product module $\supp$, such that $\rho$ becomes a surjective morphism.
\end{corollary}

\begin{corollary}
Let $\subb$ denote an {\inp} pair. Suppose we have another {\inp} pair $\supp$ and a surjective morphism $\rho:\subb\to \supp$, then \eqref{eq:Ker_rhos}, \eqref{eq:Ker_rhof} and \eqref{eq:Ker_inner_product} all hold.

On the other hand, suppose we have modules $\pmod_{\supp}$ and $\dmod_{\supp}$ over some algebra $\alg_{\supp}$, and suppose we have surjective homomorphisms $\rhon:\alg_{\subb}\to \alg_{\supp}$, $\rhos:\pmod_{\subb}\to\pmod_{\supp}$ and $\rhof:\dmod_{\subb}\to\dmod_{\supp}$, such that \eqref{eq:Ker_rhos}, \eqref{eq:Ker_rhof} and \eqref{eq:Ker_inner_product} hold. Then \eqref{eq:def_surj_on_dual_pair} and \eqref{eq:def_surj_on_inp_mod} determine an {\inp} pair $\subb$, such that $\rho$ becomes a surjective morphism.
\end{corollary}


\def\supp{\A}
\def\subb{\B}

\begin{proposition}\label{prop:inj_mor_dual_pairs}
Let $\supp$ denote a dual pair. Suppose we have another dual pair $\subb$ and an injective morphism $\rho:\subb\to \supp$, then it holds that
\begin{equation}\label{eq:injection_dual_pairs}
\pair[\rhos(\pmod_{\subb})][\rhof(\dmod_{\subb})]_{\supp} \subset \rhon(\alg_{\subb}).
\end{equation}
On the other hand, suppose we have modules $\pmod_{\subb}$ and $\dmod_{\subb}$ over some algebra $\alg_{\subb}$, and suppose we have injective homomorphisms $\rhon:\alg_{\subb}\to \alg_{\supp}$, $\rhos:\pmod_{\subb}\to\pmod_{\supp}$ and $\rhof:\dmod_{\subb}\to\dmod_{\supp}$, such that
\eqref{eq:injection_dual_pairs} holds, then by defining
\begin{equation}\label{eq:def_inj_dual_pair}
\pair[\X][\W]_{\subb} \defeq (\rhon)^{-1} \pair[\rhos \X][\rhof \W]_{\supp}
\end{equation}
we end up with a dual pair $\subb$, such that $\rho$ becomes an injective morphism.
\end{proposition}

\begin{corollary}
Let $\supp$ denote an inner product module. Suppose we have another {\inp} module $\subb$ and an injective morphism $\rho:\subb\to \supp$, then it holds that
\begin{equation}\label{eq:injection_inp_modules}
\met[\rhos(\pmod_{\subb})][\rhos(\pmod_{\subb})]_{\supp} \subset \rhon(\alg_{\subb}).
\end{equation}
On the other hand, suppose we have a module $\pmod_{\subb}$ over some algebra $\alg_{\subb}$, and suppose we have injective homomorphisms $\rhon:\alg_{\subb}\to \alg_{\supp}$ and $\rhos:\pmod_{\subb}\to\pmod_{\supp}$, such that \eqref{eq:injection_inp_modules} holds, then by defining
\begin{equation}\label{eq:def_inj_inp_mod}
\met[\X][\Y]_{\subb} \defeq (\rhon)^{-1} \met[\rhos \X][\rhos \Y]_{\supp}
\end{equation}
we end up with an {\inp} module $\subb$, such that $\rho$ becomes an injective morphism.
\end{corollary}

\begin{corollary}
Let $\supp$ denote an {\inp} pair. Suppose we have another {\inp} pair $\subb$ and an injective morphism $\rho:\subb\to \supp$, then \eqref{eq:injection_dual_pairs} and \eqref{eq:injection_inp_modules} both hold.

On the other hand, suppose we have modules $\pmod_{\subb}$ and $\dmod_{\subb}$ over some algebra $\alg_{\subb}$, and suppose we have injective homomorphisms $\rhon:\alg_{\subb}\to \alg_{\supp}$, $\rhos:\pmod_{\subb}\to\pmod_{\supp}$ and $\rhof:\dmod_{\subb}\to\dmod_{\supp}$, such that
\eqref{eq:injection_dual_pairs} and \eqref{eq:injection_inp_modules} both hold, then \eqref{eq:def_inj_dual_pair} and \eqref{eq:def_inj_inp_mod} determine an {\inp} pair $\subb$, such that $\rho$ becomes an injective morphism.
\end{corollary}

\begin{remark}\label{rem:sub_implies_inj}
Note that Proposition \ref{prop:inj_mor_dual_pairs} and its corollaries always apply in case that $\alg_{\subb}$ is a subalgebra of $\alg_{\supp}$ (of course over the same ground ring), and $\pmod_{\subb}$ and $\dmod_{\subb}$ (if present) are submodules of $\pmod_{\supp}$ and $\dmod_{\supp}$ respectively.
\end{remark}

\def\Ideal{\al I}
\def\ideal{\alg_{\B}}

\begin{definition}
Let $\A$ be a dual pair. Another dual pair $\B$ is said to be a {\em subpair}\index{subpair} of $\A$, if $\alg_{\B}$ is a subalgebra of $\A$, $\pmod_{\B}$ and $\dmod_{\B}$ are submodules of $\pmod_{\A}$ and $\dmod_{\A}$ respectively, and the pairing $\pair_{\B}:\pmod_{\B}\otimes_{\alg_{\B}} \dmod_{\B}\to \alg_{\B}$ is inherited from $\A$, i.e. $\pair[\X][\W]_{\B} = \pair[\X][\W]_{\A}$ for all $\X\in \pmod_{\B}$ and $\W\in \dmod_{\B}$. In other words, there exists an injective morphism of dual pairs from $\B$ to $\A$, which is given by inclusion maps. A subpair $\B$ of $\A$ is called an {\em ideal subpair}\index{ideal subpair} if moreover holds that $\alg_{\B}$ is an ideal of $\alg_{\A}$, and $\pmod_{\B}$ and $\dmod_{\B}$ are such that $\pair[\pmod_{\B}][\dmod_{\A}]_{\A} \subset \alg_{\B}$, $\pair[\pmod_{\A}][\dmod_{\B}]_{\A} \subset \alg_{\B}$, and furthermore $\ideal \pmod_{\A} \subset \pmod_{\B}$ and $\ideal \dmod_{\A}\subset \dmod_{\B}$. Moreover, in the particular case that $\pmod_{\B}$ is spanned by $\alg_{\B}\pmod_{\A}$ and $\dmod_{\B}$ is spanned by $\alg_{\B}\dmod_{\A}$, then we call $\B$ an {\em ordinary subpair}\index{ordinary!subpair, submodule}. 

Assume now that $\A$ is an {\inp} module. Then an {\inp} module $\B$ is said to be an {\em {\inp} submodule}\index{inner product submodule} of $\A$, if $\B$ is a subpair of $\A$ when regarded as dual pairs. This is equivalent to saying that there exists an injective isometry from $\B$ to $\A$ given by the inclusion maps. An {\inp} submodule $\B$ of $\A$ is called an {\em ideal submodule}\index{ideal submodule}, if it is an ideal subpair of $\A$ when both are regarded as dual pairs. This is equivalent with requiring that $\pair[\pmod_{\B}][\pmod_{\A}]_{\A} \subset \alg_{\B}$ and $\alg_{\B}\pmod_{\A} \subset \pmod_{\B}$. Like before, if $\pmod_{\B}$ is spanned by $\alg_{\B}\pmod_{\A}$, then $\pmod_{\B}$ is called an {\em ordinary submodule}.

Finally, if $\A$ is an {\inp} pair, then $\B$ is called an (ideal) {\em {\inp} subpair}\index{inner product subpair} of $\A$ if it is an (ideal) subpair of $\A$ if regarded from the viewpoint of dual pairs, and an (ideal) {\inp} submodule of $\A$ if regarded from the viewpoint of {\inp} modules. Again, $\B$ is called an {\em ordinary inner product subpair} if it is ordinary as a subpair.
\end{definition}

The following lemma is an obvious observation, but will be of use later on.

\begin{lemma}\label{prop:inner_product_on_subdualpair}
Suppose $\A$ is a dual pair endowed with an {\inp}, and $\B$ a subpair of $\A$, but over the same algebra as $\A$ (i.e. $\alg_{\B} = \alg_{\A}$). If $\pmod_{\A}$ duplicates $\dmod_{\A}$ and $\pmod_{\B}=\pmod_{\A}$, then $\pmod_{\B}$ also duplicates $\dmod_{\B}$, and likewise if $\dmod_{\A}$ duplicates $\pmod_{\A}$ and $\dmod_{\B}=\dmod_{\A}$, then $\dmod_{\B}$ also duplicates $\pmod_{\B}$.
\end{lemma}

\begin{p}\label{p:musical_implies_ideal_submodule_is_ideal_subpair}
A subpair $\B$ of an inner product pair $\A$ is an inner product pair itself with its inner product inherited from $\A$.
On the other hand, if $\A$ is a musical pair, then any {\inp} submodule $\B$ of the primary module of $\A$ yields an {\inp} subpair of $\A$ by setting $\dmod_{\B} = \form{(\pmod_{\B})}$ and defining its pairing through 
\begin{equation}\label{eq:subpair_submodule}
\met[\X][\form{\Y}]_{\B} = \met[\X][\Y]_{\B},
\end{equation}
for $\X,\Y\in \pmod_{\B}$. Furthermore, such an inner product submodule is an ideal submodule, if and only if its associated subpair is an ideal subpair.
\end{p}


\begin{definition}
Given a dual pair $\A$ and an ideal $\alg_{\B}\subset \alg_{\A}$, then the subpair $\B$ with primary module
$$
\pmod_{\B} \defeq \{ \X \in \pmod_{\A} : \pair[\X][\W] \in \ideal\ \text{for all}\ \W\in \dmod_{\A}\}
$$
and secondary module
$$
\dmod_{\B} \defeq \{ \W \in \dmod_{\A} : \pair[\X][\W] \in \ideal\ \text{for all}\ \X\in \pmod_{\A}\}
$$
will be called the {\em maximal ideal subpair}\index{maximal!ideal subpair} over the algebra $\alg_{\B}$. 

For an inner product module $\A$, we define the {\em maximal ideal submodule}\index{maximal!ideal submodule} to be the submodule
$$
\pmod_{\B} \defeq \{ \X \in \pmod_{\A} : \met[\X][\Y] \in \ideal\ \text{for all}\ \Y\in \pmod_{\A}\}.
$$
\end{definition}

\begin{remark}\label{rem:musical_implies_same_max_sub}
Note that a maximal ideal subpair over $\alg_{\B}$ will contain any other ideal subpair over $\alg_{\B}$. It is also clear that it satisfies $\ideal \pmod_{\A} \subset \pmod_{\B}$ and $\ideal \dmod_{\A}\subset \dmod_{\B}$, as well as $\alg_{\A} \pmod_{\B} \subset \pmod_{\B}$ and $\alg_{\A} \dmod_{\B}\subset \dmod_{\B}$. Furthermore, it should be noted that when $\A$ is a {\musical} pair, then any maximal ideal subpair also determines a maximal ideal submodule, and vice versa. 
\end{remark}

\begin{definition}\label{def:ordinary_dual_pair}
By an {\em ordinary dual pair}\index{ordinary!dual pair, inner product module}, we mean a dual pair for which any maximal ideal subpair is an ordinary subpair. We define an {\em ordinary inner product module} and an {\em ordinary inner product pair} similarly.
\end{definition}

\def\qab{\A/\B}

\begin{definition}\label{def:quotient_pair}
Given a dual pair $\A$ and an ideal subpair $\B$ of $\A$, we define the {\em quotient dual pair}\index{quotient dual pair} $\A/\B$\index{$\mathscr{A}/\mathscr{B}$} by $\alg_{\qab} \defeq \alg_{\A}/\alg_{\B}$, $\pmod_{\qab} \defeq \pmod_{\A}/\pmod_{\B}$, $\dmod_{\qab} \defeq \dmod_{\A}/\dmod_{\B}$ and 
\begin{equation}\label{eq:pairing_on_quotient}
\pair[\X+\pmod_{\B}][\W+\dmod_{\B}]_{\qab} \defeq \pair[\X][\W]_{\A} + \alg_{\B}.
\end{equation}
\end{definition}

\begin{remark}
Note that \eqref{eq:pairing_on_quotient} is well-defined, as this follows directly from Proposition \ref{prop:surjection_dual_pairs} and the requirement that $\pair[\pmod_{\B}][\dmod_{\A}]_{\A} \subset \alg_{\B}$ and $\pair[\pmod_{\A}][\dmod_{\B}]_{\A} \subset \alg_{\B}$. The quotient maps $\rhon:\alg_{\A}\to \alg_{\A}/\alg_{\B}$, $\rhos:\pmod_{\A}\to \pmod_{\A}/\pmod_{\B}$ and $\rhof:\dmod_{\A}\to \dmod_{\A}/\dmod_{\B}$, form together a surjective morphism of dual pairs $\rho:\A\to\A/\B$.
\end{remark}

\begin{proposition}\label{prop:quotient_dual_pair}
If the quotient dual pair $\A/\B$ is non-trivial, then it is non-degenerate if and only if $\B$ is a maximal ideal subpair.
\end{proposition}
\begin{proof}
Note that $\A/\B$ being non-degenerate is equivalent to saying: $\pair[\X][\W]\in\alg_{\B}$ for all $\X\in\pmod_{\A}$ implies $\W \in \dmod_{\B}$, and likewise when the roles of $\pmod_{\A}$ and $\dmod_{\A}$ are switched. But this is exactly the case when $\B$ is a maximal ideal subpair.
\end{proof}

The following corollary is obtained by applying the previous observations to the special case of inner product modules.

\begin{corollary}\label{cor:quotient_inp_mod}
Given an inner product module $\A$ and an ideal submodule $\B$ of $\A$, we define the {\em quotient inner product module $\A/\B$}\index{quotient inner product module} by $\alg_{\qab} \defeq \alg_{\A}/\alg_{\B}$ and $\pmod_{\qab} \defeq \pmod_{\A}/\pmod_{\B}$, and
$$
\met[\X+\pmod_{\B}][\Y+\pmod_{\B}]_{\qab} \defeq \met[\X][\Y]_{\A} + \alg_{\B}.
$$
The quotient maps $\rhon:\alg_{\A}\to \alg_{\A}/\alg_{\B}$ and $\rhos:\pmod_{\A}\to \pmod_{\A}/\pmod_{\B}$, form together a surjective morphism of inner product modules $\rho:\A\to\B/\A$. Furthermore, if $\A/\B$ is non-trivial, then it is non-degenerate if and only if $\B$ is a maximal ideal submodule.
\end{corollary}

\begin{remark}
Note that in case $\B$ is a maximal ideal subpair, the quotient dual pair (or quotient inner product module) $\A/\B$ is non-degenerate even if $\A$ isn't non-degenerate. Hence, by letting $\ideal = \{0\}$, we obtain a recipe for turning a degenerate dual pair (or inner product module) into a non-degenerate one.
\end{remark}

\begin{p}\label{p:quotient_inp_pair}
Similarly as in Definition \ref{def:quotient_pair} and Corollary \ref{cor:quotient_inp_mod}, if $\A$ is an {\inp} pair and $\B$ an ideal {\inp} subpair, then we may form the quotient {\inp} pair $\A/\B$, whose structure as dual pair is induced by $\A/\B$ as quotient of dual pairs, and whose structure as inner product module is induced by $\A/\B$ as quotient of inner product modules. Also, the quotient map of dual pairs $\A\to \A/\B$ now becomes a surjective morphism of inner product pairs. 
\end{p}

\begin{proposition}\label{prop:ind_dupl}
If $\A$ is a {\gen} {\inp} pair, and $\B$ is an ideal subpair of $\A$, then $\A/\B$ is a {\gen} pair as well.
\end{proposition}
\begin{proof}
For any $\Y+\pmod_{\B}\in \pmod_{\qab}$, we have that 
\begin{align*}
&\met[\X+\pmod_{\B}][\Y+\pmod_{\B}]_{\qab} = \met[\X][\Y]_{\A}+\ideal\\
& = \met[\X][\form{\Y}]_{\A}+\alg_{\B} = \pair[\X+\pmod_{\B}][\form{\Y}+\dmod_{\B}]_{\qab},
\end{align*}
for all $\X+\pmod_{\B}\in \pmod_{\qab}$, where $\form{\Y}$ denotes a duplicate of $\Y$ in $\A$. Hence $\dmod_{\qab}$ duplicates $\pmod_{\qab}$. Likewise, for an element $\W+\dmod_{\B}\in\dmod_{\qab}$, we have that $\vec{\W}+\pmod_{\B}$ yields a duplicate of $\W+\dmod_{\B}$. Hence $\qab$ is indeed {\gen}.
\end{proof}

\begin{corollary}\label{cor:quotient_is_musical}
If $\A$ is a {\gen} pair and $\B$ is a maximal ideal subpair of $\A$, then the quotient dual pair $\A/\B$ is a {\musical} pair.
\end{corollary}
\begin{proof}
This follows directly from Corollary \ref{cor:quotient_inp_mod} and Proposition \ref{prop:ind_dupl}.
\end{proof}

\subsection{Projections on inner product modules}

\begin{definition}\label{p:modulo}
Throughout the following, we adopt a convention that can be described as follows: given a dual pair $\A$ and an ideal subpair $\B$, we say that a certain property holds {\em modulo}\index{modulo an ideal subpair} $\B$, if the property holds for the inherited structure on $\qab$, and likewise in the case of inner product modules and inner product pairs. More concretely, we say for instance that a dual pair $\A$ is non-degenerate modulo an ideal subpair $\B$, if $\pair[\X][\W]_{\A} \in \alg_{\B}$ for all $\X\in \pmod_{\A}$ implies $\W\in \dmod_{\B}$, and likewise if the roles of primary and secondary module are reversed (provided that $\qab$ is non-trivial). 
\end{definition}

{
\def\subs{\A^{\!\top}}
\def\subsort{\A^{\!\bot}}
\def\pmodtan{\pmod_{\A^{\!\top}}}
\def\pmodort{\pmod_{\A^{\!\bot}}}
\def\subsB{\B^{\!\top}}
\def\proj{(\cdot)^{\top}}
\def\ortproj{(\cdot)^{\bot}}
\newcommand{\prof}[1]{#1^{\top}}
\newcommand{\ortprof}[1]{#1^{\bot}}
\begin{definition}
Let us be given an inner product module $\A$ and another submodule of $\A$ over the same algebra $\alg_{\A}$, which we formally denote by $\subs$. If a module morphism $\proj:\A\to \subs$ is such that for some ideal $\alg_{\B}$ of $\alg_{\A}$ it holds for all $\X,\Y\in \pmod_{\A}$ that
\begin{equation}\label{eq:defining equation projection}
\met[\X][\prof{\Y}] - \met[\X][\Y] \in \alg_{\B},
\end{equation}
then we call it a {\em projection}\index{projection!modulo an ideal} of $\A$ onto $\subs$ modulo $\alg_{\B}$. If $\alg_{\B} = \{0\}$, we call it simply a \emindex{projection}. Note that a projection $\proj:\A\to \subs$ modulo $\alg_{\B}$ is, in the sense of \ref{p:modulo}, a projection modulo $\B$, for any ideal submodule $\B$ having $\alg_{\B}$ as its underlying algebra. Given a projection $\proj:\A\to \subs$ modulo $\alg_{\B}$, we define the orthogonal complement to $\subs$ modulo $\alg_{\B}$ as the submodule $\subsort$ over $\alg_{\A}$ where
$$
\pmod_{\subsort} \defeq \{\X\in \pmod_{\A} : \met[\X][\Y] \in \alg_{\B} \ \text{for all}\ \Y\in \pmod_{\subs}\}.
$$
Furthermore, a projection modulo $\alg_{\B}$ like above, gives rise to an {\em orthogonal projection}\index{orthogonal projection} $\ortproj:\A\to \subsort$ modulo $\alg_{\B}$, given by
\begin{equation}\label{eq:definition_ort_proj}
\ortprof{\X} \defeq \X - \prof{\X}.
\end{equation}
It follows straightforwardly that the image of $\pmod_{\A}$ lies indeed in $\pmodort$.
\end{definition}

\begin{theorem}\label{thm:projection modulo B}
Let $\A$ be a non-degenerate inner product module, $\subs$ a submodule of $\A$ over the same algebra, and $\B$ a maximal ideal submodule of $\subs$. If there exist a projection $\proj:\A\to \subs$ modulo $\alg_{\B}$, then this projection must be unique modulo $\B$. Furthermore, this projection is a {\em retraction (also called left inverse)}\index{retraction} of the natural inclusion $\subs\to \A$, modulo $\B$.
\end{theorem}
\begin{proof}
What we mean by unique modulo $\B$, is that the composition of the projection with the quotient map $\rho: \subs \to \subs/\B$ is uniquely determined. So suppose now that we have two such projections modulo $\alg_{\B}$, and write $\phi$ for their difference (i.e. $\phi^\sharp:\pmod_{\A}\to \pmod_{\subs}$ is the difference between the two module projections, whereas $\phi^*$ is just the identity on the underlying algebra). Then it follows that $\met[\phi \X][\Y] \in \alg_{\B}$ for all $\Y\in \pmod_{\subs}$, from which follows that $\phi \X\in \pmod_{\B}$, because $\B$ was assumed to be a maximal ideal submodule of $\subs$. But then the composition of $\phi$ with $\rho$ will be the zero map, which proves the statement. 

The additional statement is equivalent to saying that $\X^{\top}-\X\in \pmod_{\B}$ for all $\X\in \pmod_{\subs}$, which is clearly the case as $\met[\X^{\top}-\X][\Y]\in \alg_{\B}$ for all $\Y\in \pmod_{\subs}$. 
%
%
\end{proof}
}

\newcommandtwoopt{\meti}[2][\cdot][\cdot]{\langle #1,#2\rangle_{\!\lra}}
\newcommandtwoopt{\metii}[2][\cdot][\cdot]{\langle #1,#2\rangle_{\!\lra/\ideal}}

\section{Derivations}\label{sec:derivations}

\def\d{\mathrm{d}}

\begin{definition}
A map $\d:\alg\to \mod M$ from a $\ring$-algebra $\alg$ to an $\alg$-module $\mod M$ is called a $\ring$-linear {\em algebra derivation}\index{algebra derivation}\index{derivation!of algebras|seealso{algebra derivation}}, if it is a $\ring$-linear map which satisfies the Leibniz rule:
\begin{equation}\label{EQ Leibniz}
\d(f g) = g\,\d f + f\,\d g.
\end{equation}
We denote the space of all $\ring$-linear algebra derivations from $\alg$ to $\mod M$ by $\Der_{\ring}(\alg,\mod M)$\index{$\mathrm{Der}_{\mathbb{K}}(\mathcal{O},\bf{M})$} or just $\Der(\alg,\mod M)$ if there is no confusion about the intended ring.
The module $\mod M$ is possibly an algebra $\al P$ for which there exists a homomorphism $\phi:\alg \to \al P$, thus making $\al P$ into an $\alg$-module. If $\mod M = \alg$, and $\phi$ is just the identity, then we call the elements of this space simply {\em derivations on $\alg$}, and we denote this space as $\Der_{\ring}(\alg)$ or simply $\Der(\alg)$.
\end{definition}

\begin{lemma}\label{lem:derivation_and_homomorphism}
Given a dual pair, then a derivation 
$$\d\in\Der(\alg,\dmod), \quad f\mapsto \d f,$$
gives rise to a map 
$$\d\in \Hom_{\alg}(\pmod,\Der(\alg)), \quad X\mapsto \d_X.$$
\end{lemma}
\begin{proof}
For $X\in\pmod$, let us define $\d_X:\alg\to \alg$ by $\d_X f \defeq \pair[X][\d f]$. This expression is certainly $\alg$-linear in the argument $X$, and it is moreover a derivation as we have that\begin{multline*}
\d_X (f g) = \pair[X][\d(fg)]
 = \pair[X][f \d g+g \d f]\\ = f\pair[X][\d g]+g\pair[X][\d f] = f \d_X g + g\d_X f.
\end{multline*}
\end{proof}

\begin{definition}\label{def:derivationOnModules}
More generally, if $\mod N$ is an $\alg$-module, then a $\ring$-linear {\em module derivation}\index{module derivation}\index{derivation!of modules|seealso{module derivation}} on $\mod N$ over an algebra derivation $\d\in \Der_{\ring}(\alg)$, is by definition a $\ring$-linear map $\D:\mod N\to \mod N$ for which the following Leibniz rule holds:
\begin{equation}\label{eq:leibnizOnModules}
\D(f V) = \d f\, V + f\, \D V,
\end{equation}
where $f\in \alg$ and $V\in \mod N$. We denote the space of all such module derivations over a derivation $\d$ by $\Der_{\d}(\mod N)$\index{$\mathrm{Der}_{\mathrm{d}}(\bf{N})$}. 
\end{definition}


\def\dif{\D}
{
\def\alghom{\phi}
\begin{proposition}
If $\mod N$ is an $\alg$-module, and $\d$ a derivation on $\alg$, then $\Der(\mod N)$ is an $\alg$-module as well, where scalar multiplication is given by
\begin{equation}\label{eq:module structure Der}
(h \dif)f \defeq h (\dif f)
\end{equation}
for $h, f\in\alg$ and $\dif\in\Der_{\d}(\mod N)$. In particular, if $\D$ is a derivation on $\mod N$ over $\d\in \Der(\alg)$, then for $g\in\alg$ we have that $g\D$ is a derivation on $\mod N$ over $g\,\d$. 
\end{proposition}
\begin{proof}
First we note that $\Der(\alg)$ is a module over $\alg$, where scalar multiplication is defined similarly as in \eqref{eq:module structure Der} (this could also be proven by first adapting the following proof to the special case that $\mod N = \alg$). The space $\Der(\mod N)$ clearly possesses an additive structure under addition of linear maps. Furthermore, if $\D$ is a module derivation on $\mod N$ over $\d\in \Der(\alg)$, then we have that
\begin{align*}
(g\dif)(f V) &= g (\dif(f V)) = g (\d f\, V + f \dif V)\\
&= g (\d f) \, V + f g \dif V = (g\,\d) f \ V + f (g \dif) V,
\end{align*}
so that $g \dif$ is a derivation on $\mod N$ over $g\,\d$.
\end{proof}
}


\begin{proposition}\label{pr:composition derivations}
Suppose we have two derivations in $\Der(\mod N)$, which we formally denote by $\D_X$ and $\D_Y$, which on its turn are derivations over $\d_X$ and $\d_Y$ in $\Der(\alg)$. Then we may compose $\D_X$ and $\D_Y$ to get another map $\D_X \D_Y:\mod N \to
\mod N$. This composition is distributive with respect to addition, and moreover satisfies
$$\D_X(f \D_Y) = (\d_X f) \D_Y + f \D_X \D_Y$$
for all $f\in\alg$ and $\D_X, \D_Y\in \Der(\mod N)$.
\end{proposition}
\begin{proof}
First, we have that 
\begin{multline}
\D_X(\D_Y+\D_Z)V = \D_X(\D_Y\! V+\D_Z\! V)\\= \D_X\D_Y\!V + \D_X\D_Z\!V= (\D_X\D_Y+\D_X\D_Z)V,\end{multline} and likewise we find that 
$$(\D_X+\D_Y)\D_Z = \D_X\D_Z + \D_Y\D_Z.$$ Next, we have that
\begin{multline}
\D_X(f\D_Y)V = \D_X(f\D_Y\! V)\\
= \d_X f\ \D_Y\! V + f \D_X\D_Y\! V = (\d_X f \D_Y + f \D_X\D_Y) V. 
\end{multline}
\end{proof}

In general, the composition of two derivations on a module fails to be a derivation itself. The following proposition tells us that the commutator of $\D_X$ and $\D_Y$ really is a derivation.

\def\X{\D_X}
\def\Y{\D_Y}
\def\Z{\D_Z}
\begin{proposition}\label{ST:existence Lie bracket}
Given $\X, \Y\in \Der(\mod N)$, let us define their commutator\index{commutator} $[\X,\Y]$ as the map from $\mod N$ to $\mod N$, given by
\begin{equation}\label{EQ:Lie bracket}
[\X,\Y]V = \X\Y V-\Y\X V.
\end{equation}
Then, for derivations $\X,\Y\in \Der(\mod N)$ over algebra derivations $\d_X,\d_Y\in \Der(\alg)$, we have that their commutator $[\X,\Y]$ is a derivation on $\mod N$ over the commutator $[\d_X,\d_Y]\in \Der(\alg)$.
\end{proposition}
\begin{proof}
We have
\begin{align*}
[\X,\Y](fV) &= \X(\d_Y f\, V + f \Y\! V)-\Y(\d_X\! f\, V+ f \X\! V) \\
&= \d_X\d_Y\! f\, V+\d_Y\! f\ \X\! V+\d_X\! f\ \Y\! V+f \X\Y\! V\\ 
&\quad -\d_Y\d_X\! f V-\d_X\! f\ \Y\! V-\d_Y\! f\ \X\! V-f \Y\X\! V\\
&= \d_X\d_Y\! f\, V+f \X\Y\! V-\d_Y\d_X\! f V -f \Y\X\! V\\
&= [\d_X,\d_Y]f\, V + f [\X,\Y]V,
\end{align*}
showing that $[\X,\Y]$ is a derivation on $\mod N$ over the algebra derivation $[\d_X,\d_Y]$.
\end{proof}

\begin{theorem}
The derivation space $\Der_{\ring}(\mod N)$, equipped with the above commutator bracket $[\cdot,\cdot]$, forms a Lie algebra over $\ring$.
\end{theorem}
\begin{proof}
We will respectively show that the Lie-bracket is $\ring$-linear, alternating and that it satisfies the Jacobi identity.

$\ring$-linear: the bracket is right-additive because 
\begin{multline*}
[\X,\Y+\Z] = \X(\Y+\Z) - (\Y+\Z)\X\\
= \X\Y - \Y\X + \X\Z - \Z\X = [\X,\Y] + [\X,\Z].
\end{multline*}
Similarly one shows that the bracket is left-additive. The bracket is obviously linear over $\ring$, because the derivations are so.

Alternating: $[\X,\X] = \X\X - \X\X = 0$. In particular this implies that $[\X,\Y] = -[\Y,\X]$.

Jacobi-identity: We have to show that $[\X,[\Y,\Z]] + [\Y,[\Z,\X]] + [\Z,[\X,\Y]]$ equals $0$. Noting that the expression is a cyclic sum, it can be written as
\begin{align*}
& \cycsum{\X,\Y,\Z} [\X,[\Y,\Z]]\\ 
&= \cycsum{\X,\Y,\Z} \X\Y\Z - \X\Z\Y -\Y\Z\X + \Z\Y\X,
\end{align*}\index{$\sum^{\mathrm{cycl}}$}
where the summations are taken cyclically over $\X$, $\Y$ and $\Z$. 

Because we have
\[
\cycsum{\X,\Y,\Z} \X\Y\Z = \cycsum{\X,\Y,\Z} \Y\Z\X\quad\mbox{and}\quad \cycsum{\X,\Y,\Z} \X\Z\Y = \cycsum{\X,\Y,\Z} \Z\Y\X,
\]
it follows that this expression equals $0$.
\end{proof}

\begin{proposition}
Given a derivation $\D_X\in \Der(\mod N)$ over $\d_X\in\Der(\alg)$, then the adjoint map\index{adjoint map} $$\ad_{\D_X}:\Der(\mod N)\to \Der(\mod N),\ \D_Y\mapsto [\D_X,\D_Y]$$\index{$\mathrm{ad}$} 
is a derivation on $\Der(\mod N)$ over $\d_X$.
\end{proposition}
\begin{proof}
Using Proposition \ref{pr:composition derivations}, we obtain:
\begin{align*}
\ad_{\X} (f \Y) &= [\X,f\Y] = \X(f\Y)-f\Y\X\\
&= \d_X f\, \Y + f\X\Y - f\Y\X \\
 &= \d_X f\, \Y + f[\X,\Y] = \d_X f\, \Y + f\ad_{\X} \Y.
\end{align*}
\end{proof}

\begin{remark}
From the above we also deduce that for derivations $\D_X$ and $\D_Y$ over $\d_X$ and $\d_Y$, and $f,g\in \alg$, we have that:
\begin{equation}\label{eq:lie bracket two functions}
[f\X,g\Y] = fg[\X,\Y]+f \d_X g\, \Y-g \d_Y f\, \X.
\end{equation}
\end{remark}

\section{Rinehart spaces}\label{sec:Rinehart spaces}

\begin{definition}\label{def:LR_algebra}
By a {\em \lrs}\index{Lie-Rinehart algebra} (cf. \cite{Hu04}), we mean an object consisting of the following data:
\begin{itemize}
\item an algebra $\alg$ over a ring $\ring$,
\item a Lie algebra $\mod M$ over $\ring$, which is also a module over $\alg$,
\item an $\alg$-linear map $\d$\index{$\mathrm{d}$} given by
$$\d:\mod M \to \Der(\alg),\ X\mapsto \d_X,$$
such that for all $X\in \mod M$ the adjoint map 
$$\ad_{X}:\mod M \to \mod M,\quad Y \mapsto [X,Y]$$
is a derivation on $\mod M$ over $\d_X$.
\end{itemize}

\end{definition}

\begin{p}
The last condition of Definition \ref{def:LR_algebra} is equivalent to saying that
\begin{equation}\label{eq:condition_LR_alg}
[X,f Y] = \d_X f\, Y+f[X,Y],
\end{equation}
for all $X,Y\in \mod{M}$ and $f\in\alg$. The map $\d:\mod{M}\to \Der(\alg)$ is commonly referred to as the {\em anchor map}\index{anchor map}. We may also regard it as a map 
$\d:\alg \to \dual{\mod{M}}, f\mapsto \d f,$
which is implicitly determined by requiring $\pair[X][\d f] = \d_X f$, where $\pair$ stands for the natural pairing between $\mod{M}$ and $\dual{\mod{M}}$. 
Regarded from this viewpoint, we refer to $\d$ as the {\em differential operator}\index{differential operator}.
\end{p}




\begin{definition}
By a {\em non-degenerate {\lrm}}\index{Rinehart space}, we mean a dual pair $\A$ whose primary and secondary module we denote by $\derx = \derx_{\A}$\index{$\mathfrak{X}_{\mathscr{A}}$} and $\derf = \derf_{\A}$\index{$\Omega_{\mathscr{A}}$} respectively, and for which there exists a $\d\in \Der(\alg,\derf)$\index{$\mathrm{d}$}, such that the associated map $\d\in\Hom_{\alg}(\derx,\Der(\alg))$ (cf. Lemma \ref{lem:derivation_and_homomorphism}) turns $\alg$ and $\derx$ into a {\lrs}. Furthermore we impose the following extra conditions:
\begin{enumerate}
\item $\derx$ and $\derf$ are 
finite-dimensional modules over $\alg$,
\item $\derf$ is spanned by the image of the map $\d:\alg\to \derf$, and
\item the pairing between $\derx$ and $\derf$ is non-degenerate.
\end{enumerate}
In the following, by a {\em Rinehart space} we will always mean a non-degenerate Rinehart space. In occasions that all but the last of the extra conditions (the one on non-degeneracy) is satisfied, we will call it explicitly a {\em degenerate Rinehart space}. 

We will use the following terminology to denote the components of a {\lrm}: we refer to $\alg$ as its {\em function space}\index{function space of Rinehart space}, $\derx$ as its {\em \xspace}\index{tangent space of a Rinehart space}, $\derf$ as its {\em \omspace}\index{cotangent space of a Rinehart space} and $\d:\alg\to \derf$ as the {\em differential operator}\index{differential operator}. Elements of $\alg$, $\derx$ and $\derf$ are respectively referred to as {\em functions}\index{function on Rinehart space}, {\em (tangent) vector fields}\index{tangent vector field on a Rinehart space} and 
{\em differential 1-forms}\index{differential 1-form on a Rinehart space}. It should be remarked that this is all just formal terminology, as in some examples the elements of $\alg$ in strict sense are not functions, and similarly for $\derx$ and $\derf$.
\end{definition}

\begin{p}\label{p:extra_conditions}
The second of the extra requirements above, means that the only linear subspace of $\derf$ containing $\d f$ for each $f\in\alg$, is the space $\derf$ itself. Recall that for any $\alg$-module, one can construct the smallest submodule containing some set $V$ by taking the submodule of finite linear combinations of elements from $V$. Hence, by the last requirement, any differential 1-form $\omega\in\derf$ can be written as a finite linear combination of differentials, i.e. $\omega = g_1 \d f_1 + \ldots + g_n \d f_n$, where $f_i,g_i\in\alg$, and $n$ is a positive integer depending on $\omega$. 

It should be noted that among the statements we are going to prove, there are quite some which do not depend on (all) the three extra requirements. 
However, although some theoretical results indeed remain valid under less restrictive conditions, the above requirements will turn out to be quite naturally satisfied for the particular examples of {\lrm s} we aim to cover.
\end{p}

\begin{definition}
A {\em morphism}\index{morphism!of Rinehart spaces} of (possibly degenerate) {\lrm s} ${\rho:\A\to\B}$, by definition consists of a morphism of the underlying dual pairs $(\derx_{\A},\derf_{\A})$ and $(\derx_{\B},\derf_{\B})$ (cf. Definition \ref{def:morphism_dual_pairs}), such that $\rhos:\derx_{\A}\to\derx_{\B}$ is moreover a Lie algebra homomorphism, and the following diagram commutes:
  \begin{center}
\begin{tikzcd}
\alg_{\A} 
\arrow[r, "\d^{\A}"]
\arrow[d, "\rhon"]
& \Omega_{\A}
\arrow[d, "\rhof"]\\
\alg_{\B} 
\arrow[r, "\d^{\B}"]
& \Omega_{\B} 
\end{tikzcd}
  \end{center}
  

\end{definition}



\begin{p}
The kernel of the differential operator $\d$ yields a subalgebra of $\alg$. We refer to the elements in this kernel as {\em constant functions}\index{constant!function on Rinehart space}. If $\alg$ is a unital algebra, then it follows that $1\in\alg$ is a constant function, for we have that $\d(1) = \d(1\cdot 1) = 1\, \d(1) + 1\, \d(1) = \d(1) + \d(1)$, so that subtracting $\d(1)$ from both sides yields $\d(1) = 0$. Note that if $\alg$ is a unital algebra over $\ring$, then $\ring$ forms a subalgebra of $\alg$, and because $\d$ is assumed to be $\ring$-linear it follows that any element of $\ring$ gives rise to a constant function in $\alg$. Regardless of the ground ring, we have that any integer multiple of $1\in\alg$ must be a constant function.  
\end{p}


\def\otim{\otimes}
\begin{definition}
Suppose we have a {\lrm} $\A$ and some module $\mod N$ over $\alg_{\A}$. An {\em affine connection}\index{affine connection} (or {\em connection}\index{connection|seealso{affine connection}} for short) on $\mod N$ over $\A$ is by definition an $\alg_{\A}$-linear map 
$$
\nabla:\derx\to \Der(\mod N), \ X\mapsto \nabla_{X},\index{$\nabla$}
$$
such that $\nabla_X$ is a derivation on $\mod N$ over $\d_{X}\in\Der(\alg)$, where $\d = \di{\A}$ denotes the differential operator of the {\lrm} $\A$. In such case, $\mod N$ is called an {\em affinely connected module over $\A$}\index{affinely connected module}. In particular, a connection satisfies the property
\begin{equation}\label{eq:connection_property}
\nabla_X (f N) = \d_X f\ N + f\nabla_X N,
\end{equation}
where $X\in \derx_{\A}$, $N\in \mod N$ and $f\in \alg_{\A}$.
\end{definition}

\def\modA{{\mod N}_{\A}}
\def\modB{{\mod N}_{\B}}
\begin{definition}
Let $\modA$ and $\modB$ be affinely connected modules over $\A$ and $\B$, with connections $\nabla^{\A}$ and $\nabla^{\B}$ respectively. We call $\Phi:\modA \to \modB$ a {\em morphism of affinely connected modules}\index{morphism!of affinely connected modules}, if there exists a associated underlying morphism of {\lrm s} $\phi:\A \to \B$ such that 
$$\Phi (\nabla^{\A}_{X_{\A}} N_{\A}) = \nabla^{\B}_{X_{\B}} N_{\B},$$
where $X_{\A}\in \derx_{\A}$, $N_{\A}\in \modA$ and $X_{\B}\in \derx_{\B}$, $N_{\B}\in \modB$ are their images under the map $\phi^{\sharp}$ and $\Phi$ respectively.
\end{definition}

The following is a straightforward observation that we will use later on.

\begin{proposition}\label{prop:inducing connection on isomorphic modules}
Let $\A$ be a {\lrm}, and $\mod N$ an affinely connected module over $\A$. If $\mod N'$ is another $\alg_{\A}$-module which is isomorphic with $\mod N$, then a connection on $\mod N'$ exists such that the isomorphism between $\mod N$ and $\mod N'$ becomes an isomorphism of affinely connected modules.
\end{proposition}



\begin{definition}\label{def:curvature tensor}
To any affinely connected module $\mod N$ over a {\lrm} $\A$, one can associate a so-called {\em curvature tensor}\index{curvature tensor}, which is an $\alg_{\A}$-linear map $R:\derx\otim\derx\otim\mod N\to \mod N$ given by\index{$R$}
\begin{equation}\label{eq:curvature_tensor}
R(X,Y)Z = [\nabla_X,\nabla_Y]Z - \nabla_{[X,Y]}Z.
\end{equation}
\end{definition}

Here, $R(X,Y):\mod N\to \mod N$ stands for the $\alg_{\A}$-linear map which more formally could be described as $(X\otimes Y)\lrcorner R$, where $R$ is regarded as a tensor map $R:(\derx\otim \derx)\otim \mod N\to\mod N$ (cf. Definition \ref{def:non-degenerate}). The claim that $R$ is a tensor actually requires a proof.

\begin{proposition}
The operator $R$ defined by \eqref{eq:curvature_tensor} is a tensor.
\end{proposition}
\begin{proof}
To see that $R$ is tensorial in $Y$ (and by a likewise argument also in $X$), note that
$$[\nabla_X,f \nabla_Y] = f[\nabla_X,\nabla_Y] + \d_X f\, \nabla_Y$$ and that
from $[X, f Y] = f[X,Y] + \d_X f\, Y$ it follows that
$$\nabla_{[X, f Y]} = f\nabla_{[X,Y]} + \d_X f\, \nabla_Y.$$
Hence we have that
$$[\nabla_X,f \nabla_Y] - \nabla_{[X, f Y]} = f\cdot([\nabla_X,\nabla_Y] - \nabla_{[X,Y]}).$$
To see that $R$ is tensorial in $Z$, recall that by Proposition \ref{ST:existence Lie bracket} we have that $[\nabla_X,\nabla_Y]$ is a derivation over $\d_{[X,Y]}$, and because $\nabla_{[X,Y]}$ is also a derivation over $\d_{[X,Y]}$, it follows that the difference between $[\nabla_X,\nabla_Y]$ and $\nabla_{[X,Y]}$ becomes $\alg_{\A}$-linear.
\end{proof}

\begin{definition}
In case the curvature tensor associated to a connection is just the zero tensor, then we call it a {\em flat connection}\index{flat connection}.
\end{definition}

\begin{definition}
A {\em covariant derivative}\index{covariant derivative} on a {\lrm} $\A$ is by definition a connection on its tangent space $\derx = \derx_{\A}$ over $\A$ itself, i.e. an $\alg_{\A}$-linear map $\nabla:\derx\to \Der(\derx)$, such that $\nabla_X$ is a derivation over $\d_X$. A covariant derivative is called {\em symmetric}\index{symmetric!covariant derivative} if
\begin{equation}
[X,Y] = \nabla_X Y - \nabla_Y X
\end{equation}
for all $X,Y\in \derx$.
\end{definition}


\def\M{\lra{M}}

\begin{definition}\label{def:compatible_with_pairing}
If $\A$ is a {\lrm} and $\B = (\pmod,\dmod)$ is a dual pair over $\alg_{\A}$, then an affine connection on $\B$ over $\A$ consists by definition of a connection on $\pmod$ over $\A$ and a connection on $\dmod$ over $\A$, both of which we denote by $\nabla = \nabla^{\B}$. We call this connection {\em compatible}\index{compatible} with the pairing, if it holds that 
$$\d_X \pair[M][U]_{\B} = \pair[\nabla_X M][U]_{\B} + \pair[M][\nabla_X U]_{\B},$$
for all $X\in \derx_{\A}$, $M\in \mod M$, $W\in \mod W$ and where $d$ denotes the differential operator associated to $\A$.
For the special case that $\pmod = \dmod$, i.e. if $\B$ is an inner product module, a similar formula applies, in which case we say that the connection is compatible with the inner product.
\end{definition}


\subsection{First examples}

\begin{example}\label{ex:real manifold}
Any smooth $n$-dimensional real manifold $M$ gives rise to a {\lrm} where $\ring = \R$, $\alg$ is the space $C^\infty(M)$ of smooth functions on $M$ and $\derx$ is the space $C^\infty(M,TM)$\index{$C^\infty$} of smooth vector fields on $M$, by which we mean the smooth sections of the tangent bundle of $M$, and $\derf$ the space $C^\infty(M,T^*\!M)$ of smooth differential 1-forms on $M$. If $M$ is an analytic manifold, one may restrict $\alg$ to the space of analytic functions $C^\omega(M)$ and $\derx$ and $\derf$ to the space of analytic vector fields $C^\omega(M,TM)$\index{$C^\omega$} and analytic differential 1-forms $C^\omega(M,T^*\!M)$, respectively. In both cases the dimension of this {\lrm} will be $n$.
\end{example}

\begin{example}\label{ex:complex manifold}
On a complex $n$-dimensional manifold $M$, we may take $\ring = \C$, $\alg = H^0(M)$\index{$H^0$}, the space of holomorphic functions on $M$, $\derx = H^0(M,TM)$, the space holomorphic vector fields on $M$, i.e. the holomorphic sections of the holomorphic tangent bundle $TM$, and $\derf = H^0(M,T^*\!M)$, the space of holomorphic differential 1-forms. This will give an $n$-dimensional {\lrm}. Regarding $M$ as a smooth real $2n$-dimensional manifold, it is also possible to let $\alg = \C\otimes C^\infty(M)$, the complex-valued smooth functions on $M$, $\derx = C^\infty(M,\C\otimes TM)$, the smooth sections of the complexified tangent bundle of $M$, and $\derf = \C^\infty(M,\C\otimes T^*\!M)$. This gives a {\lrm} which is $2n$-dimensional.
\end{example}

\begin{example}
Many other examples are obtained by varying the ground ring $\ring$, and picking a suitable class of functions with respect to the ground ring chosen. For example, in a similar way as in Example \ref{ex:complex manifold}, we may construct {\lrm s} corresponding to paraholomorphic or paracomplex manifolds. These resemble the {\lrm s} of the previous example, except that the ground ring is in this case the ring of split-complex numbers, with function space and derivation space adapted accordingly (e.g. \cite{Ch}). Perhaps another example worth mentioning is the case where $\ring$ is the field of $p$-adic numbers (e.g. \cite{Sc}). 
\end{example}


\begin{example}\label{ex:polynomial ring}
As remarked before, the `function space' $\alg$ does not need to consists of functions in a strict sense. For example, given a smooth real manifold $M$, one could consider functions that are smooth on a dense open subset of $M$, and then let $\alg$ consist of equivalence classes of functions, where two functions are regarded equivalent if they coincide on a dense open subset. Other examples include the space of multivariate formal power series or multivariate polynomials over a ring $\ring$. In the case of multivariate polynomials over $\R$ or $\C$, the resulting {\lrm} is not that different from certain {\lrm s} described in Example \ref{ex:real manifold} and \ref{ex:complex manifold}, except that the function space $\alg$ is now restricted to only polynomial functions. But for example for the finite field $\F_q$, we have that the polynomial ring $\F_q[x_1,\ldots,x_n]$ cannot be embedded in the space of functions from $\F_q^n$ to $\F_q$. As will be become clear later on, {\lrm s} also provide a framework to do Riemannian geometry with a finite field as ground ring.
\end{example}


\subsection{Riemannian Rinehart spaces}

\begin{definition}
By a \emindex{Riemannian metric} on a {\lrm} $\A$, we simply mean an {\inp} on its {\xspace} $\derx$, i.e. a symmetric bilinear map $\met:\derx\otimes_{\alg} \derx\to \alg$. We call the resulting inner product pair to which this leads, a {\em Riemmanian {\lrm}}\index{Riemannian Rinehart space}. By a {\em musical (resp. non-degenerate, \gen)}\index{musical!Riemannian metric}\index{non-degenerate!Riemannian metric}\index{duplicative!Riemannian metric} Riemannian metric, we mean a Riemannian metric which is musical (resp. non-degenerate, \gen) as inner product on the dual pair $\A$. 
\end{definition}

\begin{definition}
By a {\em morphism (or isometry) of Riemannian {\lrm s}}\index{morphism!of Riemannian Rinehart spaces|seealso{isometry}}\index{isometry!of Riemannian Rinehart spaces} $\rho:\A\to\B$, we mean a morphism of {\lrm s} preserving the Riemannian metric, i.e. $\rhos:\derx_{\A}\to\derx_{\B}$ should be an isometry.
\end{definition}

\begin{definition}
Given a musical Riemannian \lrm, we define the \emindex{gradient} as the map $\grad:\alg\to \derx$\index{$\mathrm{d}^{\sharp}$}, 
given by $\grad f \defeq (\d f)^\sharp$. By construction, it is the unique derivation in $\derx$ such that $\met[\grad f][X] = \d_Xf$ for all $X\in \derx$ (see Proposition \ref{prop:musical isomorphism} and Definition \ref{def:inverse canonical isomorphism}).
\end{definition}

\begin{p}\label{par:nabla surjective}
It should be noted that for a musical Riemannian {\lrm}, the gradients must necessarily span the entire space $\derx$, i.e. if a linear subspace of $\derx$ contains $\grad f$ for each $f\in\alg$, then it is the whole space $\derx$.
\end{p}


\begin{definition}\label{def:LC_connection}
Let $\A$ be a Riemannian {\lrm} and let $\nabla$ be a symmetric covariant derivative on $\A$. We call $\nabla$ a \emindex{Levi-Civita connection} if it is moreover compatible with the Riemannian metric on $\derx$, i.e. if
\begin{equation}\label{eq:compatible}
\d_X\met[Y][Z] = \met[\nabla_X Y][Z]+\met[Y][\nabla_X Z],
\end{equation}
for all $X,Y,Z\in \derx$. A Riemannian {\lrm} endowed with a Levi-Civita connection, we briefly call a {\em Levi-Civita {\lrm}}\index{Levi-Civita Rinehart space}.
\end{definition}

\begin{proposition}
On a musical Riemannian {\lrm}, any affine connection on the tangent space $\derx$ is naturally transferable to an affine connection on the cotangent space $\derf$, such that the musical isomorphisms become isomorphisms of affinely connected modules over $\A$, i.e.
$$\nabla_X Y^\flat = (\nabla_X Y)^\flat\qquad \text{or equivalently}\qquad \nabla_X \omega^\sharp = (\nabla_X \omega)^\sharp$$
where $X, Y\in\derx$ and $\omega\in\derf$. Moreover, if the affine connection on $\derx$ is compatible with the Riemannian metric, then the affine connections on $\derx$ and $\derf$ together are compatible with their pairing in the sense of Definition \ref{def:compatible_with_pairing}.
\end{proposition}
\begin{proof}
The first part follows directly from Proposition \ref{prop:inducing connection on isomorphic modules}.
For the second claim, we calculate
\begin{align*}
\d_X\met[Y][\omega] &= \d_X\met[Y][\vec{\omega}] = \met[\nabla_X Y][\vec{\omega}]+\met[Y][\nabla_X \vec{\omega}]\\
 &= \met[\nabla_X Y][\omega] + \met[Y][\vec{(\nabla_X \omega)}] = \met[\nabla_X Y][\omega] + \met[Y][\nabla_X \omega].
\end{align*}
\end{proof}
The following theorem is in a sense a generalization of the fundamental theorem of Riemannian geometry, which says that on every Riemannian manifold there exists a unique Levi-Civita connection.

\begin{theorem}\label{thm:fundamental theorem}
Let $\lra$ be a Riemannian \lrm, and suppose that $2$ is a unit in $\ring$. If the Riemannian metric is non-degenerate, then any Levi-Civita connection on $\lra$ is necessarily unique. If the Riemannian metric is moreover {\gen} and $(\derx_{\A},\derf_{\A})$ is isomorphic with $(\derx_{\A},\dual \derx_{\A})$ as dual pairs, then a Levi-Civita connection on $\lra$ certainly exists.
\end{theorem}
\begin{proof}
To prove the first claim, assume that we have a Levi-Civita connection $\nabla$. Using that the inner product is symmetric and that $\nabla$ is compatible with the inner product, we find the following equations:
\begin{align*}
\cycsum{X,Y,Z} \d_X\met[Y][Z] &= \cycsum{X,Y,Z} \met[\nabla_X Y][Z] + \met[\nabla_X Z][Y]\\
&= \cycsum{X,Y,Z} \met[\nabla_X Y][Z] + \cycsum{X,Y,Z} \met[\nabla_Y X][Z]\\
\cycsum{X,Y,Z} \met[X][{[Y,Z]}] &= \cycsum{X,Y,Z} \met[X][\nabla_Y Z - \nabla_Z Y]\\
&= \cycsum{X,Y,Z} \met[\nabla_X Y][Z] - \cycsum{X,Y,Z} \met[\nabla_Y X][Z].
\end{align*}
Summing the equations above, we see that
\begin{align*}
\cycsum{X,Y,Z} \Big( \d_X\met[Y][Z] + \met[X][{[Y,Z]}] \Big)
&= 2\cycsum{X,Y,Z} \met[\nabla_X Y][Z]\\ 
&= 2\Big(\met[\nabla_X Y][Z] + \met[\nabla_Y Z][X] + \met[\nabla_Z X][Y]\Big)\\
&= 2\Big(\met[\nabla_X Y][Z] + \d_Z\met[X][Y] + \met[X][{[Y,Z]}]\Big),
\end{align*}
so that
\begin{equation}\label{eq:Koszul}
\met[\nabla_X Y][Z] = \frac{1}{2}\cycsum{X,Y,Z} \Big(\d_X\met[Y][Z] + \met[X][{[Y,Z]}]\Big) -\d_Z\met[X][Y] - \met[X][{[Y,Z]}].
\end{equation}
Now if $\tilde\nabla$ would be another Levi-Civita connection, it would have to satisfy the above equation for the same reasons. Hence by the additivity of the inner product, we find that for any $X$, $Y$ and $Z$ in $\derx$, we have 
$$\met[\nabla_X Y-\tilde \nabla_X Y][Z] = 0.$$
Because $\met$ is non-degenerate, this can only hold if $\nabla_X Y=\tilde \nabla_X Y$ for any $X$ and $Y$, which proves uniqueness.
To prove existence, consider the map $\form\nabla: \derx \times \derx \to \derf$, 
implicitly defined by
\begin{equation}\label{EQ:KoszulLinForm}
\pair[\form\nabla_X Y][Z] = \frac{1}{2}\cycsum{X,Y,Z} \Big(\d_X\met[Y][Z] + \met[X][{[Y,Z]}]\Big) -\d_Z\met[X][Y] - \met[X][{[Y,Z]}],
\end{equation}
where $\form\nabla_X Y$ stands for $(\form\nabla Y)X$. A straightforward calculation shows that the above expression is $\ring$-linear in $Y$ and $\alg$-linear in $X$ and $Z$. Furthermore, one verifies that
\begin{equation*}
\pair[\form\nabla_X(fY)][Z] = f\, \pair[\form\nabla_X Y][Z] + \d_X f\, \met[Y][Z],
\end{equation*} 
i.e. we have $\form\nabla_X(fY) = f\, \form\nabla_X Y + \d_X f\ Y^\flat$. Now define $\nabla$ to be the composition of $\form\nabla$ with the musical isomorphism $(\cdot)^\sharp: \derf\to \derx$, i.e. 
$$\nabla_X Y = (\form\nabla_X Y)^\sharp$$
for all $X,Y\in\derx$. The map $\nabla: \derx \times \derx \to \derx$ then satisfies \eqref{eq:connection_property}, and hence can be regarded as a covariant derivative $\nabla: \derx \to \Der(\derx)$. By its very construction it also satisfies \eqref{eq:Koszul}, so that it becomes straightforward to verify that $\met[\nabla_X Y][Z] + \met[\nabla_X Z][Y] = \d_X\met[Y][Z]$ and
$\met[\nabla_X Y][Z] - \met[\nabla_Y X][Z] = \met[{[X,Y]}][Z]$, implying that $\nabla_X Y - \nabla_Y X = [X,Y]$.
\end{proof}

\begin{remark}
The condition that $2$ is a unit in $\ring$ is not strictly necessary for the existence of a Levi-Civita connection. As we will see soon, there are {\lrm s} for which $2$ is not a unit in $\ring$, but which do posses a Levi-Civita connection. 
We do not know, however, if Theorem \ref{thm:fundamental theorem} still holds if the condition that $2$ is a unit in $\ring$ is dropped altogether. Another rightful question unanswered here, is whether the condition that $(\derx_{\A},\derf_{\A})$ is isomorphic with $(\derx_{\A},\dual \derx_{\A})$ could be loosened.
\end{remark}

\begin{definition}
We refer to the curvature tensor associated to a Levi-Civita connection (see \ref{def:curvature tensor}) as the \emindex{Riemann tensor}. We say that a {\lrm} with a Levi-Civita connection $\nabla$ has {\em constant sectional curvature}\index{constant!sectional curvature} $c\in\ring$, if its Riemann tensor satisfies
\begin{equation}\label{eq:constant sectional curvature}
R(X,Y)Z = c\,(\met[Y][Z] X - \met[X][Z] Y).
\end{equation}
for any $X,Y,Z\in \derx$. A {\lrm} with constant sectional curvature is called a \emindex{space form}.
\end{definition}

\subsection{Quotient Rinehart spaces}\label{sec: LR submanifolds}

\begin{p}
In this part we will give a description of submanifolds in the framework of {\lrm s}. In this setting, we regard an embedding of manifolds from the perspective of the associated homomorphism between their algebras of functions, namely the restriction of functions on the ambient manifold to functions on the submanifold. If we restrict the function space of the submanifold to those functions that are thus obtained from a function on the ambient space, then the interrelating homomorphism is surjective. Recall that when we have a surjective homomorphism of $\ring$-algebras $\rhon:\alg_{\A} \to \alg_{\c}$, then the isomorphism theorem of algebra homomorphisms says that $\alg_{\B} \defeq \Ker \rhon$ is an ideal, and that $\alg_{\c}$ is isomorphic with $\alg_{\A}/\alg_{\B}$. Hence, in the setting of {\lrm s}, the study of submanifolds is translated to the study of quotient homomorphisms between their function spaces, as well as the related morphisms between their (co)tangent spaces.
\end{p}

\begin{theorem}\label{thm:quotient_lrm}
Suppose $\A$ is a (possibly degenerate) {\lrm} and $\B$ 
a maximal ideal subpair of $\A$, such that $\derx_{\B}$ is a Lie ideal of $\derx_{\A}$ (but possibly with different underlying algebra) and $\d \alg_{\B} \subset \derf_{\B}$. Then by defining the differential operator $\d^{\qab}:\alg_{\qab}\to \derf_{\qab}$ as
\begin{equation}\label{eq:dif_op_on_quotient}
\d^{\A/\B} (f+\alg_{\B}) \defeq \d^{\A} f + \derf_{\B},
\end{equation}
the quotient pair $\A/\B$ becomes a {\lrm}, and the quotient map of dual pairs $\rho:\A\to \A/\B$ becomes a surjective morphism of (possibly degenerate) {\lrm s}.
\end{theorem}
\begin{proof}
Because $\derx_{\B}$ is a Lie ideal of $\derx_{\A}$, this means that $[\derx_{\B},\derx_{\A}]\subset \derx_{\B}$, and thus on the primary module $\derx_{\A}/\derx_{\B}$ of $\qab$, a Lie algebra structure is induced by
$$[X+\derx_{\B},Y+\derx_{\B}]_{\qab} = [X,Y]_{\A}+\derx_{\B}.$$
By taking \eqref{eq:dif_op_on_quotient} as definition for the differential operator, one verifies that a.o. \eqref{eq:condition_LR_alg} will be satisfied, and hence turns $\qab$ into a {\lrm}, such that the quotient map is a surjective morphism.
\end{proof}

The commutative diagram displayed here, will be helpful in discussing the notion of quotient {\lrm s}.

\def\Bnn{\B^{0}_0} 
\def\Bno{\B^0_{\bot}}
\def\Btn{\B^{\top}_{0}}
\def\Bt{\Btn}
\def\Bo{\Bno}
\def\eBnn{(\derx_{\B}^{0},\derf_{\B}^0)} 
\def\eBno{(\derx_{\B}^{0},\derf_{\B}^{\bot})} 
\def\eBtn{(\derx_{\B}^{\top},\derf_{\B}^0)}

\def\Ann{\A^0_0}
\def\Atp{\A^{\top}}
\def\Apo{\A_{\bot}}
\def\Att{\A^{\top}_{\top}}
\def\Aoo{\A^{\bot}_{\bot}}
\def\App{\A}
\def\At{\Atp}
\def\Ao{\Apo}
\def\eAtt{(\derx_{\A}^{\top},\derf_{\A}^{\top})} 
\def\eAoo{(\derx_{\A}^{\bot},\derf_{\A}^{\bot})} 
\def\eAtp{(\derx_{\A}^{\top},\derf_{\A})} 
\def\eApo{(\derx_{\A},\derf_{\A}^{\bot})} 
\def\eApp{(\derx_{\A},\derf_{\A})}

\def\Coo{\c^{\bot}_{\bot}}
\def\Cpp{\c}
\def\Ctt{\c^{\top}_{\top}}
\def\Ct{\Ctt}
\def\Co{\Coo}
\def\eCoo{(\derx_{\c}^{\bot},\derf_{\c}^{\bot})} 
\def\eCpp{(\derx_{\c},\derf_{\c})} 
\def\eCtt{(\derx_{\c}^{\top},\derf_{\c}^{\top})}

\begin{table}[!ht]
  \begin{center}
\begin{tikzcd}
\Bt = \eBtn 
\arrow[r, hookrightarrow] 
& \Ao = \eApo 
\arrow[r, "\rho", twoheadrightarrow] 
\arrow[d, hookrightarrow]
& \Co = \eCoo \\
\Bnn = \eBnn 
\arrow[r, hookrightarrow] 
\arrow[d, hookrightarrow] 
\arrow[u, hookrightarrow]
& \A = \eApp 
\arrow[r, "\rho", twoheadrightarrow] 
& \c = \eCpp \\
\Bo = \eBno 
\arrow[r, hookrightarrow] 
& \At = \eAtp 
\arrow[u, hookrightarrow]
\arrow[r, "\rho", twoheadrightarrow]
& \Ct = \eCtt 
\end{tikzcd}
  \end{center}
 \caption{}
 \label{tab:cdi}
\end{table}

Let $\A$ be a {\lrm} and $\alg_{\B}$ an ideal of $\alg_{\A}$. By $\derf_{\A}^{\bot}$\index{$\Omega_{\mathscr{A}}^{\bot},\Omega_{\mathscr{A}}^0$} we denote the span of $\d \alg_{\B}$ over the algebra $\alg_{\A}$, i.e. the span of all differentials $\d f$ with $f\in \alg_{\B}$. We let $\derx_{\A}^\top$\index{$\mathfrak{X}_{\mathscr{A}}^{\top},\mathfrak{X}_{\mathscr{A}}^0$} denote the submodule of $\derx_{\A}$ given by
$$
\derx_{\A}^\top \defeq \{ X\in \derx_{\A} : \pair[X][\derf_{\A}^{\bot}]\subset \alg_{\B}\}. 
$$
In other words, $\derx^{\top}_{\A}$ is given by all tangent vector fields $X\in \derx_{\A}$ such that $\d_X \alg_{\B}\subset \alg_{\B}$. Geometrically, one can think of $\derx_{\A}^\top$ as those tangent vector fields on the surrounding space, which at the variety defined by $\alg_{\B}$ are tangent. Furthermore we define
\begin{align*}
\derx_{\A}^0 &\defeq \{ X\in \derx_{\A} : \pair[X][\derf_{\A}]\subset \alg_{\B}\}, \\
\derf_{\A}^0 &\defeq \{ \omega\in \derf_{\A} : \pair[\derx_{\A}][\omega]\subset \alg_{\B}\}.
\end{align*}
Geometrically, one can think of these two modules as the modules of vector fields and differential 1-forms that vanish at the variety. 

Furthermore, we let $\derx_{\B}$, $\derx_{\B}^{\top}$, $\derx_{\B}^0$, $\derf_{\B}$, $\derf_{\B}^{\bot}$ and $\derf_{\B}^0$ be the same modules as the ones with subscript $\A$, but now regarded as modules over the algebra $\alg_{\B}$. 
By this construction, it is easy to verify that $\Bt = \eBtn$, $\Bnn = \eBnn$ and $\Bo = \eBno$ are maximal ideal subpairs of $\Ao = \eApo$, $\A = \eApp$ and $\At = \eAtp$ respectively. Hence we obtain three non-degenerate quotient dual pairs over the algebra $\alg_{\c} = \alg_{\A}/\alg_{\B}$, which we denote by $\Coo = \eCoo = \Apo/\Btn$, $\c = \eCpp = \App/\Bnn$ and $\Ctt = \eCtt = \Atp/\Bno$. Here, the primary module $\derx_{\c}^{\bot}$\index{$\mathfrak{X}_{\mathscr{A}}^{\bot}$}\index{$\Omega_{\mathscr{A}}^{\top}$} equals $\derx_{\A}/\derx_{\B}^{\top}$, and similarly with the other modules that thus appear.

Now if two tangent vector fields $X, Y\in \derx_{\A}$ are such that $\d_X$ and $\d_Y$ both map the ideal $\ideal$ into itself, then evidently their commutator also maps this ideal into itself. Hence it follows that $\derx^{\top}_{\A}$ is a Lie-subalgebra of $\derx_{\A}$. So with the inherited pairing from $\A = \eApp$, we have that $\At = \eAtp$ becomes a (possibly degenerate) {\lrm}, and the natural injective morphism of dual pairs $\iota:\At\to\A$ yields an injective morphism of (possibly degenerate) {\lrm s}. Because $\derx_{\B}^0$ (which equals $\derx_{\A}^0$, except for the underlying algebra) is a Lie ideal of $\derx_{\A}^{\top}$ and by construction $\d \alg_{\B}\subset \derf_{\B}^{\bot}$, it follows from Theorem \ref{thm:quotient_lrm} that $\Ct = \eCtt$ is a {\lrm}, which we call a {\em quotient {\lrm}}\index{quotient Rinehart space}. 


\subsubsection{Riemannian Quotient Rinehart spaces}

\def\derxA{\derx_{\A}}
\def\derxAt{\derx_{\A}^{\top}}
\def\derfA{\derf_{\A}}
\def\derxB{\derx_{\B}}
\def\derfB{\derf_{\B}}

\begin{p}
Throughout most of this subsection, we will be concerned with the situation that the {\lrm} $\A$ in Table \ref{tab:cdi} is a musical Riemannian {\lrm}. In that case, it follows by Remark \ref{rem:sub_implies_inj} that the other dual pairs in the first two columns of Table \ref{tab:cdi} (i.e. all but the quotient dual pairs on the right) are turned into inner product pairs as well, whose inner products are canonically derived from the one on $\A$. Moreover, because $\B = \eBnn$ and $\Bo = \eBno$ are easily seen to be ideal inner product pairs of $\A = \eApp$ and $\At = \eAtp$ respectively, the quotient dual pairs $\c = \eCpp$ and $\Ct = \eCtt$ become inner product pairs as well, and we have that the connecting maps among these inner product pairs are all morphisms of inner product pairs.
\end{p}

\begin{p}\label{p:inner product derz}
The assumption that $\A$ has a musical Riemannian metric also implies that $\pair[X][\derf_{\A}] = \met[X][\derx_{\A}]$ for any $X\in \derx_{\A}$, and thus it follows that $\derx_{\A}^0$ can alternatively be defined by:
\begin{equation}\label{eq:characterization_derxAo}
\derx_{\A}^0 = \{ X\in \derx_{\A} : \met[X][\derx_{\A}]\subset \alg_{\B}\}.
\end{equation}
Moreover, from \ref{p:musical_implies_ideal_submodule_is_ideal_subpair} and Corollary \ref{cor:quotient_is_musical} it follows that the quotient {\inp} pair $\c$ is a musical pair as well. (This does not hold for $\Ctt$ in general, as can be seen from the case where $\A$ is a {\lrm} associated to a pseudo-Riemannian manifold, and the ideal $\alg_{\B}$ describes a light-like submanifold).
\end{p}


\def\Amod{\A^{\!\top}}
\def\tander{\derx_{\Amod}}
\def\phin{\phi^*}
\def\phis{\phi^\sharp}
\def\phif{\phi^\flat}
\def\dderx{\derx^\vee}
\def\kerderf{\derf_{0}}

\def\subm{\B^{\!\top}}
\def\extm{\B^{\!+}}
\def\derfo{\derf_{\ideal}^{\bot}}
\def\mis{\derx_{\B}^0}
\def\X{X}
\def\Y{Y}
\def\nablaA{\nabla^{\A}}
\def\nablaAt{\nabla^{\A^{\!\top}}}
\def\nablaC{\nabla^{\c}}
\def\nablaCt{\nabla^{\c^{\!\top}}}
\def\quotx{X_{\c}} \def\quoty{Y_{\c}} \def\quotz{Z_{\c}}
\def\quotf{f_{\c}}
\def\prex{X_{\A}} \def\prey{Y_{\A}} \def\prez{Z_{\A}}
\def\pref{f_{\A}}

\def\derz{\derx_{\A}^{0}}

In the following proposition, by an ordinary {\lrm} we mean a {\lrm} which is ordinary as a dual pair (cf. Definition \ref{def:ordinary_dual_pair}).   

\begin{proposition}\label{prop:connection derz}
In Table \ref{tab:cdi}, suppose that $\A$ has the structure of an ordinary {\musical} {\lrm} with Levi-Civita connection $\nablaA$. If $X\in \derx_{\A}^0$, or if $X\in \derx_{\A}^{\top}$ and $Y\in \derx_{\A}^{0}$, then we have that $\nablaA_X Y\in \derx_{\A}^{0}$. Moreover, a connection $\nablaC$ over $\Ctt$ is induced on $\derx_{\c}$, such that $\rhos$ and $\rhon$ constitute a morphism of affinely connected modules, i.e.
\begin{equation}\label{eq:quotient_preserves_connection}
\nablaC_{\quotx} \quoty = \rhos(\nablaA_{\prex} {\prey}),
\end{equation}
where $\prex$ and $\prey$ are preimages of $\quotx\in \derx_{\c}^{\top}$ and $\quoty\in\derx_{\c}$ under the quotient map $\rhos$. This connection is compatible with the inner product on $\derx_{\c}$.
\end{proposition}
\begin{proof}
First we consider the case that $X\in \derz$. It is clear that if $X = f X'\in \alg_{\B} \derx_{\A}$, then $\nablaA_X Y = f\nablaA_{X'} Y \in \alg_{\B} \derx_{\A} \subset \derx_{\A}^{0}$. Because we assumed that $\A$ is ordinary, we have that any $X\in \derz$ can be written as a linear combination of elements from $\alg_{\B} \derx_{\A}$, so that by linearity it follows that $\nablaA_X Y \in \derx_{\A}^{0}$ for all $Y\in \derx_{\A}$. 

Now we prove the case that $X\in \derx_{\A}^{\top}$ and $Y\in \derz$. By \eqref{eq:compatible}, we have that
$$\d_X\met[Y][Z] = \met[\nablaA_X Y][Z]+\met[Y][\nablaA_X Z]$$ 
for all $Z\in\derx_{\c}$. From \eqref{eq:characterization_derxAo} we know that $\met[Y][\nablaA_X Z]\in \ideal$, and also that $\met[Y][Z] \in \ideal$ so that $\d_Y\met[X][Z]\in \ideal$ as well. Hence $\met[\nablaA_Y X][Z]\in \ideal$ for all $Z\in \derx$, so that \eqref{eq:characterization_derxAo} now yields that $\nablaA_Y X\in \derz$. 

Because $\derx_{\A}^0 = \Ker \rhos$, we have that \eqref{eq:quotient_preserves_connection} is well-defined. Because the pairing, the differential operator and the inner product of $\A$ and $\c$ are all preserved under the quotient map $\rho:\A\to \c$, it follows from a straightforward verification that \eqref{eq:quotient_preserves_connection} defines a connection which is compatible with the metric.
\end{proof}

\begin{p}\label{p:projection_for_C}
Although we have natural morphisms (in fact exact sequences) 
$$\derx_{\c}^{\top} \hookrightarrow \derx_{\c} \twoheadrightarrow \derx_{\c}^{\bot}$$
and 
$$\derf_{\c}^{\bot}  \hookrightarrow \derf_{\c} \twoheadrightarrow \derf_{\c}^{\top},$$
because these morphisms go in opposite directions, there are a priori no relating morphism among the dual pairs $\Ct$, $\c$ and $\Co$. This changes however, in case we have a projection $(\cdot)^\top:\derx_{\A}\to \derx_{\A}^\top$ modulo $\alg_{\B}$. First, it gives us a surjective morphism $\A\to \At$, which is the identity on the cotangent space. 
Moreover, a unique projection $(\cdot)^{\top}:\derx_{\c}\to \derx_{\c}^{\top}$ is induced by
$$
(\rhos X)^{\top} \defeq \rhos(X^{\top}),
$$
which is well-defined, because this projection will map $\derx_{\B}^0$ onto itself (as follows from the last statement in Theorem \ref{thm:projection modulo B}). Together with the quotient map $\derf_{\A}/\derf_{\B}^0 \to \derf_{\A}/\derf_{\B}^{\bot}$, this induces a surjective morphism from $\c$ to $\Ct$.
\end{p}

\begin{p}
\def\Aop{\A^{\bot}}
\def\Apt{\A_{\top}}

Under the additional assumption that $\A$ is musical, there is actually more that can be said. First, let us define
$$
\derx_{\A}^{\bot} \defeq \{ X\in \derx_{\A} : \pair[X][\derx_{\A}^{\top}]\subset \alg_{\B}\}, 
$$
and 
$$
\derf_{\A}^{\top} \defeq \{ \omega\in \derf_{\A} : \pair[\derx_{\A}^{\bot}][\omega]\subset \alg_{\B}\}. 
$$
Now given the above projection modulo $\alg_{\B}$, we have that the quotient module $\derx_{\A}/\derx_{\B}^{\top}$ is isomorphic with $\derx_{\A}^{\bot}/\derx_{\B}^0$ through the map $$X + \derx_{\B}^{\top} \mapsto X^{\bot} + \derx_{\B}^0,$$
where $(\cdot)^{\bot}:\derx_{\A}\to \derx_{\A}^{\bot}$ stands for the orthogonal projection modulo $\alg_{\B}$, as defined by \eqref{eq:definition_ort_proj}. Likewise, $\derx_{\A}/\derx_{\B}^{\bot}$ can be identified with $\derx_{\A}^{\top}/\derx_{\B}^0$. Furthermore, the modules $\derx_{\A}$, $\derx_{\A}^{\top}$, $\derx_{\A}^{\bot}$ and $\derx_{\A}^{0}$ can be identified through musical isomorphism with $\derf_{\A}$, $\derf_{\A}^{\top}$, $\derf_{\A}^{\bot}$ and $\derf_{\A}^{0}$ respectively, and similarly for the spaces with subscript $\B$. This, on its turn, also induces isometries between $\derx_{\A}/\derx_{\B}^{\bot}$ and $\derf_{\A}/\derf_{\B}^{\bot}$, and between $\derx_{\A}^{\bot}/\derx_{\B}^{0}$ and $\derf_{\A}^{\bot}/\derf_{\B}^{0}$. We may summarize all this by stating that in case $\A$ is a musical Riemannian {\lrm},   for which there exists a projection $(\cdot)^{\top}:\derx_{\A} \to \derx_{\A}^{\top}$, then we have the situation as shown in Table \ref{tab:cdii}.

\begin{table}[!ht]
  \begin{center}
\begin{tikzcd}
\Bnn = \eBnn 
\arrow[r, hookrightarrow] 
& \Aoo = \eAoo 
\arrow[r, twoheadrightarrow, "\rho"] 
\arrow[d, hookrightarrow, shift right=1ex]
& \Coo = \eCoo 
\arrow[d, hookrightarrow, shift right=1ex]
\\
\Bnn = \eBnn 
\arrow[r, hookrightarrow] 
\arrow[d, equal] 
\arrow[u, equal]
& \A = \eApp 
\arrow[r, twoheadrightarrow, "\rho"]
\arrow[u, twoheadrightarrow, "\bot" swap, shift right=1ex]
\arrow[d, twoheadrightarrow, "\top", shift left=1ex]
& \c = \eCpp 
\arrow[u, twoheadrightarrow, "\bot" swap, shift right=1ex]
\arrow[d, twoheadrightarrow, "\top", shift left=1ex]\\
\Bnn = \eBnn 
\arrow[r, hookrightarrow] 
& \Att = \eAtt 
\arrow[u, hookrightarrow, shift left=1ex]
\arrow[r, twoheadrightarrow, "\rho"]
& \Ct = \eCtt 
\arrow[u, hookrightarrow, shift left=1ex]
\end{tikzcd}
  \end{center}
 \caption{}
 \label{tab:cdii}
\end{table}

All the arrows in this diagram represent morphisms of inner product pairs. Moreover, for the two injective morphisms in the middle column, the morphisms $(\cdot)^{\top}:\A\to \Att$ and $(\cdot)^{\bot}:\A\to \Aoo$ are retractions modulo $\Bnn$ (which are given by projections modulo $\alg_{\B}$ on each of their modules), and the injective morphisms on the right column have retractions $(\cdot)^{\top}:\c\to \Ctt$ and $(\cdot)^{\bot}:\c\to \Coo$ (given by projections on each of their modules).
\end{p}

\def\supcon{\nablaA} \def\subcon{\nabla^{\Ctt}} 
\def\supx{X} \def\supy{Y} \def\supz{Z}
\def\subx{\rhos \supy} \def\suby{\rhos \supy} \def\subz{\rhos \supz}

\begin{theorem}\label{thm:inherited connection}
In Table \ref{tab:cdii}, suppose that $\A$ has the structure of an ordinary musical {\lrm} with Levi-Civita connection $\nablaA$. 
Then a Levi-Civita connection $\subcon$ is induced on $\Ct$, which is uniquely determined by the equation:
\begin{equation}\label{eq:inherited connection}
\subcon_{\quotx} \quoty \defeq (\nablaC_{\quotx} \quoty)^{\top}.
\end{equation}
where $\quotx, \quoty\in \derx_{\c}^{\top}$.
\end{theorem}
\begin{proof}
We first prove that \eqref{eq:inherited connection} yields a covariant derivative, and then that it is really a Levi-Civita connection. Note that the projection on the right-hand side of \eqref{eq:inherited connection} stands for the projection $(\cdot)^{\top}:\derx_{\c}\to \derx_{\c}^{\top}$, which is uniquely determined by the projection $(\cdot)^\top:\derx_{\A}\to \derx_{\A}^{\top}$ (cf. Theorem \ref{thm:projection modulo B}). Here, $\quotx$ and $\quoty$ can be regarded as elements of $\derx_{\c}$ through the earlier mentioned inclusion $\derx_{\c}^{\top} \hookrightarrow \derx_{\c}$. 
Because the projection is $\alg_{\c}$-linear, it immediately follows that $\subcon_{\quotx} {\quoty}$ is $\alg_{\c}$-linear with respect to $\quotx$. In the following, let $\quotf\in \alg_{\c}$ and let $\prex, \prey\in \derx_{\A}^{\top}$ and $\pref\in \alg_{\A}$ be preimages of $\quotx$, $\quoty$ and $\quotf$ under the quotient morphism $\rho$. Then
\begin{align*}
\subcon_{\quotx} (\quotf \quoty) &= (\nablaC_{\quotx} (\quotf\quoty))^{\top}
= (\d_{\quotx} \quotf\ \quoty+ \quotf \nablaC_{\quotx} \quoty)^{\top}\\
&= \d_{\quotx} \quotf\ \quoty + \quotf \subcon_{\quotx} \quoty,
\end{align*}
so $\subcon_{\quotx}$ is indeed a derivation on $\derx_{\c}^{\top}$ over $\d_{\quotx}$.
The see that the connection is symmetric, first note that
\begin{align*}\label{eq:commutator_connection}
& \nablaC_{\quotx} \quoty - \nablaC_{\quoty}\quotx\\
&= \rhos(\nablaA_{\prex} \prey - \nablaA_{\prey}\prex)\\
&= \rhos [\prex,\prey]\\
&= [\quotx,\quoty] \in \derx_{\c}^{\top} \subset \derx_{\c}
\end{align*}
where we used that $\derx_{\A}^{\top}$ is closed under Lie brackets. Hence we have
\begin{align*}
& \subcon_{\quotx} \quoty - \subcon_{\quoty} \quotx\\
&= (\nablaC_{\quotx} \quoty - \nablaC_{\quoty}\quotx)^{\top}\\
&= ([\quotx,\quoty])^{\top} = [\quotx,\quoty].
\end{align*}
Furthermore we have
\begin{align*}
\d_{\quotx}\met[\quoty][\quotz] &= \rhon(\d_{\prex}\met[\prey][\prez])\\
&= \rhon(\met[\supcon_{\prex} \prey][\prez] + \met[\prey][\supcon_{\prex} \prez])\\
&= \met[\nablaC_{\quotx} \quoty][\quotz] + \met[\quoty][\nablaC_{\quotx} \quotz]\\
&= \met[(\nablaC_{\quotx} \quoty)^{\top}][\quotz] + \met[\quoty][(\nablaC_{\quotx} \quotz)^{\top}]\\
&= \met[\subcon_{\quotx} \quoty][\quotz] + \met[\quoty][\subcon_{\quotx} \quotz],
\end{align*}
showing that the connection is compatible with the Riemannian metric on $\Ct$, so that it is indeed a Levi-Civita connection. Note that the introduction of the projection is allowed because we were taking the inner product with an element from $\derx_{\c}^{\top}$ (a tangent vector field). Also note that there is no strict need to distinguish between $X_{\A}$ and $X_{\c}$ as elements of $\derx_{\A}^{\top}$ and $\derx_{\c}^{\top}$ or as elements of $\derx_{\A}$ and $\derx_{\c}$, but that there is a clear distinction between $\nablaC$ and $\subcon$.
\end{proof}

\begin{remark}
It follows from Theorem \ref{thm:fundamental theorem}, that at least when the induced Riemannian metric on $\Ct$ is non-degenerate and $2$ a unit in the ground ring $\ring$, then the above Levi-Civita connection induced by a projection $(\cdot)^\top:\derx_{\A}\to \derx_{\A}^{\top}$ modulo $\alg_{\B}$, is the only Levi-Civita connection that exists on $\Ct$.
\end{remark}

\begin{theorem}\label{thm:second_fun_form}
Under the conditions of the previous theorem, the so-called {\em second fundamental form}\index{second fundamental form} $h:\derx_{\c}^{\top}\otimes \derx_{\c}^{\top} \to \derx_{\c}^{\bot}$, which is given by the formula
\begin{equation}\label{eq:second_fun_form}
h(\quotx,\quoty) \defeq (\nablaC_{\quotx} \quoty)^{\bot}
\end{equation}
for $\quotx, \quoty\in \derx_{\c}^{\top}$, is a symmetric bilinear map.
\end{theorem}
\begin{proof}
First we show that $h$ is bilinear. As we saw in the previous theorem, we have for $\quotx, \quoty\in \derx_{\c}^{\top}$ that
$$\nablaC_{\quotx} \quoty - \nablaC_{\quoty} \quotx = [\quotx,\quoty] \in \derx_{\c}^{\top},$$
and hence
\begin{align*}
& h(\quotx,\quoty) - h(\quoty,\quotx)\\
&= (\subcon_{\quotx} \quoty - \subcon_{\quoty} \quotx)^{\bot}\\
&= ([\quotx,\quoty])^{\bot}\\
&= 0 + \derx_{\B}^{\bot} \in \derx_{\c}^{\bot}.
\end{align*}
To see that $h$ is bilinear, note that $h$ is clearly linear in its first argument, and because it is symmetric, it is also linear in its second argument. 
\end{proof}

We end this section with a lemma that we will use in Example \ref{ex:space_form}.

\def\pr{f}
\def\tander{\derx_{\A}^{\top}}
\def\derz{\derx_{\A}^0}
\begin{lemma}\label{lem: principal}
Let $\lra$ be a musical Riemannian {\lrm} and $\ideal$ a principal ideal of $\alg_{\A}$ generated by some non-constant $\pr\in\alg_{\A}$. Let us write $N$ for the gradient of $\pr$. Then the following holds:
\begin{enumerate}
\item\label{en:tander and metXN in ideal} $\derx_{\A}^{\top}$ consists exactly of those $X\in \derx$ such that $\met[X][N]\in \ideal$.
\item\label{en:projection} If there exists an element $q\in \alg$ such that $1-q\met[N][N]\in\ideal$, then we have a module endomorphism $(\cdot)^\top:\derx_{\A}\to \derx_{\A}^{\top}$, given by
$$X^\top \defeq X - q \met[X][N]N,$$
which is a projection modulo $\ideal$.
\item\label{en:subspace spanned by N} If $(\cdot)^\top:\derx_{\A}\to \derx_{\A}^{\top}$ is a projection modulo $\ideal$, and $X$ belongs to the submodule spanned by $N$, then $X^\top\in \derz$.
\end{enumerate} 
\end{lemma}
\begin{proof}

(1) By definition, $\derx_{\A}^{\top}$ consists exactly of those $X\in \derx_{\A}$ such that $\met[X][\Omega^{\bot}_{\A}]\subset \ideal$. Because any element in $\Omega^{\bot}_{\A}$ is a multiple of $\d \pr = \form{N}$, it thus follows that $\derx_{\A}^{\top}$ consists of those $X$ such that $\met[X][N]\in \ideal$.

(2) For $Y\in\derx_{\A}$ we have
\begin{align*}
\met[Y^\top][N] &= \met[Y - q\met[Y][N]N][N]\\
&= (1-q\met[N][N])\met[Y][N]\in\ideal,
\end{align*}
so $Y^\top\in\tander$. Furthermore, because $Y-Y^\top$ is a scalar multiple of $N$, it follows from \eqref{en:tander and metXN in ideal} that $\met[X][Y-Y^\top]\in \ideal$ for all $X\in \tander$, so we have indeed a projection modulo $\ideal$.

(3) For any $Y\in \derx$ we have that
$$\met[N^\top][Y] = \met[N^\top][Y^\top] = \met[N][Y^\top]\in \ideal,$$
so that from \eqref{eq:characterization_derxAo} it follows that $N^{\top}\in \derz$. Of course, any multiple of $N$ will then also be projected to $\derz$.


\end{proof}


\subsection{Euclidean Rinehart spaces}

For a general {\lrm} $A$, there does not always exist a basis for the module $\derx_{\A}$. And even if it exists, such a basis may or may not satisfy certain convenient properties. In the following definition, we will provide some names for a couple of situations that may arise. 

{
\def\myalg{\alg}
\def\myderx{\derx}
\def\myderf{\derf}
\begin{definition}
Given a {\lrm}, 
if $\myalg$ is a unital algebra, and the tangent space $\myderx$ and cotangent space $\derf$ have bases which are dual to each other, i.e. there are bases $\{\corv_i\}$ and $\{\omega_i\}$ of tangent and cotangent elements respectively, such that $\pair[\corv_i][\omega_j]$ equals $1$ if $i=j$ and $0$ otherwise, then we call our {\lrm} a {\em locally affine Rinehart space}\index{locally affine Rinehart space}. By an \emindex{affine Rinehart space}, we mean a locally affine {\lrm} for which the dual bases can be chosen in such a way, that the differential 1-forms $\{\omega_i\}$ are all exact, i.e. $\omega_i = \d\cor_i$ for some $\cor_i\in\myalg$. In this case, we refer to $\cor_i$ as the {\em coordinates}\index{coordinate} on $\lra$, and we refer to the dual elements $\corv_i$ of $\d\cor_i$ as the {\em coordinate vector fields}\index{coordinate!vector field}.
\end{definition}

\begin{proposition}
Any locally affine {\lrm} $(\derx_{\lra},\derf_{\lra})$ is isomorphic with $(\derx_{\lra},\dual \derx_{\lra})$.
\end{proposition}
\begin{proof}
This follows immediately from the observation that to any $\eta\in\Hom(\derx,\alg)$ one can associate the 1-form:
$$
\eta = \eta(\corv_1)\,\omega_1 + \ldots + \eta(\corv_n)\,\omega_n.
$$
\end{proof}


\begin{proposition}\label{prop:vanishing lie bracket}
If $\lra$ is an affine {\lrm} with coordinate vector fields $\corv_i$, then $[\corv_i,\corv_j]=0$ for any indices $i$ and $j$.
\end{proposition}
\begin{proof}
It is easily seen that $\d_{[\corv_i,\corv_j]}\cor_k = \d_{\corv_i} \d_{\corv_j} \cor_k - \d_{\corv_j} \d_{\corv_i} \cor_k = 0$ for any indices $i$,$j$ and $k$. Because by assumption the differentials $\d\cor_i$ span all of $\derf$, we can write any differential 1-form $\omega$ as a finite linear combination of differential 1-forms $\d\cor_i$. As a consequence, we have $\pair[{\lbrack \corv_i,\corv_j\rbrack}][\omega] = 0$ for all $\omega\in\derf$, and so by the assumed non-degeneracy of the pairing between $\derx$ and $\derf$, it follows that $[\corv_i,\corv_j] = 0$ for any $i$ and $j$.
\end{proof}

If we are given a locally affine {\lrm} with basis $\{\corv_i\}$, we can put a Riemannian metric on it by setting 
\begin{equation}\label{eq:Euclidean metric}
\met[\corv_i][\corv_j] = \dird{i}{j}
\end{equation}
where $\dird{i}{j}$ is the Kronecker delta function, and extend this linearly over $\alg$. This inner product is clearly non-degenerate, and it is also {\gen}, for if $\omega\in \derf$ is such that $\pair[\corv_i][\omega] = f_i$, we may set $Z = f_1 \corv_1 +\ldots +f_n \corv_n$, where $n$ is the dimension of the tangent (and cotangent) space, so that $\pair[Y][\omega] = \met[Y][Z]$ for all $Y\in \derx$. 

\begin{definition}\label{def:Euclidean metric}
Given a locally affine {\lrm} with basis $\{\corv_i\}$
, we call the metric that satisfies \eqref{eq:Euclidean metric} the \emindex{Euclidean metric} associated to the basis $\{\corv_i\}$. A locally affine {\lrm} with a Euclidean metric is called a {\em locally Euclidean {\lrm}}\index{locally!Euclidean Rinehart space}. An affine {\lrm} with a Euclidean metric, is called a {\em Euclidean {\lrm}}\index{Euclidean Rinehart space}. 
\end{definition}

\begin{proposition}\label{prop:locally_Euclidean_is_musical}
The inner product on a locally Euclidean {\lrm} is musical, where the musical isomorphism is given by
\begin{equation}\label{eq:Euclidean_musical_isomorphism}
\form{(f_1 \corv_1 +\ldots +f_n \corv_n)} = f_1 \omega_1 +\ldots +f_n \omega_n,
\end{equation}
the differential is given by
\begin{equation}\label{eq:Euclidean_differential}
\d f = \d_{\corv_1} f\ \omega_1 +\ldots +\d_{\corv_n} f\ \omega_n,
\end{equation}
and the gradient is given by
\begin{equation}\label{eq:gradient}
\grad f = \d_{\corv_1} f\ \corv_1 +\ldots +\d_{\corv_n} f\ \corv_n.
\end{equation}
\end{proposition}
\begin{proof}
It is evident that \eqref{eq:Euclidean_musical_isomorphism} determines an isomorphism, for which holds that $\pair[Y][\form{Z}]_{\A} = \pair[Z][\form{Y}]_{\A} = f_1 g_1 + \ldots + f_n g_n$ if $Y = f_1 \corv_1 +\ldots +f_n \corv_n$ and $Z = g_1 \corv_1 +\ldots +g_n \corv_n$. Hence it follows by the converse statement of Proposition \ref{prop:musical isomorphism} that the inner product is musical, and related to the above isomorphism through Equation \eqref{eq:relation_inp_and_musical_isom}. Equation \eqref{eq:Euclidean_differential} is verified by the observation that 
$$
\pair[\corv_j][\d_{\corv_1} f\ \omega_1 +\ldots +\d_{\corv_n} f\ \omega_n] = 
\d_{\corv_j} f = \pair[\corv_j][\d f]
$$
for any basis element $\corv_j$. Equation \eqref{eq:gradient} follows as a direct consequence of this.
\end{proof}

\begin{remark}
In the situation that our {\lrm} is Euclidean instead of locally Euclidean, it follows from Equation \eqref{eq:gradient} that 
\begin{equation}\label{eq:Euclidean_gradient_x_i}
\grad {\cor_i} = \corv_i.
\end{equation}
\end{remark}

\begin{theorem}\label{thm:locally Euclidean has LC connection}
On a locally Euclidean {\lrm} with Euclidean basis $\{X_i\}$, there exists a unique Levi-Civita connection satisfying
\begin{equation}\label{eq:Euclidean_connection}
\nabla_Y (f \corv_i) = \d_Y f\ \corv_i
\end{equation}
for all $f\in \alg_{\lra}$, $Y\in \derx_{\lra}$ and basis elements $X_i$. Moreover, this connection is flat.
\end{theorem}
\begin{proof}
First note that for a connection, satisfying \eqref{eq:Euclidean_connection} is equivalent to the requirement that 
\begin{equation}\label{eq:Euclidean_connection_equiv}
\nabla_Y \corv_i = 0
\end{equation}
for all $Y\in \derx_{\lra}$ and basis elements $X_i$. Now let us define a connection simply by linear extension of Equation \eqref{eq:Euclidean_connection}. Because $\{\corv_i\}$ is a basis for $\derx$, this determines the value of $\nabla_Y Z$ for any $Y,Z\in \derx$. Furthermore we have that
\begin{align*}
\nabla_{f \corv_i} (g \corv_j) - \nabla_{g \corv_j} (f \corv_i) &= f \derive{\corv_i} g\ \corv_j - g \derive{\corv_j} f\ \corv_i\\
&= [f \corv_i,g \corv_j],
\end{align*}
by Proposition \ref{prop:vanishing lie bracket} and \eqref{eq:lie bracket two functions}. Because this equation extends additively to the whole of $\derx$, we obtain that the connection is symmetric.
Also we have that
\begin{align*}
\derive{Y} \met[f \corv_i][g \corv_j] &= \derive{Y}(fg)\met[\corv_i][\corv_j] = (g \derive{Y}f +f \derive{Y}g)\met[\corv_i][\corv_j]\\
&= \met[\derive{Y}f\ \corv_i][g \corv_j] + \met[f \corv_i][\derive{Y}g\ \corv_j]\\
&= \met[\nabla_Y (f \corv_i)][g \corv_j] + \met[f \corv_i][\nabla_Y (g \corv_j)].
\end{align*}
Like before, this equality extends additively to the whole $\derx$, so that we may conclude that $\nabla$ is indeed a Levi-Civita connection. To prove that the Riemann tensor is identically zero, it suffices to show that $R(\corv_i,\corv_j)\corv_k = 0$ for any indices $i$,$j$ and $k$. But this is obviously the case, because we have that $\nabla_{\corv_i}\corv_j = 0$ for any indices $i$ and $j$.
\end{proof}

\begin{proposition}
A locally affine {\lrm} is an ordinary {\lrm}, and a locally Euclidean {\lrm} is an ordinary Riemannian {\lrm}.
\end{proposition}
\begin{proof}
Suppose that $\pair[X][\omega]\in \alg_{\B}$ for all $\omega\in \derf_{\A}$. Then writing $X = f_1 \corv_1 +\ldots +f_n \corv_n$, we have that $\pair[X][\omega_i] = f_i \in \alg_{\B}$ for $i = 1,\ldots,n$. Hence $X$ is clearly spanned by $\alg_{\B}\derx_{\A}$. Likewise one proves that if $\pair[X][\omega]\in \alg_{\B}$ for all $X\in \derx_{\A}$, then $\omega$ is spanned by $\alg_{\B}\derf_{\A}$. The statement about a locally Euclidean {\lrm} follows immediately from Proposition \ref{prop:locally_Euclidean_is_musical} and Remark \ref{rem:musical_implies_same_max_sub}.
\end{proof}

\subsubsection{Examples}\label{part:final examples}

\begin{example}\label{ex:real euclidean}
In Example \ref{ex:real manifold}, the manifold $M=\R^n$ gives rise to an affine {\lrm}, where the coordinates are given by standard coordinates $x_i$ on $\R^n$, and coordinate vector fields by $X_i = \frac{\partial}{\partial x_i}$. Endowing this affine {\lrm} with a Euclidean metric,
we obtain by Theorem \ref{thm:locally Euclidean has LC connection} a flat Levi-Civita connection, which by Theorem \ref{thm:fundamental theorem} must be unique. This connection corresponds with the standard derivation of vector fields on $\R^n$.
\end{example}

\begin{example}\label{ex:complex euclidean}
As discussed in Example \ref{ex:complex manifold}, dependent on whether we regard a complex manifold from a holomorphic viewpoint or a real $2n$-dimensional viewpoint, we end up with a different {\lrm} (although these can be related through an injective morphism from the first into the second). In either case, the manifold $\C^n$ gives rise to an affine {\lrm}. Dependent on our choice, the coordinates are given by the $n$ standard holomorphic coordinates $z_i$ on $\C^n$, or by $2n$ coordinates, which are typically $\{x_i,y_i\}_{1\le i\le n}$, where $x_i = \mathrm{Re}(z_i)$ and $y_i = \mathrm{Im}(z_i)$. Another option for the latter case is to take as coordinates $\{z_i, \bar z_i\}_{1\le i\le n}$, where $z_i$ are the standard holomorphic coordinates and $\bar z_i$ the complex conjugates of $z_i$.

In the first case, the associated Euclidean {\lrm} is given by the {\em complex Euclidean $n$-space}, and its Euclidean metric is the complex bilinear form given by:
$\met[\del{z_i}][\del{z_j}] = 1$ if $i=j$ and $0$ otherwise. A complex manifold with such a complex bilinear metric is called a {\em holomorphic Riemannian manifold}. As will be discussed thoroughly in the second part of this thesis, such manifolds can be used to relate the geometries of so-called Wick-related (pseudo-)Riemannian manifolds.  

In the case that we consider $\C^n$ as real $2n$-dimensional manifold, then Definition \ref{def:Euclidean metric} requires that the metric is given by $\met[\del{x_i}][\del{y_j}] = 0$ and $\met[\del{x_i}][\del{x_j}] = \met[\del{y_i}][\del{y_j}] = 1$ if $i = j$ and $0$ otherwise. This defines the Riemannian metric that corresponds to an {\em Hermitian form}. Through complex extension of this Riemannian metric, we end up with an {\inp} that is indeed bilinear over the function algebra of smooth complex-valued functions.
\end{example}

\begin{example}\label{ex:finite field}
In Example \ref{ex:polynomial ring}, the affine {\lrm} is obtained by setting $\ring = \F_q$, $\alg = \F_q[x_1,\ldots,x_n]$, $\derx = \Der(\alg)$, which is generated by a basis $\{\del{x_i}\}$ of formal derivations on the polynomial ring, and $\derf = \dual \derx$, for which we have a basis of dual elements $\{\omega_i\}$ such that $\d x_i = \omega_i$. By Theorem \ref{thm:locally Euclidean has LC connection}, we can endow this {\lrm} with a Euclidean metric together with a flat Levi-Civita connection associated to it. Note that this is possible even if $q$ is a power of $2$, in which case $2$ is not a unit in $\ring$.
\end{example}

\subsection{Non-flat space forms}

We will now apply the theory we have developed so far to show the following: suppose $\lra$ is a Euclidean {\lrm} with coordinates $\{\cor_i\}$, for which the ground ring $\ring$ contains $2$ as a unit. Then for any unit $c\in\ring\subset \alg$ such that $\cor_1^2+\ldots+\cor_n^2-c^{-1}$ is not invertible, there exists a quotient Rinehart space of $\lra$ which together with its uniquely induced Levi-Civita connection has constant sectional curvature $c$, i.e. it is a non-flat space form.

\def\ci{c\,} \def\cinv{c^{-1}}
\def\u{x}
\def\U{X}
\def\su{y}
\def\sU{Y}

\begin{example}\label{ex:space_form}
Let $\lra$ be a Euclidean {\lrm} of dimension $n$, with coordinates given by $\{\cor_i\}$, and suppose further that $2$ is a unit in $\ring$. Now let $c$ be a unit in $\ring$, regarded as constant element of $\alg$, such that $(\cor_1^2+\ldots +\cor_n^2 - \cinv)$ determines a proper principal ideal. For the sake of convenience, let us consider ${\pr} = \frac{1}{2}(\cor_1^2+\ldots +\cor_n^2 - \cinv)$ as generator of $\ideal$, and by \eqref{eq:gradient} we calculate that $N = \grad {\pr}$ is given by $\cor_1 \corv_1 + \ldots + \cor_n \corv_n$. We may verify the following facts:
\begin{align*}
\d_N \cor_i &= \cor_i\\
\nabla_X N &= X \quad \mbox{for all}\ \ X\in \derx\\
\met[\corv_i][N] &=\cor_i\\
\met[N][N] &= \cor_1^2+\ldots +\cor_n^2.
\end{align*}
So our ideal $\ideal$ is alternatively given by $(\met[N][N]_{\A}-\cinv)$. By letting $q = \ci$, we find that $1-q\met[N][N]_{\A} = \ci (\cinv-\met[N][N]_{\A})\in\ideal$. Hence by Lemma \ref{lem: principal} \eqref{en:projection} we have a map $(\cdot)^\top:\derxA\to \derxAt$ which is a projection modulo $\ideal$. In this case, it is given by
$$\ort X = X - \ci\met[X][N]_{\A} N.$$
In particular we have 
$$\ort \U_i = \U_i - \ci \u_i N.$$
We calculate that
\begin{align}
\meti[\ort \U_i][\ort \U_j] &= \meti[\U_i- \ci \u_i N][\U_j-\ci \u_j N]\notag\\
&= \dird{i}{j} - 2 \ci \u_i \u_j + \ci^2 \u_i \u_j \meti[N][N]\label{eq:inpAm}
\end{align}
where $\dird{i}{j}$ denotes the Kronecker delta function. Further we have
\begin{equation}\label{eq:derAm}
\d_{\ort \U_i} \u_j = \d_{\U_i - \ci \u_i N} \u_j = \dird{i}{j} - \ci \u_i \u_j
\end{equation}
and 
\begin{align}\label{eq:conAm}
\notag\nabla_{\ort \corv_i} \ort \corv_j &= \nabla_{\ort \corv_i} (\U_j - \ci \u_j N)
= -c\,\nabla_{\ort \corv_i}(\u_j N)\\ &= - \ci \u_j \ort \corv_i - \ci(\dird{i}{j} - \ci \u_i \u_j) N
\end{align}
for all $i$ and $j$.

\def\subcon{\nabla^{\Ct}}
\renewcommandtwoopt{\metii}[2][\cdot][\cdot]{\langle #1,#2\rangle_{\c}}
\def\derI{\derx_{\c}^{\top}}
We will now show that the {\lrm} $\Ct = \eCtt$ associated to the ideal $\ideal$ has constant sectional curvature. By \ref{p:projection_for_C}, we know that the projection $(\cdot)^\top:\derxA\to \derxAt$ modulo $\ideal$ induces a surjective morphism  $\c\to \Ct$. By Theorem \ref{thm:inherited connection}, there exists a Levi-Civita connection $\subcon$ on $\Ct$, determined by \eqref{eq:inherited connection}. 
Let us now write 
$\sU_i$ for $(\rhos \U_i)^{\top} = \rhos(\U_i^{\top})$, and $\su_i$ for $\rhon(\u_i)$. 
Because $\{\U_i\}$ forms a basis of $\derx_{\A}$, we know that the derivations $\{\sU_i\}$ span the whole space $\derx_{\c}^{\top}$. Equations \eqref{eq:inpAm} up to \eqref{eq:conAm} now give us:
\begin{align}
\label{eq:inp}\metii[\sU_i][\sU_j] &= \dird{i}{j} - \ci\su_i \su_j\\
\label{eq:der}\d_{\sU_i} \su_j &= \dird{i}{j} - \ci\su_i \su_j\\
\label{eq:con}\subcon_{\sU_i} \sU_j &= - \ci \su_j \sU_i
\end{align}
for all $i$ and $j$. Note that in the first equation, we use that $\rhon\meti[N][N] = \cinv \in \alg_{\c}$. By the symmetry of the induced Levi-Civita connection, we moreover have that
\begin{equation}\label{eq:lie}
[\sU_i,\sU_j] = \ci(\su_i \sU_j - \su_j \sU_i).
\end{equation}
Now we calculate that
\begin{align*}
\cinv\subcon_{\sU_i} \subcon_{\sU_j} \sU_k &= \subcon_{\sU_i} (- \su_k \sU_j) = (-\dird{i}{k} + \ci \su_i \su_k) \sU_j + \ci \su_j \su_k \sU_i\\
\cinv\subcon_{\sU_j} \subcon_{\sU_i} \sU_k &= \subcon_{\sU_j} (- \su_k \sU_i) = (-\dird{j}{k} + \ci \su_j \su_k) \sU_i + \ci \su_i \su_k \sU_j\\
\cinv\subcon_{[\sU_i,\sU_j]} \sU_k &= \subcon_{\su_i \sU_j - \su_j \sU_i} \sU_k = -\ci \su_i \su_k \sU_j + \ci \su_j \su_k \sU_i.
\end{align*}
From this, we deduce that
\begin{align}
\cinv R(\sU_i,\sU_j)\sU_k &= \cinv \subcon_{\sU_i} \subcon_{\sU_j} \sU_k -\cinv \subcon_{\sU_j} \subcon_{\sU_i} \sU_k -\cinv \subcon_{[\sU_i,\sU_j]} \sU_k\\
&= (\dird{j}{k} - \ci\su_j \su_k) \sU_i - (\dird{i}{k} - \ci\su_i \su_k) \sU_j\\
&= \metii[\sU_j][\sU_k] \sU_i - \metii[\sU_i][\sU_k] \sU_j.
\end{align}
and so
$$R(\sU_i,\sU_j)\sU_k = \ci\,(\metii[\sU_j][\sU_k] \sU_i - \metii[\sU_i][\sU_k] \sU_j).$$
Because $R$ is a tensor, and the $\sU_i$ span the whole space $\derI$, we have that this equality extends linearly to any $X,Y,Z\in \derI$, which gives us 
$$R(X,Y)Z = \ci\,(\metii[Y][Z] X - \metii[X][Z] Y).$$
This shows that $\Ctt$ is a space form whose constant sectional curvature equals $c$.
\end{example}

\begin{remark}
As a final remark, let us consider what the previous result implies for a couple of instances of the {\lrm} $\lra$. If $\lra$ is the {\lrm} corresponding to $\R^n$ with its ring of smooth functions (cf. Example \ref{ex:real euclidean}), then $(x_1^2+\ldots+x_n^2-\cinv)$ will only determine a proper ideal if $c>0$, in which case it corresponds with the real $n-1$-dimensional sphere. In case $\lra$ corresponds with the complex Euclidean $n$-space (cf. Example \ref{ex:complex euclidean}) then we may obtain a so-called {\em holomorphic Riemannian sphere} (cf. \cite{PeVa}) for any $c\ne 0$. Something unusual happens when we choose for our function algebra a polynomial ring like $\alg = \R[x_1,\ldots,x_n]$ (cf. Example \ref{ex:polynomial ring}). In that case $(x_1^2+\ldots+x_n^2-\cinv)$ determines a proper ideal for both positive and negative values of $c$, and hence we also obtain a quotient Rinehart space with a negative constant curvature. Here we see a phenomenon that one also encounters in real algebraic geometry, namely that a non-unital real polynomial with no real solutions does determine a proper ideal in the polynomial ring, and hence a valid object of study, even though its corresponding variety is just the empty set.

It is also worth noting that the different holomorphic Riemannian spheres, whose sectional curvature ranges over all non-zero complex values, can all be related to each other through a dilation on the ambient space (as will be discussed more detailed in the next chapter). In the setting of {\lrm s}, such a dilation gives rise to an automorphism $\phi$ on the ambient {\lrm}, for which $\phi^{\sharp}$ and $\phi^{\flat}$ are the identity, and $\phi^* x_i = \alpha x_i$ for all coordinates $x_i$ and some fixed non-zero $\alpha\in \ring$. Moreover, on constant functions is $\phi^*$ the identity. Because the function space, generally speaking, does not only contain polynomials, some extra requirements would actually be needed to make this determine the unique algebra homomorphism $\phi^*$ that corresponds geometrically with a dilation on the ambient space, but for the current discussion it is not relevant to make this mathematically precise.

In the real case, the (virtual) spheres of positive and negative sectional curvature cannot be related to each other through a dilation, as there is no real value that squares to a negative number. Something similar happens when the ambient space is a Euclidean {\lrm} over a finite field, as described in Example \ref{ex:finite field}. In that case we get a family of spheres whose constant sectional curvature is a square in the finite field, and an equally large class of spheres whose constant sectional curvature is not a square. Only two such spheres within the same class can be related through a dilation.
\end{remark}
}

}

\cleardoublepage

%

%
  \graphicspath{{\chaptersdir/chapter2/image/}}%
  \chapter{\titletwo}\label{ch:chapter2}

The structure of this chapter, based on joint work with Joeri Van der Veken \cite{PeVa}, is as follows. In the Sections \ref{SC Preliminaries}, \ref{SC Complexifications} and \ref{SC Hol Riem submanifold theory}, we rephrase and extend the theory of holomorphic Riemannian manifolds from respectively a {\em complex-linear} perspective, a complex analytic perspective and subsequently the viewpoint of submanifold theory. This theory is then used in Section \ref{SC applications}, the final section of this chapter, where we demonstrate by three different examples our method to relate certain kinds of submanifolds in one (pseudo-)Riemannian space, to submanifolds with corresponding geometric properties in other so-called Wick-related spaces.

{

\newcommand{\Id}{\operatorname{Id}}
\newcommand{\R}{\mathbb{R}}
\newcommand{\C}{\mathbb{C}}
\newcommand{\Z}{\mathbb{Z}}
\newcommand{\hm}{g}
\newcommand{\nm}{g'}
\newcommand{\tm}{g''}
\newcommand{\del}{\partial}
\newcommand{\holsec}{H^0}
\newcommand{\smsec}{C^\infty}
\newcommand{\rasec}{C^\omega}
\newcommand{\seqref}[1]{{\it (\ref{#1})}}
\newcommand{\remarksymbol}{}
\newcommandtwoopt{\met}[2][\cdot][\cdot]{\langle #1,#2\rangle}
\newcommandtwoopt{\meti}[2][\cdot][\cdot]{\langle #1,#2\rangle'}
\newcommandtwoopt{\metii}[2][\cdot][\cdot]{\langle #1,#2\rangle''}
\newcommand{\emdex}[1]{{\em #1}\index{#1}}
\renewcommand{\d}{\mathrm{d}}

\section{Preliminaries on complex Riemannian geometry}\label{SC Preliminaries}

\subsection{Complex vector spaces with holomorphic inner product}

In the following, we will assume all vector spaces to be finite dimensional. We adopt the following two definitions on subspaces of a complex vector space (cf. \cite{Bo}).

\begin{definition}
Let $V$ be a complex vector space. We call a real linear subspace $W\subset V$ {\em totally real}\index{totally real!subspace} if $W\cap i\,W = \{0\}$ and {\em generic}\index{generic!subspace} if $W+ i\,W = V$.
\end{definition}

Observe that if $W$ is both generic and totally real, then its real dimension equals the complex dimension of $V$. 

\begin{definition}
By a {\em holomorphic inner product}\index{holomorphic!inner product} on a complex vector space $V$, we mean a complex bilinear form $\met$ on $V$. A {\em holomorphic inner product space}\index{holomorphic!inner product space} is simply a complex vector space with a holomorphic inner product. Unlike in the previous chapter, we will in this chapter assume inner products to be non-degenerate by default.
\end{definition}

Given a holomorphic inner product on $V$ one can always choose an orthonormal basis. This means that any $n$-dimensional holomorphic inner product space can be identified with $\C^n$, equipped with the standard holomorphic inner product
\begin{equation}\label{EQ standard inner product}
\met[X][Y] = \ X_1 Y_1+\ldots+X_n Y_n,
\end{equation}
where $X=(X_1,\ldots,X_n)\in\C^n$ and $Y=(Y_1,\ldots,Y_n)\in \C^n$. The name `holomorphic inner product' comes from the fact that, unlike the often seen sesquilinear inner product, the above inner product is a holomorphic function from $\C^n\times \C^n$ to $\C$. 

We note that, given a complex-linear operator $A:V \to V$, there does not necessarily exist a basis of $V$ consisting of eigenvectors of $A$, not even when $A$ is symmetric with respect to a holomorphic inner product on $V$. However, if there exists a basis of eigenvectors of a symmetric operator $A$, there also exists an orthonormal basis of eigenvectors of $A$.

\begin{definition}\label{DF real slice hol in prod space}
Given a holomorphic inner product space $(V,\met)$, we use the term {\em real slice}\index{real slice!of inner product space} to denote a real linear subspace $W\subset V$, for which $\met$ restricted to $W$ is non-degenerate and real valued, i.e., $\met[X][Y]\in \R$ for all $X,Y\in W$.
\end{definition}

\begin{example}\label{EX pseudo-Euclidean}
  Consider the standard holomorphic inner product space $(\C^n,\met)$ and let $\{e_1,\ldots,e_n\}$ be the standard basis. The real linear subspace
  \begin{equation}\label{EQ pseudo-Euclidean} 
  \R^n_k := \mathrm{span}_{\R}\{ie_1,\ldots,ie_k,e_{k+1},\ldots,e_n\}  
  \end{equation}
  is a real slice. Indeed, the restriction of the inner product $\met$ to $\R^n_k$ is real valued. In particular, it is the standard pseudo-Euclidean metric of signature $k$ on $\R^n$, as the notation $\R^n_k$ suggests.
  If $z_j = x_j + iy_j$ are the standard coordinates on $\C^n$, the space $\R^n_k$ is given by $x_1 = \ldots = x_k = y_{k+1} = \ldots = y_n = 0$. Hence, $(y_1,\ldots,y_k,x_{k+1},\ldots,x_n)$ are natural real coordinates on $\R^n_k$.  
\end{example}
  
  Since any $n$-dimensional holomorphic inner product space $(V,\met)$ can be identified with the standard holomorphic inner product space $(\C^n,\met)$ by choosing an orthonormal basis, we can always find a real slice in $V$ for any signature $k$, with $0 \leq k \leq n$ , namely the one corresponding to $\R_k^n \subset \C^n$. One easily verifies that all real slices of $V$ come from such an identification, in particular, they are all related by the action of the complex orthogonal group $\mathrm{O}(n,\C)$ on $V$.

\begin{proposition}\label{ST real slices hol in prod space}
The real slices of a holomorphic inner product space are totally real subspaces.
\end{proposition}
\begin{proof}
Let $W$ be a real slice of a holomorphic inner product space $(V,\met)$, and let $X\in W\cap i\,W$. We need to show that $X=0$. Since $X\in iW$, there exists a vector $X'\in W$ such that $X=i X'$. Then we have that for all $Y\in W$, both $\met[X][Y]$ and $\met[X'][Y]=-i\met[X][Y]$ are real. Hence $\met[X][Y]=0$ for all $Y\in W$, so that the non-degeneracy of $\met$ implies $X=0$.
\end{proof}

\subsection{Complex-linear differential geometry}

In the following, we are concerned with complex manifolds. 
Importantly, we will regard complex manifolds from the perspective of what we call {\em complex-linear differential geometry}, which may be characterized by the property that all linear algebra involved is assumed to have $\C$ rather than $\R$ as its underlying field. This means, for instance, that a vector basis is usually assumed to be linearly independent over $\C$, rather than over $\R$, and tensors over the holomorphic tangent bundle are assumed to have $\C$ as their underlying field.

Indeed, the tensors (of any mixed type) we will be considering, all turn out to be tensorial over $\C$ in each of their components. In particular, we have that those tensors are all {\em pure} in the sense of \cite{YaAk}. Note that such assumptions cannot be imposed in the context of Hermitian manifolds without additionally complexifying them, because an Hermitian metric is not even complex bilinear with respect to the intrinsic complex structure of the manifold.

Apart from being $\C$-linear, the tensor fields we will be concerned with are usually also holomorphic. We recall the following definition of holomorphic vector fields in a general holomorphic vector bundle:

\begin{definition}
Let $\pi:V\to M$ be a holomorphic vector bundle. A section $X$ of $V$ is called {\em holomorphic}\index{holomorphic!section of vector bundle}, if $X$ is holomorphic as a map between the complex manifolds $M$ and $V$. The space of holomorphic sections of $V$ is denoted by $\holsec(M,V)$\index{$H^0$}.
\end{definition}

It follows directly that a vector field is holomorphic if and only if it has holomorphic component functions with respect to any local complex coordinates. For the special case of $V = T_l^k M$\index{$T_l^k M$}, by which we mean the tensor bundle of type $(l,k)$ over the holomorphic tangent bundle $TM$, this means that a smooth section $F$ of $T_l^k M$ is holomorphic if and only if the component functions $F_{i_1\ldots i_l}^{j_1\ldots j_k}$ are holomorphic relative to any local holomorphic coordinates $\{z_1,\ldots,z_n\}$ of $M$.

The space of holomorphic tensors fields is closed under the usual tensor operations such as addition, tensor multiplication and contraction.

\subsection{Complex and holomorphic Riemannian geometry}\label{SC complex Riemannian geometry}

In this section we will review some basic theory on complex Riemannian manifolds, which we will rephrase here from the viewpoint of the complex-linear setting.

We will define a {\em complex Riemannian manifold}\index{complex Riemannian manifold} as a complex manifold $M$ endowed with a {\em complex Riemannian metric}\index{complex Riemannian metric} $\met$: a smooth section of $T^0_2M$, which is non-degenerate and symmetric (cf. \cite{Le}). Hence, $\met$ is a complex bilinear form on every tangent space.

Like in the case of ordinary (pseudo-)Riemannian manifolds, any complex Riemannian manifold can be equipped with a {\em Levi-Civita connection}\index{Levi-Civita!connection} $\nabla$, which is the unique affine connection that is symmetric and compatible with the complex Riemannian metric, i.e.,
\begin{align}
[X,Y] &= \nabla_X Y - \nabla_Y X,\\
X \met[Y][Z] &= \met[\nabla_X Y][Z]+\met[Y][\nabla_X Z].
\end{align}
for all vector fields $X$, $Y$ and $Z$. Note that in this chapter we use the notation $X f$ to denote the derivation of $f$ along the vector field $X$, which more formally is written as $\mathrm{d}_X f$. In the following, we will focus our attention to complex Riemannian manifolds for which the complex metric is holomorphic, which are called {\em holomorphic Riemannian manifolds}\index{holomorphic Riemannian metric} (cf. \cite{DuZe}). 

\begin{remark}
As has been pointed out in \cite{BoFrVo}, there is a direct correspondence between $n$-dimensional holomorphic Riemannian manifolds and anti-K\"ahler manifolds (also called K\"ahler-Norden manifolds): by restricting the holomorphic metric $\met$ to its real part $\meti$, our manifold may be regarded as a real $2n$-dimensional pseudo-Riemannian manifold with canonical complex structure $J$ and a Levi-Civita connection $\nabla$ induced by $\meti$, satisfying:
$$\meti[JX][JY]=-\meti[X][Y],\quad \mbox{and}\quad \nabla_X JY = J \nabla_X Y,$$
which are the characteristic properties of an anti-K\"ahler manifold. 
On the other hand, any anti-K\"ahlerian metric (also called Norden metric) can be obtained from a uniquely determined holomorphic metric in this manner. The metric $\met$ may be written as a direct sum of its real and imaginary part as:
\begin{equation}
  \met = \meti + i \metii
\end{equation}
where $\meti$ is the (primary) Norden metric\index{Norden metric} and $\metii$ is the secondary Norden metric (also called the {\em twin metric}). These Norden metrics are related to each other by the relation:
\begin{equation}
\meti[X][Y] = \metii[X][JY],\quad \mbox{and}\quad \metii[X][Y] = -\meti[X][JY]
\end{equation}
for $X,Y\in T M$. \remarksymbol
\end{remark}

The easiest example of a holomorphic Riemannian manifold is the complex Euclidean $n$-space $\C^n$. The tangent space in a point being naturally identified with the vector space $\C^n$, the standard holomorphic metric on $\C^n$ is given by \eqref{EQ standard inner product}.

\section{Complexifications and real slices}\label{SC Complexifications}

\subsection{Complexifications of totally real submanifolds}

\begin{definition}
Let $Q$ be a real analytic manifold and $M$ a complex manifold. An immersion $f:Q\to M$ is called {\em totally real}\index{totally real!immersion} at a point $p\in Q$, if $TQ|_p$ is a totally real subspace of $TM|_p$ (i.e. $TQ|_p \cap i\,TQ|_p = \{0\}$). It is called {\em generic}\index{generic!immersion} at a point $p\in Q$, if $TQ|_p$ is a generic subspace of $TM|_p$ (i.e. $TQ|_p + i\,TQ|_p = TM|_p$). A totally real immersion is an immersion that is totally real at all points, and likewise for a generic immersion.
\end{definition}

In the above definition and further on, we identify $TQ|_p$ and $\d f(TQ|_p)$. It should be noted that, although the above definition of totally real is common in the area of complex analysis (cf. \cite{Bo}), it differs from the meaning it is usually given in the area of submanifold theory, where one commonly sees the additional requirement that for each point $p\in Q$, the subspace $TQ|_p$ is perpendicular to the subspace $iTQ|_p$ (cf. \cite{ChOg}).

\begin{definition}
Let $Q$ be a real analytic manifold. If a connected complex manifold, which we will denote by \emdex{$\mathbb{C} Q$}, is such that $Q$ is a generic totally real submanifold of $\C Q$, then we call $\C Q$ a \emdex{complexification} of $Q$. Furthermore, if $\C Q$ is a complexification of a real analytic manifold $Q$, and $f:Q\to M$ is a totally real analytic immersion of $Q$ into a complex manifold $M$, then a holomorphic map $\C f:\C Q \to M$\index{$\mathbb{C} f$} such that the restriction of $\C f$ to $Q$ equals $f$, is called a {\em holomorphic extension}\index{holomorphic extension} of $f$.
\end{definition}

Note that in the above definition, $f$ being an analytic immersion means that it is real-analytic with respect to the implicit real-analytic structure of the complex manifold $M$. The following theorem will be useful for the proof of later results (cf. \cite{WhBr}).

\begin{theorem}\label{ST unextendible extension}
Let $f:Q\to M$ be a totally real analytic immersion of a real analytic manifold $Q$ into a complex manifold $M$. Then there exists a complexification $\C Q$ of $Q$ for which a holomorphic extension $\C f:\C Q\to M$ of $f$ exists. Such a holomorphic extension is unique in the sense that any two such extensions are related through a biholomorphic map on a neighbourhood of $Q$.
\end{theorem}

The complex space $\C Q$ in the above theorem may be related to the notion of complexifications of an abstract manifold $Q$, as for instance in \cite{Ku}. However, it should be noted that in our situation, one may let the space $\C Q$ depend on the given immersion $f:Q\to M$.

\begin{example}
Consider the following immersions:
$$f_1:\R\to \C^2: t\mapsto (t,0),\quad f_2:\R\to \C^2: t\mapsto (\cosh t,\sinh t).$$ Although $f_1$ and $f_2$ have domains both homeomorphic to $\R$, the following holomorphic extensions 
$$\C f_1:\C\to \C^2: z\mapsto (z,0),\quad \C f_2:\C/(2\pi i\Z)\to \C^2: z\mapsto (\cosh z,\sinh z)$$ 
have complexified domains that are not homeomorphic.
\end{example}

Theorem \ref{ST unextendible extension} implies the following result, which we will need in the rest of this chapter.

\begin{corollary}\thlabel{ST analyt cont vector field}
Let $\pi:V\to M$ be a holomorphic vector bundle, and let $X:Q\to V|_Q$ be a real-analytic section over a real-analytic totally real submanifold $Q\subset M$. Then there exist a complexification $\C Q\subset M$ of $Q$, and a holomorphic extension of $X$ to a uniquely determined section $\C X\in \holsec(\C Q,V|_{\C Q})$.
\end{corollary}
\begin{proof}
First we observe that $X:Q\to V$ is a totally real immersion, for suppose that there is a $p\in Q$ and $v,w\in TQ|_p$, such that $\d X(v)=i\,\d X(w)$. We have that $\pi \circ X = \Id$, hence $\d\pi \circ \d X = \Id$, the identity operator on $TQ$. Using furthermore that $\d\pi$ is complex-linear, we see that
$$v = (\d\pi\circ \d X)(v) = \d\pi(i\,\d X(w)) = i\,(\d\pi\circ \d X)(w) = i\,w,$$ 
from which it follows that $v=w=0$, and thus $\d X(v)=\d X(w)=0$, proving that $X$ is a totally real immersion into $V$. Hence, by Theorem \ref{ST unextendible extension}, there exists a complexification $\C Q$ of $Q$ and a unique holomorphic extension $\C X: \C Q\to V$ of $X$. Now, it remains to prove that $\C X$ is a section over $\C Q$, i.e. $\pi\circ \C X = \Id|_{\C Q}$. This follows from the fact that $\pi\circ \C X$ equals the identity on $Q$, so that both $\pi\circ \C X$ and $\Id|_{\C Q}$ are a holomorphic extension of $\Id|_Q$. Hence, by the unicity part of Theorem \ref{ST unextendible extension}, it follows that $\pi\circ\C X = \Id|_{\C Q}$.
\end{proof}

\subsection{Real slices}

\begin{definition}\label{DF real slice}
Given a complex manifold $M$ with complex Riemannian metric $\met$, we call a submanifold $N\subset M$ a {\em real slice}\index{real slice!of complex manifold} of $(M,\met)$, if at any point $p\in N$ we have that $TN|_p$ is a real slice of $(TM|_p,g)$ in the sense of Definition \ref{DF real slice hol in prod space}.
\end{definition}

The above definition says that the metric $\met$ restricted to the real tangent bundle of $N$, is real valued. Hence, the restriction $\met|_{TN}$ turns $N$ into a (pseudo-)Riemannian manifold.
This notion of real slices is in accordance with the concept of real slices as introduced in \cite{Le}. It is easily seen that a submanifold $N$ is a real slice if and only if the secondary Norden metric $\metii|_M$ vanishes entirely on $M$, implying that the induced complex metric $\met|_{TN}$ coincides with the induced primary Norden metric $\meti|_{TN}$. 

The following two propositions are immediate consequences of Definitions \ref{DF real slice hol in prod space} and \ref{DF real slice} and Proposition \ref{ST real slices hol in prod space}.

\begin{proposition}\label{ST real slice}
A real slice of a complex Riemannian manifold is a totally real submanifold.
\end{proposition}

\begin{proposition}\label{ST submanifold of real slice}
Let $M$ be a complex Riemannian manifold and $P$ a real slice of $M$. Then any real submanifold $Q\subset P$ is also a real slice of $M$. 
\end{proposition}

This last proposition brings along the following useful corollary.

\begin{corollary}\label{ST intersection is real slice}
Let $M$ be a complex Riemannian manifold, $N$ a complex submanifold of $M$ and $P$ a real slice of $M$. Then the intersection $N\cap P$ is a real slice of $N$, provided it is non-empty.
\end{corollary}

It should be noted that for a given complex manifold $M$, only very specific (pseudo-)Riemannian manifolds may appear as generic real slices of it. For if $Q$ is a generic real slice of $M$, then a complexification $\C Q$ is an open subset of $M$. This fully determines the complex Riemannian space $M$ in a neighborhood around $Q$. We will make the following definition (cf. \cite{Br}).

\begin{definition}\label{DF Wick-related}
Two (pseudo-)Riemannian manifolds $P$ and $Q$ are said to be \emdex{Wick-related} if there exists a holomorphic Riemannian manifold $(M,\met)$, such that $P$ and $Q$ are embedded as real slices of $M$.
\end{definition}

As will be discussed later in more detail, when two manifolds are Wick-related, their geometries are intimately related through their common complexified ambient space. 

As many geometric properties can be expressed as the vanishing of a certain tensor, the following two lemmas will turn out to be useful for the applications in Section \ref{SC applications}. The first one is the following very trivial observation.

\begin{lemma} \label{ST from complex-analogous to pseudo-Riem}
If $N$ is a complex $n$-dimensional submanifold of a holomorphic Riemannian manifold $M$ for which a certain holomorphic tensor field $F$ vanishes, then for any real slice $P$ of $M$ for which the intersection $Q=N\cap P$ is a real $n$-dimensional submanifold, we have that $F|_Q$ vanishes on $Q$ as well.
\end{lemma}

The second lemma follows from \thref{ST analyt cont vector field}, and is in some sense a converse to Lemma \ref{ST from complex-analogous to pseudo-Riem}.

\begin{lemma} \label{ST from pseudo-Riem to complex-analogous}
Let $M$ be a holomorphic Riemannian manifold, and let $P$ be a generic real slice of $M$. Then given a real-analytic (pseudo-)Riemannian submanifold $Q\subset P$ on which a particular real-analytic tensor $F$ vanishes, there exists a complexification $N=\C Q\subset M$ for which the holomorphic extension $\C F$ of $F$ vanishes on $N$.
\end{lemma}

\section{Holomorphic Riemannian submanifold theory}\label{SC Hol Riem submanifold theory}

We will start this section with a brief overview of {\em holomorphic Riemannian submanifold theory}, the holomorphic Riemannian counterpart of \mbox{(pseudo-)}\linebreak[0]Riemannian submanifold theory. As we will see shortly, in the setting of holomorphic Riemannian submanifold theory, all tensors and operators of interest turn out to be holomorphic in nature, meaning that they are respectively either holomorphic tensors, or operators that result in a holomorphic tensor if all of their arguments are holomorphic. This is an important property, as it allows for holomorphic extensions of any vector field constructed by them.

It is well-known that on a complex manifold, given holomorphic vector fields $X,Y\in \holsec(M,TM)$, we have that $[X,Y]$ is again a holomorphic vector field, and $[X,cY]=[cX,Y]=c[X,Y]$ for any $c\in\C$ (cf. \cite{KoNo}). 

On a holomorphic Riemannian manifold $(M,\met)$, the metric $\met$ is a holomorphic tensor by definition, and it is a complexification of its restriction to a (pseudo-)Riemannian metric $\met|_Q$, for any real slice $Q\subset M$. It follows from Koszul's formula for the holomorphic Riemannian Levi-Civita connection that, given $X,Y\in \holsec(M,TM)$, $\nabla_X Y$ is a holomorphic vector field as well. Furthermore, for a holomorphic function $f$ and $c\in \C$, we have that $\nabla_{f X} Y = f \nabla_X Y$ and $\nabla_{X} c Y = c \nabla_X Y$ (cf. \cite{Sl}).

The holomorphic Riemannian curvature tensor is defined by
$$R(X,Y)Z = \nabla_X\nabla_Y Z - \nabla_Y\nabla_X Z - \nabla_{[X,Y]}Z.$$
It follows from a simple calculation that this indeed defines a $\C$-linear tensor. This tensor is moreover holomorphic by the previous observations on the holomorphy of the Lie bracket and the Levi-Civita connection. 

We can then use $R$ to define a holomorphic Riemannian analogue of the sectional curvature\index{holomorphic Riemannian sectional curvature} as follows:
\begin{equation}\label{EQ sectional curvature}
K: \holsec(M,T^2M^{\mathrm{n.d.}})\to \C: X \otimes Y \mapsto \frac{\met[R(X,Y)Y][X]}{\met[X][X]\met[Y][Y]-\met[X][Y]^2},
\end{equation}
where $T^2M^{\mathrm{n.d.}}$ denotes the subbundle of $(2,0)$-tensors $X\otimes Y$ in $T^2M$, for which $X$ and $Y$ span a non-degenerate complex plane, i.e., for which the denominator in \eqref{EQ sectional curvature} is non-zero. Then $K$ is a holomorphic map and this definition corresponds to the definition of complex sectional curvature in \cite{DoKoPe}. 

In similar ways, we can define the holomorphic Riemannian analogues of the Ricci-curvature, the scalar curvature, etc. As long as these holomorphic Riemannian analogues of common (pseudo-)Riemannian quantities are defined through operations that preserve holomorphy, the result will be holomorphic as well.

We will now continue with defining some extrinsic tensors and operators for submanifolds in the holomorphic Riemannian setting. Let $(M,\met)$ be a holomorphic Riemannian manifold and $N$ a holomorphic Riemannian submanifold. We denote the Levi-Civita connection of the ambient space $M$ by $D$, and the induced Levi-Civita connection on the submanifold $N$ by $\nabla$.

\begin{proposition} \label{prop:topandbot}
The projections $(\cdot)^\top$ and $(\cdot)^\bot$ of $TM|_N =TN\oplus (TN)^\bot$ to $TN$ and $(TN)^\bot$ respectively, are both holomorphic tensor fields.
\end{proposition}
\begin{proof}
For any point $p\in N$, it is possible through a Gram-Schmidt orthonormalization process to construct a holomorphic orthonormal frame $\{X_1,\ldots,X_n\}$ on some neighborhood $U$ of $p$ in $N$. 
Then for any section $Y \in \holsec(U,TM|_U)$, one has
\begin{equation*} 
Y^\top = \sum_{i=1}^n \met[X_i][Y]X_i. 
\end{equation*}
From this expression, it is seen that the tangent projection $(\cdot)^\top$ is a holomorphic tensor, and consequently, that the orthogonal projection $(\cdot)^\bot=\Id-(\cdot)^\top$ is a holomorphic tensor as well.
\end{proof}

An important consequence of Proposition \ref{prop:topandbot} is the following.

\begin{proposition}
The holomorphic Riemannian analogues of the second fundamental form $h$ and the shape operator $A$ are holomorphic tensors.
\end{proposition}
\begin{proof}
Let $M$ and $N$ be as above, and let $X$ and $Y$ be vector fields tangent to $N$ and $\xi$ a vector field normal to $N$. In the holomorphic Riemannian setting, the formulas of Gauss and Weingarten read identical to the ones in the (pseudo-)Riemannian setting, which is:
\begin{align}\label{EQ Gauss and Weingarten}
D_X Y &= \nabla_X Y + h(X,Y),\\
D_X \xi &= -A_\xi X + \nabla_X^\bot \xi,
\end{align}
where $\nabla^{\bot}$ stands for the orthogonal connection induced on the submanifold $N$. Hence $h$ and $A$ can be expressed as: 
\begin{equation} \label{eq:handA}
h(X,Y) = (D_X Y)^\bot,\qquad A_{\xi} X = -(D_X \xi)^\top.
\end{equation}
From the previous proposition, it immediately follows that $h$ and $A$ are holomorphic, and tensoriality is verified by a straightforward calculation.
\end{proof}
Note that the orthogonal counterparts of the equations \eqref{eq:handA} are
\begin{equation*}
\nabla_X Y = (D_X Y)^\top,\qquad \nabla_X^\bot \xi = (D_X \xi)^\bot,
\end{equation*}
from which can be seen that the induced connection on $N$, which is of course the Levi-Civita connection of $(N,\met)$, and the normal connection are $\C$-linear holomorphic operators. 

Based on these fundamental extrinsic operators, many other extrinsic quantities can be defined. An important example of this is the {\em holomorphic Riemannian mean curvature}\index{holomorphic Riemannian mean curvature}
$$ H: N \to TM: p \mapsto \frac{1}{n} \mathrm{Tr}(h|_p) = \frac 1n \sum_{i=1}^n h(X_i,X_i) \in (TN)^{\bot}|_p, $$ 
where $X_1,\ldots,X_n$ is a orthonormal basis of $TN|_p$. Again, it is seen immediately that this defines a holomorphic tensor field.

\begin{remark}\label{RM holomorphic extension}
All of the holomorphic Riemannian tensors and operators discussed in this section, coincide with their (pseudo-)Riemannian counterparts when we restrict the (ambient) space to a generic real slice. By this we mean that if we take as input arguments some holomorphic extensions of real tensor fields over $Q$, then the output will be a holomorphic extension of the (pseudo-)Riemannian output on these real tensor fields. 
Indeed, the expressions through which these operators are defined are exactly the same as the way their (pseudo-)Riemannian counterparts would be defined on a real slice $Q$, except that the inputs are required to be holomorphic tensor fields. 
So in this sense, all the holomorphic Riemannian tensors and operators discussed above, are holomorphic extensions of their corresponding (pseudo-)Riemannian operators.
\end{remark}

\begin{remark}
Many well-known equations from (pseudo-)Riemannian submanifold theory, 
go through unaltered in the setting of holomorphic Riemannian submanifold theory. In particular we have the equations of Gauss and Ricci, i.e.
\begin{align}
\label{EQ Gauss equation}\met[\tilde R(X,Y)Z][W] &= \met[R(X,Y)Z][W]+\met[h(Y,W)][h(X,Z)]\\
 					&\qquad -\met[h(X,W)][h(Y,Z)],\nonumber\\
\label{EQ Codazzi equation}\met[\tilde R(X,Y)\xi][\eta] &= \met[R(X,Y)\xi][\eta]-\met[{[A_\xi,A_\eta]X}][Y],
\end{align}
where $\tilde R$ and $R$ denote the Riemann tensor on respectively the ambient manifold and the submanifold. These equations are identical to the equations in the ordinary (pseudo-)Riemannian case, as is the algebraic derivation through which these equations may be obtained (cf. \cite{Do}).
\end{remark}

\subsection{The standard holomorphic Riemannian space forms}

As mentioned above, an important example of a holomorphic Riemannian manifold is the {\em complex Euclidean $n$-space}\index{complex Euclidean space} $\C^n$. The curvature of this manifold vanishes and its isometry group is $E(n,\C) = \C^n \rtimes \mathrm O(n,\C)$, where $\C^n$ acts by translations. Some other important transformations of $\C^n$ are given by scalar multiplications with complex numbers $\alpha \neq 0$. It is clear that $L_\alpha: \C^n\to \C^n: z \mapsto \alpha z$ is a conformal transformation with conformal factor $\alpha^2$, i.e.,
$$\met[\d L_{\alpha}(X)][\d L_{\alpha}(Y)]|_{\alpha z} = \alpha^2 \met[X][Y]|_{z}$$
for all $z \in \C^n$ and all $X,Y\in T\C^n|_z$. In the particular cases that $\alpha=\pm 1$ and $\alpha = \pm i$, this gives respectively an isometry and an anti-isometry on the entire space. All dilations together with the group $E(n,\C)$ of isometries, generate the group of all conformal orthogonal transformations of the space $\C^n$, i.e. affine transformations which preserve orthogonality. In the following we will call such transformations {\em similarities}\index{similarity}, and submanifolds that can be related to each other by a similarity are called {\em similar}.

Clearly, the generic real slices of $\C^n$ as a holomorphic inner product space (cf. Example \ref{EX pseudo-Euclidean}) are also generic real slices of $\C^n$ as a complex Euclidean space. Moreover, any translation of such a real slice $\R_k^n\subset \C^n$ over a complex vector $z\in\C^n$ results in an affine subspace $z+\R_k^n$, which is a real slice with the same signature. Of course, this space only differs from the former one if $z \notin V$. We will refer to such real slices of $\C^n$ as the {\em generic affine real slices}\index{generic affine real slice}.

Another important example of a holomorphic Riemannian manifold is the {\em holomorphic Riemannian sphere}\index{holomorphic Riemannian sphere} $\C S^n(\alpha)$ of radius $\alpha$, where $\alpha$ can be any non-zero complex number:
$$\C S^n(\alpha) = \{ (z_1,\ldots,z_{n+1}) \in \C^{n+1} \ | \ z_1^2+\ldots+z_{n+1}^2 = \alpha^2 \}.$$ 
The holomorphic metric is induced by the holomorphic metric on $\C^{n+1}$. In the following, we will simply write $\C S^n$ for $\C S^n(1)$. All the $\C S^n(\alpha)$ are conformal since the dilation $L_\alpha:\C^{n+1}\to \C^{n+1}$ maps $\C S^n$ to $\C S^n(\alpha)$. The manifold $\C S^n(\alpha)$ has constant sectional curvature $1/\alpha^2$. In particular, the manifold $\C S^n(\alpha)$ is equal to the manifold $\C S^n(-\alpha)$, and anti-isometric to the manifolds $\C S^n(\pm i \alpha)$. As is suggested by the notation, the manifold $\C S^n$ is indeed a complexification of the sphere $S^n$ in $\R^{n+1}$.

As follows from Corollary \ref{ST intersection is real slice}, the intersection of $\C S^n$ with a real slice of $\C^{n+1}$ will be a real slice of $\C S^n$. If we consider the generic real slice $\R_k^{n+1}$ (cf. Example \ref{EX pseudo-Euclidean}) of $\C^{n+1}$ through the origin of $\C^{n+1}$, we obtain as real slice of $\C S^n$ the indefinite sphere (also called generalized de Sitter space)
\begin{equation*}
S^n_k = \{ (y_1,\ldots,x_{n+1}) \in \R^{n+1}_k \ | \ -y_1^2-\ldots-y_k^2+x_{k+1}^2+\ldots+x_{n+1}^2 = 1 \}.
\end{equation*}
For $k=0$, this is simply the $n$-dimensional round unit sphere $S^n$. 

It is known that the indefinite sphere $S_k^n$ is anti-isometric to the indefinite hyperbolic $n$-space (also called generalized anti-de Sitter space) $H_{n-k}^n$, where $H_k^n$ is defined by:
\begin{equation*}\label{EQ indef hyperbolic}
H^n_{k} = \{ (y_1,\ldots,x_{n+1}) \in \R^{n+1}_{k+1} \ | \ -y_1^2-\ldots-y_{k+1}^2+x_{k+2}^2+\ldots+x_{n+1}^2 = -1 \}.
\end{equation*}
That $S_k^n$ and $H_{n-k}^n$ are anti-isometric, can now be seen from the fact that by applying the anti-isometry $L_i$ to $\C^{n+1}$ (which also induces an anti-isometry between submanifolds and their image under $L_i$), we have that the submanifold $S_k^n$ is turned into $L_i S_k^n = L_i (\C S^{n} \cap \R_k^{n+1}) = L_i \C S^{n} \cap L_i \R_k^{n+1} = \C S^n(i) \cap L_i \R_k^{n+1}$. Now it follows from \eqref{EQ pseudo-Euclidean} that
$$L_i \R_k^{n+1} = \mathrm{span}_{\R}\{e_1,\ldots,e_k,i e_{k+1},\ldots,e_{n+1}\} \cong \R_{n-k+1}^{n+1},$$ so that $\C S^n(i) \cap L_i \R_k^{n+1}$ corresponds (after a rearrangement of coordinates) to $H_{n-k}^n$. Note that for $k=n$, we obtain the ordinary hyperbolic $n$-space.

So we have seen that on the one hand all the indefinite spheres are Wick-related, and on the other hand that all the indefinite hyperbolic spaces are Wick-related. Also, the indefinite spheres can be related with the indefinite hyperbolic spaces trough a dilation $L_i$ on the ambient space $\C^{n+1}$, which defines an anti-isometry between $S_k^n$ as a real slice of $\C S^n$, and $H_{n-k}^n$ as a real slice of $\C S^n(i)$. Hence, although the $n$-sphere is not directly Wick-related to hyperbolic $n$-space, they are in a sense Wick-related up to an anti-isometry. This means that a geometric problem in the $n$-dimensional sphere can be translated to an analogous problem in $n$-dimensional hyperbolic space, provided that the properties of interest are preserved under anti-isometric mappings. An easy example of this will be given in Section \ref{SC surfaces in CS^n}.

\section{Applications}\label{SC applications}

In this section we will demonstrate by some elementary examples in submanifold theory, how the above theory can be applied to reveal geometric correspondences among Wick-related manifolds. The main idea is as follows. As we know from Proposition \ref{ST real slice}, real slices are always totally real submanifolds. Hence, by \thref{ST analyt cont vector field} we have that any real-analytic tensor field over a generic real slice can be extended to a holomorphic tensor field on an open set of the ambient space. If this ambient space also contains other generic real slices, then by the uniqueness of the analytic extension, we have established a one-to-one correspondence between such tensor fields on these (pseudo-)Riemannian manifolds. As many geometric properties can be defined by the vanishing of a certain tensor, Lemmas \ref{ST from complex-analogous to pseudo-Riem} and \ref{ST from pseudo-Riem to complex-analogous} allow us to translate knowledge about such geometric properties among Wick-related manifolds (cf. Definition \ref{DF Wick-related}). Moreover, by proving a statement for the ambient holomorphic Riemannian manifold, we obtain similar statements simultaneously for all the entailed real-analytic real slices of this manifold. 

In the next three examples, we consider minimal, totally geodesic and parallel submanifolds in respectively $\C^3$, $\C S^3$ and $\C^3$. Since these geometric properties are all invariant under similarities of the ambient space, the submanifolds only need to be characterized up to similarity. 

\subsection{Minimal surfaces of revolution in holomorphic \texorpdfstring{$\C^3$}{C3}}\label{EX catenoid}

It is well-known that in $\R^3$, any minimal surface of revolution is similar to either a plane, or to the catenoid $C\subset \R^3$, given by
$$ C = \{ (x_1,x_2,x_3)\in\R^3 \ | \ x_1^2+x_2^2=\cosh^2 x_3 \}. $$
From this equation, it is indeed easily seen that $C$ is invariant under an $SO(2,\R)$ action, where $SO(2,\R)$ is represented as the subgroup of $SO(3,\R)$ given by the matrices of the form
\begin{equation}\label{EQ matrix revol}
\begin{pmatrix}
\cos \theta & -\sin \theta & 0 \\ 
\sin \theta & \cos \theta & 0 \\ 
0 & 0 & 1
\end{pmatrix}
\end{equation}
where $\theta\in\R$. Now by Lemma \ref{ST from pseudo-Riem to complex-analogous}, we obtain a minimal holomorphic Riemannian surface $\C C\subset \C^3$, given by
\begin{equation}\label{EQ catenoid in C3 1}
\C C = \{ (z_1,z_2,z_3)\in\C^3 \ | \ z_1^2+z_2^2=\cosh^2 z_3 \}. 
\end{equation}
Here and in the following, \emph{minimal} means that the (holomorphic Riemannian) mean curvature vector field $H$ vanishes everywhere. The surface \eqref{EQ catenoid in C3 1} is also easily seen to be invariant under the subgroup $SO(2,\C)\subset SO(3,\C)$, formed by the matrices \eqref{EQ matrix revol}, but now for $\theta\in \C$. This group contains both the Euclidean rotation group $SO(2,\R)$ as well as the indefinite group $SO(1,1,\R)$, given by all matrices of the form
$$\pm\begin{pmatrix}
\cosh \theta & \sinh \theta \\ 
\sinh \theta & \cosh \theta
\end{pmatrix}$$
where $\theta\in \R$. However, this group does not contain (conjugates of) the special kind of rotations in $SO(3,\C)$ around a light-like vector such as $e_1+i e_3$. Note that the following four complex surfaces are all similar to $\C C$ in $\C^3$:
\begin{align}
& \{ (z_1,z_2,z_3)\in\C^3 \ | \ z_1^2+z_2^2=-\sinh^2 z_3 \}, \label{EQ catenoid in C3 2}\\
& \{ (z_1,z_2,z_3)\in\C^3 \ | \ z_2^2+z_3^2=\cosh^2 z_1 \}, \label{EQ catenoid in C3 3}\\
& \{ (z_1,z_2,z_3)\in\C^3 \ | \ z_1^2+z_2^2=-\sin^2 z_3 \}, \label{EQ catenoid in C3 4}\\
& \{ (z_1,z_2,z_3)\in\C^3 \ | \ z_2^2+z_3^2=-\sin^2 z_1 \}. \label{EQ catenoid in C3 5}
\end{align}
The surface \eqref{EQ catenoid in C3 2} is obtained from \eqref{EQ catenoid in C3 1} by a translation over a distance $\pi i/2$ in the direction of the $z_3$ axis and the surface \eqref{EQ catenoid in C3 4} is then obtained from \eqref{EQ catenoid in C3 2} by performing a dilation $L_i$ (which is an anti-isometry) on the ambient space $\C^3$. The surfaces \eqref{EQ catenoid in C3 3} and \eqref{EQ catenoid in C3 5} are obtained from \eqref{EQ catenoid in C3 1} and \eqref{EQ catenoid in C3 4} respectively, by permuting the coordinates $z_1$ and $z_3$. Hence \eqref{EQ catenoid in C3 3} and \eqref{EQ catenoid in C3 5} are rotation symmetric with respect to the $z_1$-axis. All these surfaces are mutually similar in the holomorphic Riemannian space $\C^3$, so there is no need to distinguish between them in $\C^3$. However, as we will see next, after intersecting them with a certain real slice, we may obtain surfaces that are not similar with respect to the similarities of this real slice (i.e., the group of affine transformations that preserve orthogonality within the real slice), thus giving rise to essentially different pseudo-Riemannian surfaces.

Of the five similar surfaces above, only \eqref{EQ catenoid in C3 1} and \eqref{EQ catenoid in C3 3} have non-empty intersection with the real slice $\R^3: y_1=y_2=y_3=0$, and in both cases we recover the classical catenoid. 

By intersecting the five surfaces with the real slice $\R_1^3:x_1=y_2=y_3=0$ on the other hand, we obtain the following five minimal surfaces (for the naming we follow \cite{An} and \cite{Wo}):
\begin{enumerate}
\item the Lorentzian hyperbolic catenoid, given by $-y_1^2+x_2^2=\cosh^2 x_3$,\label{LI Cat lor hyp}
\item the Lorentzian hyperbolic catenoid of the second kind, given by \linebreak[4] $y_1^2-x_2^2=\sinh^2 x_3$,\label{LI Cat lor hyp 2}
\item the Lorentzian elliptic catenoid, given by $x_2^2+x_3^2=\cos^2 y_1$,\label{LI Cat lor ell}
\item the spacelike hyperbolic catenoid, given by $y_1^2-x_2^2=\sin^2 x_3$,\label{LI Cat spa hyp}
\item the spacelike elliptic catenoid, given by $x_2^2+x_3^2=\sinh^2 y_1$.\label{LI Cat spa ell}
\end{enumerate}

Surfaces \eqref{LI Cat lor ell} and \eqref{LI Cat spa ell} are invariant under the group of Euclidean rotations (and reflections) $O(2,\R)\subset O(2,\C)$ whose axis of symmetry is the time-like $y_1$-axis. Surfaces \eqref{LI Cat lor hyp},\eqref{LI Cat lor hyp 2} and \eqref{LI Cat spa hyp} are invariant under the hyperbolic rotations (and reflections) $O(1,1,\R)\subset O(2,\C)$. Note that using only the property that the classical catenoid is a minimal surface in $\R^3$, we have obtained all of the other surfaces above, and their minimality is guaranteed without the need for any calculations.

\begin{remark}
In this and the following examples, the submanifolds are all described by implicit equations. This is done for convenience, but it should be noted that the possibility to use parametrizations is certainly there. For example, the six rotationally invariant minimal surfaces above (the one in $\R^3$ and the other five in $\R_1^3$), may respectively be parametrized as follows:
\begin{enumerate}\setcounter{enumi}{-1}
\item $L(u,v) = (\cos u \cosh v, \sin u \cosh v, v)$,\label{LI par Cat}
\item $L(u,v) = (i \sinh u \cosh v, \cosh u \cosh v, v)$,\label{LI par Cat lor hyp}
\item $L(u,v) = (i \cosh u \sinh v, \sinh u \sinh v, v)$,\label{LI par Cat lor hyp 2}
\item $L(u,v) = (i\, v, \cos u \cos v, \sin u \cos v)$,\label{LI par Cat lor ell}
\item $L(u,v) = (i \cosh u \sin v, \sinh u \sin v, v)$,\label{LI par Cat spa hyp}
\item $L(u,v) = (i\, v, \cos u \sinh v, \sin u \sinh v)$.\label{LI par Cat spa ell}
\end{enumerate}
It is worth noting how for example case \eqref{LI par Cat} is turned into case \eqref{LI par Cat lor ell} after multiplying the coordinate $u$ by $i$, a so-called Wick rotation. Just with some additional translations on the coordinates and suitable coordinate permutation on the ambient space, many of the other cases are obtained likewise.\remarksymbol
\end{remark}

\subsection{Totally geodesic surfaces in the holomorphic sphere \texorpdfstring{$\C S^n$}{CSn}}\label{SC surfaces in CS^n}

In this example we will see how the method of Wick-relations can be applied to situations where we have an ambient space other than $\C^n$. Before we come to this, it is good to realize that for a generically chosen holomorphic Riemannian manifold, one cannot even expect to find real slices at all (cf. \cite{Br}). Nevertheless, examples of holomorphic Riemannian manifolds that do contain real slices are easily constructed, as follows from Corollary \ref{ST intersection is real slice}: if $M$ is a holomorphic Riemannian submanifold of $\C^n$, which has a non-empty intersection with a certain real slice $\R_k^n$ of $\C^n$, then this intersection is itself a real slice of $M$.  We have already seen an application of this before, in obtaining $S_k^n$ as a real slice of $\C S^n$, and it can be applied to many other situations as well. For example, if we consider the complex surface $\C C \subset \C^3$ from the previous subsection, then the catenoid in $\R^3$ and the five minimal surfaces in $\R_1^3$ listed above are real slices of $\C C$, and thus Wick-related.

Just as a simple example, let us focus on 
$$ \C S^3 = \{ (z_1,z_2,z_3,z_4) \in \C^4 \ | \ z_1^2+z_2^2+z_3^2+z_4^2=1 \} $$
and consider the totally geodesic submanifold $\C S^2$ given by $z_1^2+z_2^2+z_3^2=1$ and $z_4=0$. We know this is a totally geodesic embedding by Lemma \ref{ST from pseudo-Riem to complex-analogous}, for if we intersect $\C S^2$ with the real slice $S^3\subset \C S^3$, given by $y_1=y_2=y_3=y_4=0$ (and hence $x_1^2+x_2^2+x_3^2+x_4^2=1$), we obtain the totally geodesic surface $S^2\subset S^3$, given by $x_1^2+x_2^2+x_3^2=1$ and $x_4=0$. From this, we know for sure that by intersecting $\C S^2$ (or submanifolds similar to it, such as $z_1^2+z_2^2+z_4^2=1$ and $z_3=0$ etc.) with one of the other real slices of $\C S^3$, we obtain totally geodesic submanifolds in those respective spaces as well. In this way, we find the following totally geodesic surfaces:
\begin{enumerate}
\item $S^2$ and $S_1^2$ in $S_1^3$, or equivalently $H_2^2$ and $H_1^2$ in $H_2^3$,
\item $S_1^2$ and $S_2^2$ in $S_2^3$, or equivalently $H_1^2$ and $H^2$ in $H_1^3$,
\item $S_2^2$ in $S_3^3$, or equivalently $H^2$ in $H^3$.
\end{enumerate} 
Note that the equivalent statement about embeddings of indefinite hyperbolic spaces is obtained through applying the anti-isometry $L_i$ on the ambient space $\C S^3$, which links each totally geodesic embedding of an indefinite $2$-sphere in an indefinite $3$-sphere to a totally geodesic embedding of some indefinite hyperbolic plane in an indefinite hyperbolic $3$-space. In particular, the totally geodesic embedding of the hyperbolic plane $H^2 : x_1^2+x_2^2-y_4^2=-1$ inside the hyperbolic space $H^3 :x_1^2+x_2^2+x_3^2-y_4^2=-1$ is thus obtained. Importantly, this simple example demonstrates how Wick-relationships can be used not only to reveal correspondences between pseudo-Riemannian spaces (possibly of different signature), but even to reveal correspondences between different kinds of ordinary Riemannian spaces (in this case the $n$-dimensional sphere and $n$-dimensional hyperbolic space).

\subsection{Parallel surfaces in holomorphic \texorpdfstring{$\C^3$}{C3}}\label{EX Par}

In the following, we will classify surfaces in $\C^3$ with parallel second fundamental form, i.e., surfaces for which the holomorphic tensor $\nabla h$ vanishes identically. These are called parallel surfaces for short. 

\begin{lemma}\label{ST forms complex matrix}
Let $A$ be a linear operator on a complex two-dimensional vector space $V$, which is symmetric with respect to a holomorphic inner product $\met$ on $V$. Then there exists an orthonormal basis $\{e_1,e_2\}$ of $V$ such that, with respect to $\{e_1,e_2\}$, $A$ takes one of the following forms:
\begin{eqnarray}
A=\begin{pmatrix}
\alpha & 0 \\
0 & \beta \end{pmatrix}\quad\mbox{or}\label{EQ diagonalizable}\\
A = \begin{pmatrix}
\alpha+1 & i \\
i & \alpha-1 \end{pmatrix},\label{EQ undiagonalizable}
\end{eqnarray}
where $\alpha$ and $\beta$ are complex numbers.
\end{lemma}
\begin{proof}
Let $\{u_1,u_2\}$ be a any orthonormal basis of $V$. Since $A$ is symmetric, it can be written in the form 
$$A = 
\begin{pmatrix}
m+a & b \\
b & m-a 
\end{pmatrix}$$ 
with respect to $\{u_1,u_2\}$, where $a$, $b$ and $m$ are all complex. Now consider a general orthonormal basis given by $e_1 = \cos t\,u_1 + \sin t\,u_2$ and $e_2 = -\sin t\,u_1 + \cos t\,u_2$, where $t$ is a complex number. We want to choose $t$ in such a way that $A$ takes one of the two forms above.

If $b=0$, we are done immediately, so let us assume $b\neq 0$.

Provided that $a \neq \pm b i$, we may choose $t$ a complex number such that $\cot(2t) = a/b$, and one verifies that $A$ takes the form \eqref{EQ diagonalizable} with respect to $\{e_1,e_2\}$. (Note that $\cot t$ assumes all complex values except for $\pm i$.)

In case $a = b i$, we can choose $t$ such that $e^{-2ti}=a$, and $A$ will take the form \eqref{EQ undiagonalizable} with respect to $\{e_1,e_2\}$. In case $a = -b i$, we choose $t$ such that $e^{2ti}=a$ to get the same result.
\end{proof}

\begin{theorem}\label{ST classification hol riem}
A parallel surface $M$ in $\C^3$ is similar to an open part of one of the following four surfaces:
\begin{enumerate}
\item the complex plane $\C^2$, given by $z_3 = 0$,\label{case 1}
\item the complex sphere $\C S^2$, given by $z_1^2+z_2^2+z_3^2=1$,\label{case 2}
\item the cylinder $\C S^1\times \C$, given by $z_1^2+z_2^2=1$,\label{case 3}
\item the flat minimal surface given by $z_3 = (z_1+iz_2)^2$.\label{case 4}
\end{enumerate}
\end{theorem}
\begin{proof}
Assume $M$ is a complex submanifold with parallel holomorphic Riemannian second fundamental form $h$, and let $L:U\subset \C^2\to \C^3$ denote a local parametrization of $M\subset \C^3$. Let $e_3$ be a (holomorphic Riemannian) unit normal vector field. First suppose that with respect to a local orthonormal frame field $\{e_1,e_2\}$, the shape operator $A$ takes the form \eqref{EQ diagonalizable}. A direct computation shows that the surface is parallel if and only if $\alpha$ and $\beta$ are constant and $(\alpha-\beta)\omega_1^2=0$, where $\omega_1^2$ is a one-form defined by $\nabla_X e_1 = \omega_1^2(X) e_2$.

If $\alpha=\beta=0$, then the surface is totally geodesic, i.e. $h = 0$, which gives case \eqref{case 1}. 

If $\alpha=\beta\ne 0$, then the second fundamental form satisfies:
$$h(e_1,e_1)=\alpha e_3,\quad h(e_1,e_2)=0,\quad h(e_2,e_2)=\alpha e_3.$$
Thus it follows from the Gauss equation \eqref{EQ Gauss equation} that $K=\alpha^2$. By choosing (holomorphic) geodesic coordinates $(u,v)$ such that $ds^2 = du^2 + \cos^2(\alpha u) dv^2$, we have that
\begin{align*}
L_{uu} = \alpha e_3,\quad L_{uv} &= - \alpha \tan(\alpha u) L_v,\quad L_{vv} = \alpha \cos^2(\alpha u)e_3 + \frac \alpha 2 \sin(2 \alpha u) L_u,\\ 
&\quad D_{\del_u}e_3 = -\alpha L_u, \quad D_{\del_v}e_3 = -\alpha L_v.
\end{align*}
The solution of this system of equations is, up to translation, given by
\begin{align*}
L(u,v) &=  \frac{1}{\alpha}\left(\cos(\alpha u) \cos(\alpha v) w_1 +  \cos(\alpha u) \sin(\alpha v) w_2 + \sin(\alpha u) w_3\right)\\
e_3(u,v) &= -\cos(\alpha u) \cos(\alpha v) w_1 - \cos(\alpha u) \sin(\alpha v) w_2 - \sin(\alpha u) w_3
\end{align*}
where $\{w_1,w_2,w_3\}$ is an orthonormal basis of $\C^3$ (as holomorphic inner product space). After rescaling, we obtain case \eqref{case 2} of the theorem.

If $\alpha\neq \beta$, then $\omega_1^2 = 0.$ Hence $M$ is flat, so that by the equation of Gauss we have $\operatorname{det} A = \alpha \beta = 0$. Without loss of generality, we may assume $\beta = 0$. Now, choose coordinates $(u,v)$ on $M$ with $\del_u = e_1$ and $\del_v = e_2$. Then, the immersion $L$ satisfies
$$
L_{uu} = \alpha e_3,\quad L_{uv} = 0,\quad L_{vv} = 0,\quad D_{\del_u}e_3 = -\alpha L_u, \quad D_{\del_v}e_3 = 0.
$$
The solution of this system of equations is, up to translation, given by
\begin{align*}
L(u,v) &= \frac{1}{\alpha}\left(\cos(\alpha u) w_1 + \sin(\alpha u) w_2\right) + v\,w_3,\\
e_3(u,v) &= -\cos(\alpha u) w_1 - \sin(\alpha u) w_2,
\end{align*}
where $\{w_1,w_2,w_3\}$ is an orthonormal basis of $\C^3$. After applying a suitable similarity in $\C^3$, we obtain case \eqref{case 3} of the theorem.

Now suppose that with respect to an orthonormal frame field $\{e_1,e_2\}$, $A$ takes the form \eqref{EQ undiagonalizable}. Since the surface is parallel, we obtain that $\alpha$ is constant and $\omega_1^2 = 0$. Hence the surface is flat. But from the equation of Gauss, we obtain that the Gaussian curvature is given by $\operatorname{det} A = \alpha^2$, so that $\alpha=0$. If we take coordinates $(u,v)$ on $M$ with $\del_u = e_1$ and $\del_v = e_2$, then the formulas of Gauss and Weingarten yield the following system of equations:
\begin{align*}
\quad L_{uu} &= e_3,\quad L_{uv} = i e_3,\quad L_{vv} = - e_3,\\
\quad D_{\del_u}e_3 &= -L_u-i L_v, \quad D_{\del_v}e_3 = -i L_u+L_v.
\end{align*}
The solution of this system of equations is, up to translation, given by
\begin{align*}
L(u,v) &= \left(u-\frac 1 6 (u+i v)^{3}\right)w_1 + \left(v-\frac i 6 (u+i v)^{3}\right)w_2 + \frac 1 2 (u+i v)^{2} w_3,\\
e_3(u,v) &= -(u+i v) w_1 - (i u-v) w_2 + w_3,
\end{align*}
where $\{w_1,w_2,w_3\}$ is an orthonormal basis of $\C^3$. After applying a suitable similarity, we obtain case \eqref{case 4} of the theorem.
\end{proof}

As follows from Lemma \ref{ST from complex-analogous to pseudo-Riem} and Remark \ref{RM holomorphic extension}, any intersection of a parallel complex surface in $\C^3$ with a generic real slice results in a parallel surface in that real slice. On the other hand, from Lemma \ref{ST from pseudo-Riem to complex-analogous} it follows that all real-analytic parallel surfaces in any of the pseudo-Euclidean spaces $\R_k^n$ are obtained in this way. In this regard, it is interesting to compare our classification result above with the classification theorems of parallel surfaces in $\R^3$ and $\R_1^3$ (cf. \cite[Thm.~5.1]{ChVa}). It is worth noticing, that for the proof of Theorem \ref{ST classification hol riem}, only minor adaptions to the original proof for $\R_1^3$ were necessary.

\begin{theorem}\label{ST par in R3}
A non-degenerate parallel surface in $\R^3$ is similar to an open part of one of the following three surfaces:
\begin{enumerate}[(i)]
\item the plane $\R^2$, given by $x_3 = 0$,
\item the sphere $S^2$, given by $x_1^2+x_2^2+x_3^2=1$,
\item the flat cylinder $S^1\times \R$, given by $x_1^2+x_2^2=1$.
\end{enumerate}
\end{theorem}

These three cases in Theorem \ref{ST par in R3} above arise from intersecting the complex surfaces \eqref{case 1}, \eqref{case 2} and \eqref{case 3} from Theorem \ref{ST classification hol riem} with the Euclidean real slice $\R^3\subset\C^3$, given by $y_1=y_2=y_3=0$. The surface \eqref{case 4}, or any surface similar to it in $\C^3$, does not give a real 2-dimensional intersection with $\R^3$.

\begin{theorem}
A non-degenerate parallel surface in $\R_1^3$ (with metric $ds^2=-dy_1^2+dx_2^2+dx_3^2$) is similar to an open part of one of the following eight surfaces:
\begin{enumerate}[(i)]
\item the Euclidean plane $\R^2$, given by $y_1 = 0$,\label{LI Par R13 case 1}
\item the Lorentzian plane $\R_1^1$, given by $x_2 = 0$,\label{LI Par R13 case 2}
\item the hyperbolic plane $H^2$, given by $y_1^2-x_2^2-x_3^2=1$,\label{LI Par R13 case 3}
\item the indefinite sphere $S_1^2$, given by $-y_1^2+x_2^2+x_3^2=1$,\label{LI Par R13 case 4}
\item the flat cylinder $H^1\times \R^1$, given by $y_1^2-x_2^2=1$,\label{LI Par R13 case 5}
\item the flat cylinder $S^1\times \R_1^1$, given by $x_2^2+x_3^2=1$,\label{LI Par R13 case 6}
\item the flat cylinder $S_1^1\times \R^1$, given by $-y_1^2+x_2^2=1$,\label{LI Par R13 case 7}
\item the flat minimal Lorentzian surface given by $x_3 = (y_1-x_2)^2$.\label{LI Par R13 case 8}
\end{enumerate}
\end{theorem}

As these surfaces are all real-analytic (in fact even quadratic), they must all be obtainable from Theorem \ref{ST classification hol riem}, by taking appropriate intersections with $\R_1^3$. Cases \seqref{LI Par R13 case 1} and \seqref{LI Par R13 case 2} are obtained by intersecting the planes $z_1=0$ and $z_2=0$ (both similar to \eqref{case 1} from Theorem \ref{ST classification hol riem}) with $\R_1^3$. The cases \seqref{LI Par R13 case 3} and \seqref{LI Par R13 case 4} are obtained by intersecting with the surfaces $-z_1^2-z_2^2-z_3^2=1$ and $z_1^2+z_2^2+z_3^2=1$ (both similar to \eqref{case 2}). The cases \seqref{LI Par R13 case 5}, \seqref{LI Par R13 case 6} and \seqref{LI Par R13 case 7} are obtained by intersecting with the surfaces $-z_1^2-z_2^2=1$, $z_2^2+z_3^2=1$ and $z_1^2+z_2^2=1$ (all similar to \eqref{case 3}). Finally, \seqref{LI Par R13 case 8} is obtained by intersecting with $z_3 = (z_2+i z_1)^2$ (similar to \eqref{case 4}).

}

\cleardoublepage

%



%
  \graphicspath{{\chaptersdir/conclusion/image/}}%
  \chapter{Conclusion}\label{ch:conclusion}

In this chapter we will discuss some of the results we have obtained in this thesis, and point out in which aspects the theory and methods could be further improved. We will also mention some directions for further research.

\section*{Submanifold theory with Rinehart spaces}

In Chapter \ref{ch:chapter1} we have seen how the concept of Rinehart spaces enabled us to reformulate many classical results from the theory of Riemannian (sub)manifolds in a more general setting, such as the
fundamental theorem of Riemannian geometry, the projection of a Levi-Civita connection on a
submanifold, and the construction of non-flat space forms. As we already mentioned in the introduction, we have done so without the involvement of local or point-wise reasoning, for which the theory would have to be expanded further to include sheafs of Rinehart spaces, and theory on the spectrum of the Rinehart space. But even within the framework we have developed so far, there are many statements that could likely be improved on. An example of this is our analogue of the fundamental theorem of Riemannian geometry: the conditions that $2$ is a unit in $\mathbb{K}$ and that the cotangent space is undelimited, could likely be loosened. As has also been mentioned in \ref{p:extra_conditions}, we have focused our study to Rinehart spaces, which are the analogues of finite-dimensional Riemannian manifolds, but many of our results will likely extend, in some form or another, to more general (infinite-dimensional) Lie-Rinehart algebras.

There are also many topics which are well-studied in the context of classical Riemannian manifolds, but which we haven't even touched upon in our treatment in terms of Rinehart spaces. One can think for example of Lie groups, which among others can be regarded as isometry groups of the function space $\mathcal{O}$. Also higher degree differential forms and integration we haven't discussed. In fact, when we are given a derivation $\mathrm{d}_X$ for a tangent vector field $X$, then the concept of integration already pops up when we look for an operator $\mathrm{J}_X:\mathrm{d}_X(\mathcal{O})\to\mathcal{O}$ as right inverse to $\mathrm{d}_X$, which would be unique up to a constant. Another topic of interest is distributions, which naturally appear as the span (over the algebra of functions $\mathcal{O}$) of one or more elements from the tangent space $\mathfrak{X}$.

As the author, I kept being surprised during the writing of this thesis, how it comes that apparently for each geometric object of study, mathematics has provided an algebraic formalism which captures all of its essential properties, perhaps except only for the underlying physical intuition. So although there are still many geometric-analytic notions in Riemannian geometry and differential geometry in general, for which it may not be immediately clear how they could be conveniently described in algebraic terms, it is likely just a matter of effort to find out how. 

\section*{Submanifold theory and Wick relations}

In Chapter \ref{ch:chapter2} we have seen that for \mbox{(pseudo-)}\linebreak[0]Riemannian manifolds which are Wick-related, one can draw tight connections between the submanifolds of these spaces. Our method leans on a few different mathematical areas, among which complex analysis and complex bilinear forms, as well as holomorphic Riemannian geometry and the concept of Wick-rotations. Lemma \ref{ST from complex-analogous to pseudo-Riem} and Lemma \ref{ST from pseudo-Riem to complex-analogous} play an important role in our approach, as they enable us to infer geometric properties of a holomorphic Riemannian manifold from the geometric properties of one of its real slices, but also to infer such geometric properties for the other real slices of the holomorphic Riemannian manifold. The examples we provided were chosen with the purpose of exemplifying the methods we introduced. Of course, there are no real obstructions to apply the same principle to deal with other problems in submanifold theory as well. First of all, no restrictions are being put on the dimensions of the manifolds involved. Furthermore, when it comes to the geometric property under consideration, there is only the restriction that it is derived in an analytic manner from the main operators and tensors in submanifold theory. Hence, besides the properties we have seen in the examples above, other possible properties it could be applied to are constant mean curvature, constant sectional curvature, Einstein, quasi-minimality (marginally trapped), (semi- or \mbox{pseudo-)}\linebreak[0]symmetry, (semi- or \mbox{pseudo-)}\linebreak[0]parallelity, etc.

Despite the apparent flexibility of the method, the assumption that the submanifolds are embedded in an analytic way poses of course a limitation. Admittedly, submanifolds that are embedded smoothly but not analytically do rarely appear as the solution of an analytic differential equation, unless the solution involves a whole family of functions as a degree of freedom (like a generalized cylinder over a curve). And even in such a case, the family of analytic solutions typically forms a dense subset of the whole family of smooth solutions, so that the whole set of solutions may be easily inferred. Nonetheless, it is important to keep the required assumption of analyticity in mind, as this assumption is usually not required when studying submanifolds.

Finally, it should be noted that the main principle of extending a geometric property from a real manifold to its complexification, and from there translating it back to some of its other real slices, can likely be applied to other research areas than submanifold theory as well. The already established theory about real forms of complexified Lie algebras is an example that can be mentioned in this regard (cf. \cite{Hel}), but other applications one could think of is transferring distributions or differential equations from one real manifold to another Wick-related manifold. As we mentioned in the introduction, in fact similar techniques are (implicitly) being used in the physics literature whenever a Wick-rotation is applied, but among pure mathematicians such methods are still rarely studied, even in research areas where they could provide some benefit. Hopefully, this thesis can help to bring some more attention to it.

\cleardoublepage

%




\backmatter

\includebibliography



\printindex

\makebackcoverXII

\end{document}